\documentclass{amsart}

\usepackage{amsfonts,amssymb,verbatim,amsmath,amsthm,latexsym,textcomp,amscd}
\usepackage{latexsym,amsfonts,amssymb,epsfig,verbatim}
\usepackage{amsmath,amsthm,amssymb,latexsym,graphics,textcomp}
\usepackage{graphicx}
\usepackage{color}
\usepackage{url}
\usepackage{stmaryrd}
\usepackage{enumitem}

\usepackage{tikz}
\usepackage{mathtools}
\usetikzlibrary{positioning}
\usepackage{hyperref}
\hypersetup{
	colorlinks   = true, urlcolor     = blue, linkcolor    = black, citecolor    = black}
\begin{document}

\newtheorem{theorem}{Theorem}[section]
\newtheorem{prop}[theorem]{Proposition}
\newtheorem{lemma}[theorem]{Lemma}
\newtheorem{corollary}[theorem]{Corollary}
\newtheorem{definition}[theorem]{Definition}
\newtheorem{notation}[theorem]{Notation}
\newtheorem{convention}[theorem]{Convention}
\newtheorem{conjecture}[theorem]{Conjecture}
\newtheorem{claim}[theorem]{Claim}
\newtheorem{question}[theorem]{Question}
\newtheorem{remark}{Remark}
\newtheorem{example}{Example}
\newcommand{\map}{\rightarrow}
\newcommand{\C}{\mathcal C}
\newcommand\AAA{{\mathcal A}}
\def\AA{\mathcal A}

\def\L{{\mathcal L}}
\def\al{\alpha}
\def\A{{\mathcal A}}
\newcommand\BB{{\mathcal B}}
\newcommand\CC{{\mathcal C}}
\newcommand\DD{{\mathcal D}}
\newcommand\EE{{\mathcal E}}
\newcommand\FF{{\mathcal F}}
\newcommand\GG{{\mathcal G}}
\newcommand\GB{{\mathbb G}}
\newcommand\HH{{\mathcal H}}
\newcommand\II{{\mathcal I}}
\newcommand\JJ{{\mathcal J}}
\newcommand\KK{{\mathcal K}}
\newcommand\LL{{\mathbb L}}
\newcommand\LS{{\mathcal L}} 
\newcommand\MM{{\mathcal M}}
\newcommand\NN{{\mathbb N}}
\newcommand\OO{{\mathcal O}}
\newcommand\PP{{\mathcal P}}
\newcommand\QQQ{{\mathbb Q}}
\newcommand\QQ{{\mathcal Q}}
\newcommand\RR{{\mathbb R}}
\newcommand\SSS{{\Sigma}}
\newcommand\TT{{\mathcal T}}
\newcommand\UU{{\mathcal U}}
\newcommand\VV{{\mathcal V}}
\newcommand\WW{{\mathcal W}}
\newcommand\XX{{\mathcal X}}
\newcommand\YY{{\mathcal Y}}
\newcommand\ZZ{{\mathbb Z}}
\newcommand\hhat{\widehat}
\newcommand{\Lam}[1]{\ensuremath{\partial^{(2)}_{#1}}}
\newcommand\XH{\mathcal{X}}
\newcommand{\YH}{\mathcal{Y}}
\newcommand{\mh}{\hat{d}_1}
\newcommand{\X}{\widehat{X}}
\newcommand{\Y}{\widehat{Y}}
\newcommand{\F}{\widehat{F}}
\newcommand{\m}{\hat{d}}

\newcommand{\secref}[1]{Section~\ref{#1}}
\newcommand{\thmref}[1]{Theorem~\ref{#1}}
\newcommand{\lemref}[1]{Lemma~\ref{#1}}
\newcommand{\rmkref}[1]{Remark~\ref{#1}}
\newcommand{\propref}[1]{Proposition~\ref{#1}}
\newcommand{\corref}[1]{Corollary~\ref{#1}}
\newcommand{\probref}[1]{Problem~\ref{#1}}
\newcommand{\eqnref}[1]{~{\textrm(\ref{#1})}}
\newcommand{\boldG}{\mbox{{\bf G}}}

\def\Ga{\Gamma}
\def\Z{\mathbb Z}

\def\diam{\operatorname{diam}}
\def\dist{\operatorname{dist}}
\def\hull{\operatorname{Hull}}

\def\length{\operatorname{length}}
\newcommand\RED{\textcolor{red}}
\newcommand\GREEN{\textcolor{green}}
\newcommand\BLUE{\textcolor{blue}}
\def\mini{\scriptsize}

\def\acts{\curvearrowright}
\def\embed{\hookrightarrow}

\def\ga{\gamma}
\newcommand\la{\lambda}
\newcommand\eps{\epsilon}
\def\geo{\partial_{\infty}}
\def\bhb{\bigskip\hrule\bigskip}

\title[Relatively hyperbolic metric bundles and CT map]{Relatively hyperbolic metric bundles and Cannon-Thurston map}

\author{Swathi Krishna}
\address{UM-DAE CEBS, Mumbai, India.}
\date{\today}

\begin{abstract}
In \cite{krishna-sardar}, the authors define the notion of morphisms and pullbacks of metric (graph) bundles. Moreover, they show the existence of the Cannon-Thurston (CT) map for a pullback. In this paper, we prove a combination theorem for relatively hyperbolic metric (graph) bundles and also prove a relatively hyperbolic analogue of the result in \cite{krishna-sardar}. Further, we derive some applications of these results.
\end{abstract}

\maketitle

\section{Introduction}
\noindent
In a seminal work \cite{CT}, Cannon and Thurston proved that if $M$ is a closed hyperbolic $3$-manifold fibering over the circle with fiber $F$, and $\widetilde{F}$ and $\widetilde{M}$ denote the universal covers of $F$ and $M$ respectively, then the inclusion map $\tilde{i} : \widetilde{F} \to \widetilde{M}$ extends to a continuous map in the boundary. 
%
%
This was generalized to the context of geometric group theory by Mahan Mj. Let $X$ be a hyperbolic metric space and $G$ be a hyperbolic group acting on $X$ properly discontinuously and freely by isometries. Let $\Gamma$ denote the Cayley graph of $G$ and $i :G\to X$ be the natural orbit map. Let $\partial{X}$ (resp. $\partial{\Gamma}$) denote the Gromov boundary and $\widehat{X}$ (resp. $\widehat{\Gamma}$) denote the Gromov compactification of $X$ (resp. $\Gamma$). Then he asked the following:

\begin{question}\label{qn-mj}
	Does $i : \Gamma \to X$ extend continuously to $\widehat{i}: \widehat{\Gamma} \to \widehat{X}$?
\end{question} 	
\noindent
This was answered in negative by Baker and Riley in \cite{baker-riley} (also see \cite{matsuda-oguni}). Over the time, existence of Cannon-Thurston map has been proved for many special cases, in hyperbolic as well as relatively hyperbolic setting (see \cite{mahan-icm} for the history). Of particular interest to us is the existence of CT map for a hyperbolic group extension  proved by Mitra in \cite{mitra-ct}. This was generalized by the introduction of metric (graph) bundles by Mj and Sardar in \cite{pranab-mahan}. They proved a combination theorem for metric (graph) bundles with uniformly hyperbolic fibers, and also that the inclusion map of a fiber in the metric (graph) bundle admits CT.
This was further generalized in \cite{krishna-sardar} using the notion of pullback of metric (graph) bundles.\\
In this paper, we prove a combination theorem for relatively hyperbolic metric (graph) bundles and the existence of CT map for pullbacks in this case. This also generalize the existence of CT map for relatively hyperbolic group extensions due to Pal in \cite{Pal}.
\noindent
The following is the main result of the paper.
\begin{theorem}\label{main theorem}
	Let $(X,B,\pi)$ be a metric (graph) bundle satisfying the following.
	\begin{enumerate}[leftmargin=\parindent,align=left,labelwidth=\parindent,labelsep=1pt]
		\item $B$ is a hyperbolic metric space.
		\item Each fiber $F_b$, $b\in B$, is strongly hyperbolic relative to a collection of uniformly mutually cobounded subsets $\mathcal{H}_b=\{H_{b,\alpha}\}_{\alpha\in \Lambda}$.
		\item For each $b\in B$ the barycentre map $\partial\widehat{F}_b\to\widehat{F}_b$ is uniformly coarsely surjective.
		\item The fiber-identification maps are uniformly type preserving and satisfy qi-preserving electrocution condition.
		\item The induced coned-off metric (graph) bundle satisfies a $(\lambda_k, M_k,n_k)$-flaring condition, for each $k\geq 1$.
		\item Cone-bounded hallways strictly flare condition is satisfied.\end{enumerate}
	Then $X$ is hyperbolic relative to the family $\mathcal{C}$ of maximal cone-subbundles of horosphere-like spaces. Suppose $i:A\to B$ is a qi embedding and $Y$ is the pullback of $A$ under $i$. Then $Y$ is hyperbolic relative to the family $\mathcal{C}_Y$ of maximal cone-subbundles of horosphere-like spaces and $(Y,\mathcal{C}_Y)\to (X,\mathcal{C})$ admits the CT map.
\end{theorem}
\noindent
The following two results give a group theoretic version of the above result.  

\smallskip
\noindent
{\bf Theorem \ref{combi for short exact seqn}.} {\em Let $1\to (K,K_1)\to (G,N_G(K_1))\to (Q,Q_1)\to 1$ be a short exact sequence of pairs of finitely generated groups such that $K$ is strongly hyperbolic relative to a proper non-trivial subgroup $K_1$ and suppose $G$ preserves cusps. Let $G$ be weakly hyperbolic relative to $K_1$ and a cone-bounded strictly flare condition be satisfied by $\EE(G,K_1)$. Then $G$ is strongly relatively hyperbolic relative to $N_G(K_1)$.	 }

\smallskip
\noindent
 {\bf Theorem \ref{main application}.} {\em Suppose $1\to (K,K_1)\to (G,N_G(K_1))\to (Q,Q_1)\to 1$ is a short exact sequence of relatively hyperbolic groups
 	where $K$ is strongly hyperbolic relative to a proper non-trivial subgroup $K_1$ and $G$ preserves cusps. Suppose $G$ is weakly hyperbolic relative to $K_1$ and it is strongly hyperbolic relative to $N_G(K_1)$. Suppose $Q'$ is a qi embedded subgroup of $Q$ and $G_1=\pi^{-1}(Q')$.
 	Then $(G_1,N_{G_1}(K_1))$ is strongly relatively hyperbolic and $(G_1,N_{G_1}(K_1))\to(G,N_{G}(K_1))$ admits CT.}

\smallskip
\noindent
As in \cite{krishna-sardar}, we have a few applications. 

\smallskip
\noindent
{\bf Theorem \ref{CT leaf in Y}.} 
	{\em Suppose we have a metric graph bundle satisfying the hypotheses of \thmref{main theorem} such that $X$ is a proper metric space. Let $\gamma$ be a geodesic line in $\GG(Y,\CC_Y)$ with $\partial i_{Y,X}(\gamma(\infty))=\partial i_{Y,X}(\gamma(-\infty))$ in $\partial X^h$. Let $F=F_b$ be a fiber of $Y$. Then\\
	$(1)$ $\gamma(\pm\infty)\in\partial i_{F,Y}(\partial F^h)$.\\	
	$(2)$ There exists a uniform quasigeodesic line $\sigma$ in $F^h$ at a finite Hausdorff distance from $\gamma$ satisfying $\partial i_{F,Y}(\sigma(\pm\infty))=\gamma(\pm\infty)$.}
	
	\smallskip
	\noindent	
The converse of this result is trivially true. 
Further, we generalize the result characterizing relative quasiconvexity using CT laminations (cf. \cite[Theorem 1.3]{mj-rafi}). The hyperbolic analogue characterizing quasiconvexity using CT laminations was proved by Mj and Sardar in \cite{sardar-mj} for (closed) surface groups and free groups.

\smallskip
\noindent
{\bf Theorem \ref{ct lam rqc}.} {\em Let $K=\pi_1(S^h)$ be the fundamental group of a surface with finitely many punctures and let $K_1,\ldots,K_n$ be its peripheral subgroups. Let $1\to K\to G\stackrel{\pi}{\to} Q\to1$ and $1\to K_i\to N_G(K_i)\stackrel{\pi}{\to}Q\to 1$ be the induced short exact sequences of groups. Let $Q$ is a convex cocompact subgroup of the pure mapping class group of $S^h$. Let $Q_1<Q$ be a quasiconvex subgroup of $Q$ and let $\pi^{-1}(Q_1)=G_1$. Suppose $H$ is an infinite index relatively quasiconvex subgroup of $G_1$. Then $H$ is a relatively quasiconvex subgroup of $G$. 	}

\smallskip
\begin{remark}
	In \cite{krishna-sardar}, the results are done for length metric bundles, as geodesic metric bundles are slightly restrictive. We stick to geodesic metric bundles as most of the results in this paper for them also hold for length metric bundles.
\end{remark}

\noindent 
{\bf Outline of the paper:} 
In Section \ref{section 2} we recall basic definitions and results regarding hyperbolicity, relative hyperbolicity and Cannon-Thurston maps. In \secref{section 3}, we recall the basics of metric (graph) bundles, morphisms, pullbacks etc. We also define the notion of relatively hyperbolic metric (graph) bundles. In \secref{section 4} we prove a combination theorem for relatively hyperbolic metric (graph) bundles. We prove the existence of CT map in \secref{section 5} and in \secref{section 6} we look at some applications of the main result.

\section{Preliminaries}\label{section 2}
\noindent
We refer the reader to \cite{pranab-mahan} and \cite{krishna-sardar} for the definitions of some of the terms and some results appearing in the paper. We only recall the ones we need.
\begin{definition}\cite[Definition 1.1.1]{pranab-mahan}
	Let $(X,d_X),(Y,d_Y)$ be metric spaces. Let $\phi:X\to Y$ be a map and $\epsilon\geq 0$, $K\geq 1$.
Given an increasing function $f:[0,\infty)\to[0,\infty)$ with $\lim_{t\to \infty} f(t)=\infty$, the map $\phi$ is {\bf properly embedded as measured by $f$} if for any $x,y\in X$, $d_Y(\phi(x),\phi(y))\leq M$ for some $M\geq 0$ implies $d_X(x,y)\leq f(M)$.
\end{definition}
\begin{definition}
Let $I\subset \RR$ be a closed interval with end points in $\ZZ\cup\{\pm \infty\}$ and $J=\ZZ\cap I$ with the restricted metric from $\RR$. Then a $(K,\epsilon)$-qi embedding $\alpha: J\to X$ is called a {\bf dotted $(K,\epsilon)$-quasigeodesic}.\\
If $I$ is a finite interval, say $I=[m,n]$, then the length of $\alpha$ is given by $$\sum_{i=m}^{n-1} d_X(\alpha(i), \alpha(i+1)).$$	
\end{definition}

\begin{lemma}\label{elem-lemma1}\cite[Lemma 1.1]{pranab-mahan}
	\begin{enumerate}[leftmargin=\parindent,align=left,labelwidth=\parindent,labelsep=4pt]
		\item For every $K_1,K_2\geq 1$ and $D\geq 0$ there is $K_{\ref{elem-lemma1}}=K_{\ref{elem-lemma1}}(K_1,K_2,D)$,
		such that the following holds: 
		A $K_1$-coarsely Lipschitz map with a $K_2$-coarsely Lipschitz, $D$-coarse inverse is a 	$K_{\ref{elem-lemma1}}$-quasiisometry.
		\item  Given $K\geq 1$ and $R\geq 0$ there exists
		$D_{\ref{elem-lemma1}}=D_{\ref{elem-lemma1}}(K,R)$ such that the following holds:	
		Suppose $X,Y$ are any two metric spaces and $f:X\map Y$ is a $K$-quasiisometry which is $R$-coarsely surjective.
		Then there is a $K_{\ref{elem-lemma1}}$-quasiisometric $D_{\ref{elem-lemma1}}$-coarse inverse of $f$.
	\end{enumerate}
\end{lemma}
\begin{convention}
	$(1)$ The quasiisometric coarse inverse of a quasiisometry $\phi$ will be denoted by $\phi^{-1}$.\\
	$(2)$ For $K\geq 1$, a $(K,K)$-quasiisometric embedding will be referred to as a $K$-qi embedding and a $(K,K)$-(quasi)geodesic as a $K$-(quasi)geodesic.\\
	$(3)$ For a geodesic space $X$, a geodesic joining $x,y\in X$ will be denoted by $[x,y]_X$ or $[x,y]$ if the space in question is clear.\\
	$(4)$ For a path $\gamma$, $\gamma_-$ and $\gamma_+$ denote the initial point and terminal point respectively.
\end{convention}
\begin{lemma}\label{proving Lipschitz}\cite[Proposition 2.10]{pranab-mahan}
	Suppose $X$ is a geodesic metric space and $Y$ is any metric space. A map $f:X\to Y$ is coarsely $C$-Lipschitz for some $C\geq 0$ if for all $x_1, x_2\in X$, $d_X(x_1,x_2)\leq 1$ implies
	$d_Y(f(x_1),f(x_2))\leq C$.
\end{lemma}
\begin{lemma}\label{extend quasigeodesic}
	Given $K\geq 1,D\geq 0$, we have $K_{\ref{extend quasigeodesic}}=K_{\ref{extend quasigeodesic}}(K,D)$ such that following holds.
	Let $(X,d)$ be a geodesic space and let $\alpha$ be a $K$-quasigeodesic in $X$ joining $x_1,x_2\in X$. Let $y_1,y_2\in X$ such that $d(x_i,y_i)\leq D$, for $i=1,2$. Let $\alpha'$ be the concatenation of $\alpha$ and the geodesics $[x_1,y_1]$ and $[x_2,y_2]$. Then $\alpha'$ is a $K_{\ref{extend quasigeodesic}}$-quasigeodesic in $X$. 
\end{lemma}
\begin{proof}
	Let $p,q\in\alpha'$. If $p,q\in\alpha$, we know that $l(\alpha'|_{[p,q]})\leq Kd(p,q)+K$. So, let $p\in [x_1,y_1]$ and $q\in\alpha$. Then,\\ $l(\alpha'|_{[p,q]})=d_X(p,x_1)+l(\alpha|_{[x_1,q]})\leq D+Kd(x_1,q)+K\leq Kd(p,q)+KD+K+D.$ Finally, it is enough to check the case when $p\in[x_1,y_1]$ and $q\in[x_2,y_2]$. 
	Then $l(\alpha'|_{[p,q]})=d(p,x_1)+l(\alpha)+d_X(x_2,q)\leq Kd(x_1,x_2)+K+2D\leq Kd(p,q)+2KD+K+2D$.  
	We take $K_{\ref{extend quasigeodesic}}=K+2KD+2D$.
\end{proof}
\noindent
We recall from \cite{krishna-sardar} that given a length space $X$, we have a metric graph $Y$, a {\em (canonical) metric graph approximation of $X$}, with $V(Y)=X$ and there exists an edge between $x,y\in V(Y)$ if $d_X(x,y)\leq 1$. The identity map denoted by $\psi_X:X\to V(Y)$ is a quasiisometry and its inverse $\phi_X: Y\to X$ is defined as follows. On $V(Y)$, it is simply the inverse of $\psi_X$ and it maps any point in the interior of an edge to one of the end points of the edge. 
\begin{lemma}\label{length space qi to graph}\cite[Lemma 2.8]{krishna-sardar}
	$(1)$ $Y$ is a (connected) metric graph.\\
	$(2)$ The maps $\psi_X$ and $\phi_X|_{V(Y)}$ are coarsely $1$-surjective, $(1,1)$-quasiisometries.\\
	$(3)$ The map $\phi_X$ is a $(1,3)$-quasiisometry and it is a $1$-coarse inverse of $\psi_X$.
\end{lemma}
\begin{definition} Let $X$ be a geodesic metric space.
		Let $Y,Z\subset X$, $Y_1,Y_2\subset Y$. Then $Z$ {\bf coarsely bisects} $Y$ into $Y_1,Y_2$ in $X$ if 
		$Y=Y_1\cup Y_2$ and $Z$ coarsely disconnects $Y_1, Y_2$ in $X$.	This means that $Y_i\setminus Z\neq \emptyset$, $i=1,2$ and there exists $R\geq 0$ such that the following holds: For $y_i\in Y_i$, $i=1,2$, and any geodesic $\gamma$ in $X$ joining 
		$y_1, y_2$ we have $\gamma \cap N_R(Z)\neq \emptyset$.		
\end{definition}
\begin{definition}Let $X$ be a geodesic metric space, $A\subset X$.\\ 
	$(1)$ For $x\in X$, a point $y\in A$ is a {\bf nearest point projection} of $x$ on $A$ if for all $a\in A$, $d(x,y)\leq d(x,a)$.\\
	$(2)$ For $\epsilon>0$, $y\in A$ is an {\bf $\epsilon$-approximate nearest point projection} of $x$ on $A$ if for all $a\in A$, $d(x,y)\leq d(x,a)+\epsilon$.\\
	$(3)$ A nearest point projection map (resp. an $\epsilon$-{approximate nearest point projection map}) $P_A: X\to A$ is a map such that $P_A(a)=a$ for all $a\in A$ and 
	$P_A(x)$ is a nearest point projection (resp. an $\epsilon$-approximate nearest point projection) of $x$ on $A$ for all $x\in X\setminus A$.
\end{definition}


\subsection{Hyperbolicity}
In this subsection, we quickly recall some basic definitions and results of hyperbolic spaces. One may refer to \cite{gromov-hypgps}, \cite{GhH}, \cite{Shortetal} for more details. The following definition is due to Rips.
\begin{definition}
Let $X$ be a geodesic metric space and $\delta\geq 0$. A geodesic triangle $\Delta x_1x_2x_3$ in $X$ is $\delta$-{\em slim} if any side of the triangle is contained in the union of the $\delta$-neighbourhoods of the other two sides.\\	
Further, $X$ is said to be a $\delta$-{\bf hyperbolic} metric space if every geodesic triangle in $X$ are $\delta$-slim.
\end{definition}
\noindent
Gromov gave a more general definition for hyperbolicity by not requiring $X$ to be a geodesic space.
\begin{definition}
	Suppose $X$ is a metric space and $\delta \geq 0$.	
	Let $p\in X$. For a pair of points $x,y\in X$, their Gromov product with respect to $p$ is $$(x\cdot y)_p=\frac{1}{2}(d(x,p)+d(y,p)-d(x,y)).$$	
	Now, $X$ is $\delta$-{\em hyperbolic} if for all $x,y,z,p\in X$, $(x.y)_p\geq \,\mbox{min}\{(x.z)_p, (y.z)_p\}-\delta$.
\end{definition}
\begin{remark}
The two definitions are equivalent for a geodesic metric space \textup{(}See \cite[Section 6.3C]{gromov-hypgps}, \cite[Proposition 1.22, Chapter III.H]{bridson-haefliger}\textup{)}.
\end{remark}
\begin{lemma}\label{stab-qg}  
	{\em (Stability of quasigeodesics, \cite{GhH})} 
	For all $\delta\geq 0$ and $k\geq 1$, $\epsilon\geq 0$ there
	is a constant $D_{\ref{stab-qg}}=D_{\ref{stab-qg}}(\delta, k, \epsilon)$ such that the following holds:\\	
	Suppose $Y$ is a $\delta$-hyperbolic geodesic metric space. Then the Hausdorff distance between a geodesic and a $(k,\epsilon)$-quasigeodesic 
	joining the same pair of end points is less than or equal to $D_{\ref{stab-qg}}$.
\end{lemma}
\subsubsection{Quasiconvexity}
\begin{definition}
	Let $X$ be a geodesic metric space and $A\subseteq X$. For $K\geq 0$, $A$ is {\bf $K$-quasiconvex} in $X$ if any geodesic $\gamma$ with end points in $A$ is contained in a $K$-neighbourhood of $A$.
\end{definition}
\noindent
We list the results related to quasiconvexity that we need. For more results, one may refer to  \cite{pranab-mahan,krishna-sardar}.
\begin{lemma}\label{trivial lemma}\cite[Lemma 2.34]{krishna-sardar}
	Given $K\geq 0, R\geq 0, \delta\geq 0$ there is a constant $D=D_{\ref{trivial lemma}}(R,K,\delta)$ such that the
	following holds:\\
	Suppose $X$ is a $\delta$-hyperbolic metric space and $A$ is a $K$-quasiconvex subset of $X$.
	Suppose $x,y\in X$ and $\bar{x}, \bar{y}\in A$ respectively are their $1$-approximate nearest point projections on 
	$A$. Suppose $z\in [x,y]$ and $\bar{z}$ is a $1$-approximate nearest point projection of $z$ on $A$ satisfying $d(z,\bar{z})\leq R$. Then $d(z, [\bar{x},\bar{y}])\leq D$.
\end{lemma}
\begin{lemma}\label{big npp small nbd}\cite[Lemma 3.1]{mitra-trees}
	Given $\delta\geq 0,C\geq 1$, there exist $D_{\ref{big npp small nbd}}, M_{\ref{big npp small nbd}}$ such that the following holds:
	Let $Y$ be $C$-quasiconvex subsets of a $\delta$-hyperbolic metric space $(X,d)$. For any pair of points $z,z'$, if $d(P_Y(z),P_Y(z'))>D_{\ref{big npp small nbd}}$, then $[P_Y(z),P_Y(z')]\subset N_{M_{\ref{big npp small nbd}}}([z,z'])$.
\end{lemma}

\begin{definition} Suppose $X$ is a geodesic metric space and $A, B$ are two quasi-convex subsets. Let  $R>0$. 
	Then $A, B$ are {\bf mutually $R$-cobounded}, or simply  $R$-cobounded, if the set of all $1$-approximate nearest point 	projections of the points of $A$ on $B$ has a diameter at most $R$ and vice versa.\\	
	When $R$ is understood or is not important we just say that $A, B$ are cobounded.
\end{definition}
\subsection{Relative hyperbolicity}
We briefly recall some important definitions and basic results of relative hyperbolicity in this subsection. Let $(X,d)$ be a geodesic metric space and $\HH=\{H_\alpha\}_{\alpha \in\Lambda}$ be a collection of uniformly mutually cobounded subsets of $X$.
\begin{definition}\label{conedoff}
	Let $(X,d)$ be a geodesic metric space and $\HH=\{H_\alpha\}_{\alpha \in\Lambda}$ be a collection of uniformly mutually cobounded subsets of $X$. For each $H_\alpha\in\HH$, we introduce a vertex $\nu(H_\alpha)$ and join every element of $H_\alpha$ to the vertex by an edge of length $\frac{1}{2}$. This new space is denoted by $\widehat{X}=\EE(X,\HH)$. This new vertex $\nu(H_\alpha)$ is called a {\bf cone point} and the $H_\alpha$'s are called {\bf horosphere-like sets}. The space $\widehat{X}$ is called a {\bf coned-off space of $X$ with respect to $\HH$}.
\end{definition}
\begin{definition} $(1)$ The metric on $\widehat{X}$ is called the electric metric.\\
$(2)$ Any (quasi)geodesic in $\widehat{X}$ is called an {\bf electric (quasi)geodesic}.\\
$(3)$ A path in $\widehat{X}$ is said to be {\bf without backtracking} if it does not return to any horosphere-like subset $H_{\alpha}$, after leaving it.
\end{definition}
\begin{definition}
	Given $K\geq1$, electric $K$-quasigeodesics in $(X,H)$ are said to satisfy {\bf bounded penetration property} if for any two electric $K$-quasigeodesics without backtracking $\beta,\gamma$, joining $x,y\in X$ there exists $D = D(K)$ such that\\
	{\em Similar Intersection Patterns 1}: if one of $\{\beta,\gamma\}$, say $\beta$, meets a horosphere-like set $H_\alpha$, then the length (measured in the intrinsic path-metric on $H_\alpha$) from the entry point to the exit point of $\beta$ is at most $D$.\\
	{\em Similar Intersection Patterns 2}: if both $\{\beta,\gamma\}$ meet some $H_\alpha$ then the length (measured in the intrinsic path-metric on $H_\alpha$) from the entry (resp. exit) point of $\beta$ to that of $\gamma$ is at most $D$.	
\end{definition}
\begin{definition}
	 A geodesic metric space $X$ is {\bf strongly hyperbolic relative to a collection of subsets $\HH$} if the coned-off space $\mathcal{E}(X, \HH)$ is a hyperbolic metric space and $(X,\HH)$ satisfies a bounded penetration property.
\end{definition}
\noindent
Another equivalent definition of relatively hyperbolic spaces  was given by Gromov \cite{gromov-hypgps}, using cusped space. In \cite{manning-groves}, Manning and Groves defined cusped space for metric graphs. Similar constructions can be found in \cite{bowditch-relhyp}, \cite{cannon-cooper} etc.
\begin{definition}
Let $(H,d)$ be a geodesic metric space. Then the {\bf hyperbolic cone} of $H$, $H^h:=H\times[0,\infty)$ with the path metric $d_h$ is defined as follows:\\		
(1) For $(x,t),(y,t)\in H\times\{t\}$, $d_{h,t}((x,t),(y,t))= e^{-t}d(x,y)$, where $d_{h,t}$ is the induced path metric on $H\times \{t\}$. Paths joining $(x,t)$ and $(y,t)$ that lie in $H\times [0,\infty)$ are called horizontal paths.\\
(2) For $t,s\in[0,\infty)$ and any $x\in H$, $d_{h}((x,t),(x,s))= |t-s|.$ Paths joining such elements are called vertical paths.\\
In general, for $x,y\in H^h$, $d_{h}(x,y)$ is the path metric induced by these vertical and horizontal paths.
\end{definition}
\begin{definition}
Let $X$ be a geodesic metric space and $\HH$ be a collection of uniformly mutually cobounded subsets. For each $H\in\HH$, we attach a hyperbolic cone $H^h$ to $H$ by identifying $(x,0)$ with $x$ for all $x\in H$. This space is denoted by $X^h = \mathcal{G}(X,\HH)$. $X$ is said to be {\bf hyperbolic relative to $\HH$ in the sense of Gromov} if $\mathcal{G}(X,\HH)$ is a hyperbolic metric space.
\end{definition}
\begin{convention}
Let $X$ be a geodesic space, $A,B\subset X$.
Let $Y\in \{X,\widehat{X},X^h\}$.
\begin{enumerate}[leftmargin=\parindent,align=left,labelwidth=\parindent,labelsep=4pt]
\item If $X$ is hyperbolic relative to $\HH$, the hyperbolicity constant of $\widehat{X}$ and $X^h$ is denoted by $\delta$.	
\item For $x,y\in Y$, the distance between $x,y$ in $Y$ is denoted by $d_Y(x,y)$.
\item The Hausdorff distance between $A,B$ in $Y$ is denoted by $Hd_Y(A,B)$.
\item For $C\geq0$, $N^Y_{C}(A)$ denotes the $C$-neighbourhood of $A$ in $Y$.
\item Let $\gamma$ be a path in $X$. If $\gamma$ penetrates a horosphere-like set $H_\alpha$, we replace portions of $\gamma$ inside $H_\alpha$ by edges joining the entry and exit points of $\gamma$ in $H_\alpha$ to $\nu(H_\alpha)$. We denote the new path by $\hat{\gamma}$. If $\hat{\gamma}$ is an electric (quasi)geodesic, we call $\gamma$ a {\em relative (quasi)geodesic} in $X$.
\item A (quasi)geodesic in $X^h$ is called a {\em hyperbolic (quasi)geodesic}.
\item Let $\hat{\alpha}$ be an electric quasigeodesic without backtracking in $\widehat{X}$. For each $H_\alpha$ penetrated by $\hat{\alpha}$, let $x,y$ be the entry and exit points of $\hat{\alpha}$, respectively. We join $x$ and $y$ by a geodesic in $H^h_{\alpha}$. This gives a path in $X^h$ and we call it an {\em electro-ambient path}. This is called the {\em hyperbolization of $\hat{\alpha}$} in \cite{ps-kap}.
\end{enumerate} 
\end{convention}
\begin{remark}
Every path in the rest of the paper is without backtracking, unless specified otherwise.
\end{remark}
\begin{lemma}\label{electro ambient is quasigeodesic}\cite{pal-singh},\cite[Theorem 9.28]{ps-kap}
	Let $\hat{\alpha}$ be an electric quasigeodesic in $\widehat{X}$. Then the corresponding electro-ambient path $\alpha$ is a $K_{\ref{electro ambient is quasigeodesic}}$-quasigeodesic in $X^h$.
\end{lemma}
\begin{lemma}\label{new quasigeodesic}
Let $K\geq 1$. Let $h\in H$ and $\sigma_1$ be a $K$-quasigeodesic (ray) in $X^h$ with $h$ as an endpoint. Let $\sigma_2$ be the vertical ray $\{(h,t)\mid t\geq 0\}$ lying in $H^h$. Then the concatenation $\sigma=\sigma_1\ast\sigma_2$ is a $(2K+1)$-quasigeodesic ray in $X^h$. 
\end{lemma}
\begin{proof}
	The proof is very simple. Let $p,q\in\sigma$. If $p,q\in\sigma_1$, we know that $l(\sigma|_{[p,q]})\leq Kd(p,q)+K$ and if $p,q\in\sigma_2$, then $l(\sigma|_{[p,q]})=d(p,q)$. So, let $p\in\sigma_1$ and $q\in\sigma_2$. Let $q=(h,t)$ for some $t>0$. Then,\\ $l(\sigma|_{[p,q]})=l(\sigma_1|_{[h,p]})+t\leq Kd(h,p)+K+t$. Note that $d(p,q)\geq t$. Therefore, $l(\sigma|_{[p,q]})=(2K+1)d(p,q)+K$. \end{proof}
	
\noindent
Let $(a,t),(b,s)\in H^h$, for $t,s\in[0,\infty)$. Suppose $a\neq b$. Choose the smallest $r\geq\max\{s,t\}$ such that for $(a,r),(b,r)\in H\times\{r\}$, $d_{h,r}((a,r),(b,r))= e^{-r}d_H(a,b)\leq1$, where $d_{h,r}$ is the induced path metric on $H\times\{r\}$. 
Then, $d_{H}(a,b)\leq e^r$ and $\ln{d_H(a,b)}\leq r$. Let $\lambda_1$ and $\lambda_2$ denote the vertical paths in $H^h$ joining $(a,t)$ to $(a,r)$ and $(b,s)$ to $(b,r)$ respectively. Let $\lambda_0$ denote the horizontal path in $H^h$ joining $(a,r)$ to $(b,r)$. We denote the path $\lambda_1\ast\lambda_0\ast\lambda_2$ by $\lambda((a,t),(b,s))$.
If $a=b$, then $\lambda((a,t),(b,s))$ is simply the vertical segment in $H^h$ joining $(a,t),(b,s)$.  
\begin{lemma}\label{quasigeodesics}
The path $\lambda(a,b)$ is a $K_{\ref{quasigeodesics}}$-quasigeodesic in $H^h$.\end{lemma}
%
\noindent
The proof of the above result appears in the proof of \cite[Proposition 1.68]{ps-kap}, which actually shows that $H\to H^h$ is a proper embedding. We skip the proof of the above result as it is slightly long but prove the following.
\begin{lemma}\label{proper embedding of cones}
	There exists $f_{\ref{proper embedding of cones}}:[0,\infty)\to[0,\infty)$ satisfying the following.
	The inclusion $H\to H^h$ is a proper embedding as measured by $f_{\ref{proper embedding of cones}}$.
\end{lemma}
\begin{proof}
	Let $a,b\in H$ such that $d_{H^h}((a,0),(b,0))=N$. Then $l_{H^h}(\lambda((a,0),(b,0)))=l_{H^h}(\lambda_1)+l_{H^h}(\lambda_0)+l_{H^h}(\lambda_2)=2r+1$, where $d_H(a,b)\leq e^r$.\\	
	Now, $l_{H^h}(\lambda((a,0),(b,0)))\leq K_{\ref{quasigeodesics}}d_{H^h}((a,0),(b,0))+K_{\ref{quasigeodesics}}\leq K_{\ref{quasigeodesics}}N+K_{\ref{quasigeodesics}}$. Therefore, $r\leq \frac{1}{2}(K_{\ref{quasigeodesics}}N+K_{\ref{quasigeodesics}}-1)$. Thus, $d_H(a,b)\leq e^{\frac{1}{2}(K_{\ref{quasigeodesics}}N+K_{\ref{quasigeodesics}}-1)}$. So, for $f_{\ref{proper embedding of cones}}(m)=e^{\frac{1}{2}(K_{\ref{quasigeodesics}}m+K_{\ref{quasigeodesics}}-1)}$, we have the result.	
\end{proof}
\begin{lemma}\label{proper embedding in cone space}
	There exists $f_{\ref{proper embedding in cone space}}:[0,\infty)\to[0,\infty)$ satisfying the following.
The inclusion $X\to X^h$ is a proper embedding as measured by $f_{\ref{proper embedding in cone space}}$.
\end{lemma}
\begin{proof}
Let $x,y\in X$ such that $d_{X^h}(x,y)=N$. Let $\alpha$ be a geodesic in $X^h$ joining $x,y$. Suppose $H_1^h,\ldots,H_k^h$ are the horoball-like subsets penetrated by $\alpha$. Clearly $k\leq N$. We can write $\alpha=\alpha_1\ast\beta_1\ast\alpha_2\ast\cdots\beta_{k}\ast\alpha_{k+1}$, where $\alpha_i$ are the maximal subpaths of $\alpha$ lying outside the horoball-likes subsets and for $1\leq i\leq k$, let $\beta_i$ be the maximal subpath of $\alpha$ lying in $H^h_i$. Then $$l_{X^h}(\alpha)=\sum_{i=1}^{k+1}l_{X^h}(\alpha_i)+\sum_{i=1}^{k}l_{X^h}(\beta_i)=\sum_{i=1}^{k+1}l_{X}(\alpha_i)+\sum_{i=1}^{k}l_{H^h}(\beta_i).$$ By \lemref{proper embedding of cones}, for $1\leq i\leq k$, $d_{H}((\beta_i)_-,(\beta_i)_+)\leq f_{\ref{proper embedding of cones}}(l_{H^h}(\beta_i))$.\\
Now, $d_X(x,y)\leq \sum_{i=1}^{k+1}d_{X}((\alpha_i)_-,(\alpha_i)_+)+\sum_{i=1}^{k}d_{X}((\beta_i)_-,(\beta_i)_+)$. Then, $$d_X(x,y)\leq\sum_{i=1}^{k+1}l_{X}(\alpha_i)+\sum_{i=1}^{k}d_{H}((\beta_i)_-,(\beta_i)_+)\leq N+kf_{\ref{proper embedding of cones}}(N)\leq N(1+f_{\ref{proper embedding of cones}}(N)).$$ So for $f_{\ref{proper embedding in cone space}}(m)=m(1+f_{\ref{proper embedding of cones}}(m))$, we have the proof.
\end{proof}
\begin{remark}\label{identification}
Suppose $X$ is strongly hyperbolic relative to a collection of subsets $\HH$. Then the space  $\EE(\GG(X,\HH),\HH^h)$ is quasiisometric to $\EE(X,\HH)$. This is because $\EE(X,\HH)$ is isometrically embedded in $\EE(\GG(X,\HH),\HH^h)$ and $\EE(\GG(X,\HH),\HH^h)$ lies in a 1-neighbourhood of the image of $\EE(X,\HH)$.
\end{remark}
\begin{lemma}\label{lemma from limint}
	Given $R>0$, there exists $K_{\ref{lemma from limint}}=K_{\ref{lemma from limint}}(R)>0$ such that the following holds.
	Let $X$ be a geodesic metric space strongly hyperbolic relative to a collection of uniformly mutually cobounded subsets $\HH = \{H_{\alpha}\}_{\alpha \in \Lambda}$. Let $\gamma$ be a geodesic in $X^h$ with endpoints lying outside the horoball-like sets. If $x\in X$ such that $x\in N_R^{X^h}(\gamma)$, then $x\in N^X_{K_{\ref{lemma from limint}}}(\gamma\cap X)$.
\end{lemma}
\begin{proof} 
Let $y\in\gamma$ such that $d_{X^h}(x,y)\leq R$. If $y\in\gamma\cap X$, then we are done. So suppose $y\in\gamma\cap H^h$, for some $H\in\HH$. Let $\gamma_1$ denote the geodesic segment $\gamma|_{[a,b]}$, where $a$ denotes the entry point of $\gamma$ into $H$ and $b$ denotes the exit point of $\gamma$ from $H$. Let $t>0$ such that $d_{h,t}(a,b)\leq 1$ from the construction of $\lambda(a,b)$. By stability of quasigeodesics, $\textup{Hd}_{X^h}(\gamma_1,\lambda)\leq D_{\ref{stab-qg}}$, where $D_{\ref{stab-qg}}= D_{\ref{stab-qg}}(\delta,K_{\ref{quasigeodesics}},K_{\ref{quasigeodesics}})$. Since $y\in\gamma_1$, there exists $z=(h,s)\in \lambda$, for $0\leq s\leq t$, such that $d_{X^h}(y,z)\leq K_{\ref{stab-qg}}$ and we have, $d_{X^h}(x,z)\leq R+K_{\ref{stab-qg}}$. Depending on whether $z\in\lambda_1,\lambda_2$ or $\lambda_0$, $z=(a,s),(b,s)$ or $(h,t)$ for some $h\in H$, respectively.
So, $s\leq R+K_{\ref{stab-qg}}$ and $d_{X^h}(z,\{(a,0),(b,0)\}) \leq s+1\leq R+K_{\ref{stab-qg}}+1$. Thus, $d_{X^h}(x,\{(a,0),(b,0)\})\leq 2(R+K_{\ref{stab-qg}})+1$. As $X$ is properly embedded in $X^h$, $d_{X}(x,y)\leq f_{\ref{proper embedding in cone space}}(2(R+K_{\ref{stab-qg}})+1)$.
%
%
%
%
\end{proof}
\begin{lemma}\label{coned-off hyp}\cite[Proposition 4.6]{farb-relhyp}
	Given $\delta,C,D>0$, there exists $\Delta$ such that if X is a $\delta$-hyperbolic metric space with a collection $\HH$ of $C$-quasiconvex sets, then $\X$ is $\Delta$-hyperbolic.
\end{lemma}
\begin{lemma}\label{electric-hyp track}\cite{farb-relhyp}
$(1)$ Electric geodesics in $\EE(X,\HH)$ and relative geodesics in $X$ joining the same pair of points in $X$ have similar intersection patterns with $H$ for all $H\in\HH$, i.e. they track each other off horosphere-like sets.\\
$(2)$ Electric geodesics in $\EE(X,\HH)$ \textup{(after identifying with $\EE(\GG(X,\HH),\HH^h)$)}, and hyperbolic geodesics in $\GG(X,\HH)$ joining the same pair of points in $\GG(X,\HH)$ have similar intersection patterns with $H^h$ for all $H^h\in\HH^h$, i.e. they track each other off horoball-like sets.
\end{lemma}
\begin{lemma}\label{coned-off strong hyp}\cite{brahma-ibdd,farb-relhyp}
	Suppose $X$ is a $\delta$-hyperbolic metric space with a collection $\HH$ of $C$-quasiconvex, $D$-mutually cobounded subsets. Then $(X,\HH)$ is strongly relatively hyperbolic.
\end{lemma}
\subsection{Electric Projection}

Suppose $X$ is hyperbolic relative to the collection $\HH$.
\begin{definition}\cite{mahan-pal}
	Let $i :X^h \to\EE(\GG(X,\HH),\HH^h)$ be the inclusion map. We identify $\EE(\GG(X,\HH),\HH^h)$ with $\widehat{X}$. Let $\hat{\alpha}$ be an electric geodesic in $\widehat{X}$ and $\alpha$ be the corresponding electro-ambient quasigeodesic. Let ${\pi}_{\alpha}$ be a nearest point projection from $X^h$ onto $\alpha$.
		The {\bf electric projection} is the map $\hat{\pi}_{\hat{\alpha}} : \widehat{X} \to \hat{\alpha}$ defined as follows.		
		If $x\in X$, $\hat{\pi}_{\hat{\alpha}}(x)=i({\pi}_{\alpha}(x))$. 	
		If $x$ is a cone point of a horosphere-like set $H_{\alpha}\in \HH$, choose some $z\in H_{\alpha}$ and define  $\hat{\pi}_{\hat{\alpha}}(x) = i({\pi}_{\alpha}(z))$.
\end{definition}
\noindent
The following lemma shows that the electric projection is coarsely well defined.
\begin{lemma}\label{electric proj coarsely unique}\cite[Lemma 1.16]{mahan-pal}
	Given $\delta>0$, there exists $P_{\ref{electric proj coarsely unique}}>0$ depending $\delta$, $D_{\ref{big npp small nbd}}$ and $M_{\ref{big npp small nbd}}$ such that the following holds.\\
	For $H\in\HH$, $x,y\in H$ and a geodesic ${\hat{\alpha}}$ in $\widehat{X}$, we have $d_{\widehat{X}}(i({\pi}_{\alpha}(x)),i({\pi}_{\alpha}(y)) \leq P_{\ref{electric proj coarsely unique}}$.
\end{lemma}
\noindent
Further, as in the case of the nearest point projection maps in hyperbolic metric spaces, we have the following two results.
\begin{lemma}\cite[Lemma 1.17]{mahan-pal}\label{elec proj bound}
Given $\delta>0$, there exists $P_{\ref{elec proj bound}}=P_{\ref{elec proj bound}}(\delta)$ such that the following holds.	
For $x,y\in X$ and an electric geodesic $\hat{\lambda}$ in $\widehat{X}$, we have $d_{\widehat{X}}(\hat{\pi}_{\hat{\alpha}}(x)),\hat{\pi}_{\hat{\alpha}}(y))\leq P_{\ref{elec proj bound}}d_{\widehat{X}}(x,y)+P_{\ref{elec proj bound}}$.
\end{lemma}

\begin{lemma}\cite[Lemma 1.18]{mahan-pal}\label{almost commute}
Given $\delta>0$, $K\geq 1$, there exists $P_{\ref{almost commute}}=P_{\ref{almost commute}}(\delta,K)$ such that the following holds.\\
Suppose $\phi:X\to X$ is a type-preserving $K$-quasiisometry. Let $a,b\in\widehat{X}$ and $\widehat{\mu}_1$ be a quasigeodesic in $\widehat{X}$ joining $a,b$. Suppose the induced map $\hat{\phi}:\widehat{X}\to\widehat{X}$ is also a quasiisometry. Let $\widehat{\mu}_2$ be a quasigeodesic in $\widehat{X}$ joining $\hat{\phi}(a)$ and $\hat{\phi}(b)$. If $x\in\widehat{X}$, then $d_{\widehat{X}}(\hat{\pi}_{\hat{\mu}_2}(\hat{\phi}(x)),\hat{\phi}(\hat{\pi}_{\hat{\mu}_1}(x)))\leq P_{\ref{almost commute}}$.
\end{lemma}
\subsection{Partial Electrocution}\label{pel}
In \cite{brahma-amalgeo}, Mj modified the coning-off of horosphere-like sets in the following way.
Let $X$ be a geodesic metric space and $\HH$ be a collection of subsets as before. Let $\PP=\{\PP_\alpha\}_{\alpha}$ be a collection of uniformly hyperbolic metric spaces such that for each $\alpha$, there exists a uniform retraction $g_\alpha:H_\alpha\to\PP_\alpha$. Let $\GG=\{g_\alpha\}_{\alpha}$. 
\begin{definition}\label{pel space}
A {\bf partially electrocuted space} with respect to the quadruple $(X,\HH,\PP,\GG)$ is a space obtained gluing the mapping cylinders of $g_\alpha$ to $H_\alpha$ in $X$ by identifying each $H_{\alpha}\times\{0\}$ to $H_{\alpha}$ naturally. 
\end{definition} 
 \begin{convention}
 Given $(X,\HH,\PP,\GG)$ as in Definition \ref{pel space}.\begin{enumerate}[leftmargin=\parindent,align=left,labelwidth=\parindent,labelsep=4pt]
 \item The partially electrocuted space is denoted by $X_{pel}=\mathcal{PE}(X,\HH,\PP,\GG)$.
 \item The partially electrocuted metric is denoted by $d_{pel}$.
 \item The (quasi)geodesics in $X_{pel}$ are called the partially electrocuted (quasi)geodesics.
 \end{enumerate}
 \end{convention}
\noindent
Suppose $X$ is strongly hyperbolic relative to $\HH$.
\begin{lemma}\label{brahma-amalgeo}\cite[Lemma 7.1]{brahma-amalgeo}
 The space $(X_{pel},d_{pel})$ is a hyperbolic metric space and the sets $\PP_\alpha$ are uniformly quasiconvex.
\end{lemma}
\begin{lemma}\label{hyp-pel-track}\cite[Lemma 2.11]{mahan-reeves}
Given $K,\epsilon >0$, there exists $C>0$ such that the following holds:\\
Let 
$\gamma_{pel}$ and $\gamma^h$ denote respectively a $(K,\epsilon)$-partially electrocuted quasigeodesic in
$(X_{pel},d_{pel})$ and a hyperbolic $(K,\epsilon)$-quasigeodesic in $(X^h,d_h)$ joining $a,b$. Then $\gamma_{pel}$ and $\gamma^h$ have similar intersection patterns, i.e., $\gamma^h\cap X$ lies in a (hyperbolic) $C$-neighborhood of (any representative of) $\gamma_{pel}$. Further, outside of a $C$-neighborhood of the horoballs that $\gamma^h$ meets, $\gamma^h$ and $\gamma_{pel}$ track each other.
\end{lemma}
\begin{remark}\label{coned off-pel qi}
	The space obtained by coning-off $\PP_\alpha$'s in $(X_{pel},d_{pel})$ is same as $\EE(X,\HH)$. The map $j:\EE(X_{pel},\PP)\to\EE(X,\HH)$ which is identity on $X$ and sends $\PP_\alpha$ and the cone point of $\PP_\alpha$ to the cone point of the corresponding $H_\alpha$ is a quasiisometry.  	
\end{remark}
%
\begin{definition}
Let $(X_i,\HH_i)$ be strongly relatively hyperbolic metric spaces, for $i=1,2$, where each $\HH_i$ is a collection of subsets of $X_i$.  Let $\phi:X_1\to X_2$ be a quasiisometry.\\
$(1)$ $\phi$ is {\bf type preserving} if there exists $R\geq 0$ such that for every $H_1\in\HH_1$, there exists $H_2\in\HH_2$ such that $\phi(H_1)$ lies in an $R$-neighbourhood of $H_2$ in $X_2$ (resp. $\phi^{-1}(H_2)$ lies in an $R$-neighbourhood of $H_1$ in $X_1$).\\
$(2)$ $\phi$ satisfies {\bf qi-preserving electrocution condition} if the induced map $\widehat{\phi}:\X_1\to\X_2$ is also a quasiisometry.
\end{definition}
\begin{lemma}\label{type preserve implies qi}\cite[Lemma 1.2.31]{Pal-thesis}
Let $(X_i,\HH_i)$ be strongly relatively hyperbolic metric spaces, for $i=1,2$. Let $\phi:X_1\to X_2$ be a type preserving quasiisometry. Then the map $\phi^h:X_1^h\to X_2^h$ induced by $\phi$ is a quasiisometry.
\end{lemma}	
\noindent
See also \cite{sisto-mackay}. An immediate consequence of the above result is the following.
\begin{lemma}\label{for mg approx}  
Let $(X,d_X)$ be a geodesic metric space strongly hyperbolic relative to a collection of subsets $\HH'$. Let $Y$ be a metric graph approximation of $X$. Then $Y$ is strongly hyperbolic relative to $\HH$, where $\HH=\{\psi_X(H')\mid H'\in\HH'\}$.
\end{lemma}

\subsection{Boundaries of relatively hyperbolic spaces and CT maps}
We quickly go over Gromov boundary and properties first. For an extensive survey, one may refer to \cite{kapovich-benakli}.
\begin{definition}
Let $X$ be a hyperbolic geodesic metric space. Let $p\in X$.
\begin{enumerate}[leftmargin=\parindent,align=left,labelwidth=\parindent,labelsep=4pt]
\item The {\em geodesic boundary} $\partial_g{X}$ (resp. {\em quasigeodesic boundary} $\partial_{q}{X}$) of $X$ is an equivalence class of geodesic (resp. quasigeodesic) rays in $X$. Two geodesic (resp. quasigeodesic) rays $\alpha$ and $\beta$ are equivalent if $Hd(\alpha,\beta)<\infty$.
\item A sequence $\{x_n\}$ in $X$ is said to converge to infinity if $\lim_{i,j\map \infty}(x_i.x_j)_p=\infty$. Any two such sequences $\{x_n\}$ and $\{y_n\}$ are equivalent if $\lim_{i,j\map \infty}(x_i.y_j)_p=\infty$.
The {\em sequential boundary} $\partial_s{X}$ is set of all equivalence classes of sequences $\{x_n\}$ converging to infinity. 
\end{enumerate}
\end{definition}
\noindent
Note that the sequential boundary is independent of the choice of $p\in X$. For the three boundaries $\partial_g X,\partial_q X$ and $\partial_s X$, there is a bijection between them if $X$ is a proper geodesic hyperbolic space.  
We denote the boundary generally by $\partial X$ and $\overline{X}=X\cup\partial X$.

\begin{definition}
	$(1)$ A  sequence of points $\{\xi_n\}$ in $\partial X$, converges to $\xi\in\partial X$ if the following holds:\\
	 Let $\xi_n=[\{x^n_k\}_k]$ and $\xi=[\{x_k\}]$. Then 
	$\lim_{n\map \infty}(\liminf_{i,j\map \infty}(x_i.x^n_j)_p)=\infty$.\\
	$(2)$ The {\bf limit set of a subset $Y$} of $X$, $\Lambda Y=\{\xi\in \partial X:\exists\,\{y_n\}\subset Y\mbox{with}\,y_n\to\xi\}$. 
\end{definition}  
\noindent
Now suppose $X$ is a proper geodesic space strongly hyperbolic relative to a collection $\HH$. We identify $\EE(X^h,\HH^h)$ and $\widehat{X}=\EE(X,\HH)$.
\begin{lemma}\cite[Proposition 2.10]{hamenstadt}\label{hamenstadt}
The boundary $\partial X^h$ is homeomorphic to $\partial\widehat{X}\cup\bigcup_{H\in\HH}\partial H^h$.
\end{lemma}
\begin{remark}\label{remark about bdry}
The proof of \cite[Proposition 2.10]{hamenstadt} shows that for a geodesic ray $\gamma$ in $\partial X^h$, if $\gamma$ is unbounded in $\widehat{X}$ then $\gamma(\infty)$ gives a unique point in $\widehat{X}$. If $\gamma$ is bounded in $\widehat{X}$, then $\gamma(\infty)\in\bigcup_{H\in\HH}\partial H^h$. Also see \cite[Corollary 6.4,Lemma 6.12]{manning} and \cite[Theorem 3.2]{dowdall-taylor}, which together give Lemma \ref{hamenstadt}.
\end{remark}
\noindent
Note that for $X_{pel}=\PP\EE(X,\HH,\PP,\GG)$ as in \ref{pel}, $\EE(X_{pel},\PP)$ and $\EE(X,\HH)$ are also quasiisometric. Then we have the following.
\begin{lemma}\label{all boundaries}
$\partial\widehat{X}$ and $\partial\widehat{X}_{pel}$ are homeomorphic.
\end{lemma}
\begin{lemma}\label{lem: bdry defn}\cite[Lemma 2.45]{krishna-sardar}
	Let $x\in X$ be any point.
	Suppose $\{x_n\}$ is any sequence in $X$ and $\beta_{m,n}$ is a $k$-quasigeodesic joining $x_m$ to $x_n$ for all $m,n\in\NN$. Suppose $\alpha_n$ is a $k$-quasigeodesic joining $x$ to $x_n$.
	Then\\	
	$(1)$ $\{x_n\}\in \mathcal S$ if and only if $\lim_{m,n\map\infty} d(x,\beta_{m,n})=\infty$ if and only if there is a constant $D$ such that for all $M>0$ there is $N>0$ with 
	$Hd(\alpha_m\cap B(x;M),\alpha_n\cap B(x;M))\leq D$ for all $m,n\geq N$.\\	
	$(2)$ Suppose, moreover, $\xi\in \partial_s X$ and $\gamma_n$ is a $k$-quasigeodesic in $X$ joining $x_n$ to $\xi$ for all $n\in \NN$ and $\alpha$ is a $k$-quasigeodesic joining $x$ to $\xi$.	
	Then $x_n\to\xi$ if and only if $d(x,\gamma_n)\to\infty$ if and only if there is constant $D>0$ such that for all
	$M>0$ there is $N>0$ with $Hd(\alpha\cap B(x;M),\alpha_n\cap B(x;M))\leq D$ for all $n\geq N$.
\end{lemma}
\noindent
More generally, we have the following.
\begin{lemma}\label{convergence explained}\cite[Lemma 2.50]{krishna-sardar}
	Let $x\in X$ be any point. Suppose $\{\xi_n\}$ is any sequence of points in $\partial_s X$. 
	Suppose $\beta_{m,n}$ is a $k$-quasigeodesic line joining $\xi_m$ to $\xi_n$ for all $m, n\in \NN$ 
	and $\alpha_n$ is a $k$-quasigeodesic ray joining $x$ to $\xi_n$
	for all $n\in \NN$. Then\\	
	$(1)$ $\lim_{m,n\to\infty} d(x,\beta_{m,n})= \infty$ if and only if there is a constant $D=D(k,\delta)$ such that for all $M>0$ 
	there is $N>0$ with $Hd(\alpha_m\cap B(x;M),\alpha_n\cap B(x;M))\leq D$ for all $m,n\geq N$ and in this case,
	$\{\xi_n\}$ converges to some point of $\partial_s X$.\\	
	$(2)$ Suppose, moreover, $\xi\in\partial_s X$, $\gamma_n$ is a $k$-quasigeodesic ray in $X$ joining $\xi_n$ to $\xi$
	for all $n$, and $\alpha$ is a $k$-quasigeodesic ray joining $x$ to $\xi$. Then $\xi_n\map \xi$ if and only if $d(x,\gamma_n)\to\infty$ if and only if there is constant $D'=D'(k,\delta)$ such that for all
	$M>0$ there is $N>0$ with $Hd(\alpha\cap B(x;M),\alpha_n\cap B(x;M))\leq D$ for all $n\geq N$. In this case $\lim_{m,n\to\infty} d(x,\beta_{m,n})= \infty$.
\end{lemma}

\begin{lemma}\label{CT-limset}\cite[Lemma 2.60]{krishna-sardar}
	Suppose $X,Y$ are hyperbolic metric spaces and $f:Y\to X$ is a metrically proper map that admits CT. Then we have $\Lambda(f(Y))=\partial f(\partial Y)$ if either $Y$ is a proper metric space or $f$ is a qi embedding.
\end{lemma}	
\begin{lemma}\label{ideal triangles are slim}{\bf (Ideal triangles are slim)}\cite[Proposition 2.43]{krishna-sardar}
	Suppose $X$ is a $\delta$-hyperbolic metric space.
	for $x,y,z\in\overline{X}$, a $k$-quasigeodesic triangle formed by joining each pair of points from $\{x,y,z\}$ by $k$-quasigeodesics
 is an $R=R_{\ref{ideal triangles are slim}}(\delta, k)$-slim triangle.\\	
	In particular, if $\gamma_1,\gamma_2$ are two $k$-quasigeodesic rays with $\gamma_1(0)=\gamma_2(0)$
	and $\gamma_1(\infty)=\gamma_2(\infty)$, then $Hd(\gamma_1, \gamma_2)\leq R$. 
\end{lemma}
\noindent
Now we generalize \cite[Lemma 2.5]{krishna-sardar}. We skip the proof as it is almost verbatim.
\begin{lemma}\label{quasigeod criteria}
	Let $X$ be a geodesic metric space, $x,y\in\partial X$ (resp. $x,y\in X$), let $\gamma$ be a $k$-quasigeodesic in $X$ joining $x,y$ and 
	$\alpha:\RR\to X$ (resp. $\alpha:[0,l]\to X$) be a (dotted) coarsely $L$-Lipschitz path joining $x,y$. If $\alpha$ is a proper embedding as measured by a function 
	$f:[0,\infty)\to[0,\infty)$ and there exists $D>0$ such that $Hd(\alpha,\gamma)\leq D$, then $\alpha$ is (dotted) $K_{\ref{quasigeod criteria}}=K_{\ref{quasigeod criteria}}(k,f,D,L)$-quasigeodesic in $X$.
\end{lemma}
\begin{definition}
	Let $X$ be a geodesic space. For distinct points $x,y,z\in\overline{X}$ and a $k$-quasigeodesic triangle with vertices $x,y,z$, a point $x_0\in X$ such that $N_{r}(x_0)$ intersects all the three quasigeodesics, for some $r>0$, is called an {\bf $r$-barycenter} of the triangle. 
\end{definition}
\begin{lemma}{\em (Barycenters of ideal triangles, \cite[Lemma 2.46]{krishna-sardar})}\label{barycenter}
There exists $r_0>0$ such that given a hyperbolic geodesic space $X$, distinct points $x,y,z\in\overline{X}$, there exists an $r_0$-barycenter $x_0\in X$ of a $k$-quasigeodesic triangle with vertices $x,y,z$.
Moreover, there exists $L_0\geq 0$ such that any two barycenters are at most $L_0$ distance apart.\end{lemma}
\noindent
Thus, we have a coarsely well-defined map $\partial^3 X\to X$, i.e., the {\em barycenter map}. 
\begin{definition}\label{bowditch boundary}
Let $X$ be a geodesic metric space strongly hyperbolic relative to a collection of subsets $\HH$. Then the {\bf Bowditch boundary} of $(X,\HH)$ is the Gromov boundary of the hyperbolic space $X^h=\GG(X,\HH)$.
\end{definition}
\subsection{Cannon-Thurston maps}
\begin{definition}\cite{mitra-trees}\label{ct-defn}
A map $f:Y\to X$ of hyperbolic metric spaces is said to admit a {\bf Cannon-Thurston (CT) map}, if $f$ gives rise to a continuous map $\partial{f}:\partial{Y}\to\partial{X}$ in the following sense:\\	
	Given any $\xi\in\partial Y$ and a sequence of points $\{y_n\}$ in $Y$ converging to $\xi$, the sequence $\{f(y_n)\}$ 
	converges to a point of $\partial X$ independent of the choice of $\{y_n\}$ and the resulting map $\partial{f}:\partial{Y}\to \partial{X}$ is continuous. 
\end{definition}
\begin{lemma}{\em({\bf Mitra's criterion}, \cite[Lemma 2.1]{mitra-trees})}\label{CT existence}
	Suppose $X$, $Y$ are hyperbolic geodesic spaces and $f:Y\to X$ is a proper embedding. Then $f$ admits CT if the following holds:\\	
	(*) Let $y_0\in Y$. There exists a function $\tau:\RR_{\geq 0}\map \RR_{\geq 0}$, with the property that $\tau(n)\map \infty$ as 
	$n\map \infty$ such that for all geodesic segments $[y_1,y_2]_Y$ lying outside the $n$-ball around $y_0\in Y$, 
	any geodesic segment $[f(y_1),f(y_2)]_X$ lies outside the $\tau(n)$-ball around $f(y_0)\in X$.
\end{lemma}
\noindent
The definition of CT maps in the case relative hyperbolicity is the following.
\begin{definition}\cite[Definition 1.27]{mahan-pal} \label{ct-defn-1}
	Suppose, for $i=1,2$, $X_i$ is a geodesic space strongly hyperbolic relative to a collection of subsets $\HH_i$. A type preserving map $f:(X_1,\HH_1)\to(X_2,\HH_2)$ admits {\em Cannon-Thurston map}, if the induced map $f^h:X_1^h\to X^h_2$ admits CT.
\end{definition}
\noindent
Now we recall the version of Mitra's criterion for relatively hyperbolic spaces. The proof follows by Lemma \ref{CT existence} and Lemma \ref{hyp-pel-track}.
\begin{lemma}{\em({\bf Mitra's criterion for relatively hyperbolic spaces})}\label{CT existence-1}
	Let $X,Y$ be geodesic spaces strongly hyperbolic relative to a collection of subsets $\HH_X,\HH_Y$ respectively. Let $f:(X,\HH_X)\to(Y,\HH_Y)$ be a type preserving proper embedding. Then $f$ admits CT if the following holds:\\	
	(**) Let $x_0\in X$. There exists a function $\tau:\RR_{\geq 0}\map \RR_{\geq 0}$, with the property that $\tau(n)\to\infty$ as 
	$n\to\infty$ such that for $x,y\in X$ and any partially electrocuted geodesic $\alpha_{pel}$ joining $x,y$ in $X_{pel}$, if $\alpha_{pel}\cap X$ lies outside the $n$-ball around $x_0$ in $X$, then for any partially electrocuted geodesic $\beta_{pel}$ joining $f(x),f(y)$ in $Y_{pel}$, $\beta_{pel}\cap Y$ lies outside the $\tau(n)$-ball around $f(x_0)$ in $Y$.
\end{lemma}

\section{Metric bundles}\label{section 3}
\noindent
We start by recalling the definition of metric (graph) bundles.
\begin{definition}\label{defn-mbdle}\cite[Definition 1.2]{pranab-mahan}
	Suppose $(X,d)$ and $(B, d_B)$ are geodesic metric spaces; let $c\geq 1$ and let $\eta:[0,\infty) \rightarrow [0,\infty)$ be a function. Then $X$ is an $(\eta,c)$-{\bf metric bundle} over $B$ if there is a surjective $1$-Lipschitz
	map $\pi:X\rightarrow B$ such that the following conditions hold:\\
	{\bf (1)} For each $b\in B$, $F_b:=\pi^{-1}(b)$ is a geodesic metric space with respect to the path metric $d_b$ induced from $X$. The inclusion maps
	$i_b: (F_b,d_b) \rightarrow X$ are uniformly metrically proper as measured by $\eta$. \\
	{\bf (2)}  Suppose $b_1,b_2\in B$ with $d_B(b_1,b_2)\leq 1$ and let $\gamma$ be a geodesic in $B$ joining them. Then for any $x\in F_{b_1}$ there is a path $\tilde{\gamma}:[0,1]\to \pi^{-1}(\gamma)\subset X$
	of length at most $c$ such that $\tilde{\gamma}(0)=x$, $\tilde{\gamma}(1)\in F_{b_2}$.
\end{definition}
\begin{definition}\label{defn-mgbdl}\cite[Definition 1.5]{pranab-mahan}
	Suppose $X$ and $B$ are metric graphs. Let $\eta:[0,\infty) \rightarrow [0,\infty)$ be
	a function. We say that $X$ is an $\eta$-{\bf metric graph bundle}
	over $B$ if there exists a surjective simplicial map $\pi:X\rightarrow B$  such that:\\
	$1.$ For each $b\in V(B)$, $F_b:=\pi^{-1}(b)$ is a connected subgraph of $X$ and the inclusion maps
	$i_b: F_b\rightarrow X$ are uniformly metrically proper as measured by $\eta$ for the  path metric $d_b$ induced on $F_b$.\\
	$2.$ For adjacent vertices $b_1,b_2\in V(B)$, each $x\in V(F_{b_1})$ is connected by an edge to an element of $V(F_{b_2})$.
\end{definition}
\noindent
We recall the required results pertaining to metric (graph) bundles from \cite{krishna-sardar}. We follow the notations and conventions from \cite{krishna-sardar}. In both the above definitions, $F_b$ is called a fiber. A geodesic in a fiber is called a {\em fiber geodesic}. The spaces $X$ and $B$ are called the {\em total space} and the {\em base space} of the bundle respectively. 
\begin{definition}
	Suppose $\pi: X\to B$ is a metric (graph) bundle.\\	$(1)$ Let $A\subset B$ and $k\geq 1$. A {\bf $k$-qi section} over $A$ is a $k$-qi embedding
	$s:A\to X$ (resp. $s:V(A)\to X$) such that $\pi\circ s=I$, where $A$ has the restricted metric from $B$ and $I$ is the identity map of $A$ (resp. $V(A)$).\\	
	$(2)$ Given a metric space (resp. graph) $Z$ and any qi embedding $f: Z\to B$ (resp. $f:V(Z)\to V(B)$)
	a {\bf $k$-qi lift} of $f$ is a $k$-qi embedding $\tilde{f}:Z\to X$ (resp. $\tilde{f}:V(Z)\to V(X)$) such that $\pi\circ\tilde{f}=f$.
\end{definition}
\noindent
We mostly refer to the image of a qi section (resp. qi lift) as the qi section (resp. qi lift). 
\begin{convention}	
	Suppose $\gamma:I\to B$ is a (quasi)geodesic and $\tilde{\gamma}$ is a qi lift of $\gamma$. For $b=\gamma(t)$, $t\in I$, $\tilde{\gamma}(b)$ also denotes $\tilde{\gamma}(t)$. 
\end{convention}
\noindent
We recall some results from \cite{pranab-mahan} and \cite{krishna-sardar}. One may refer to the papers for proofs.
\begin{lemma}{\em(Path lifting lemma-I)}\cite[Lemma 3.8]{krishna-sardar}\label{lifting geodesics}
	Let $\pi:X\to B$ be an $(\eta,c)$-metric bundle or an $\eta$-metric graph bundle. Suppose $b_1,b_2\in B$. Suppose $\gamma:[0,L]\to B$ is a continuous, rectifiable, arc length parametrized path (resp. an edge path) in $B$ joining $b_1$ to $b_2$. For any $x\in F_{b_1}$, there is a path $\tilde{\gamma}$ in $\pi^{-1}(\gamma)$
		such that $l(\tilde{\gamma})\leq(L+1)c$ (resp $l(\tilde{\gamma})=L)$ joining $x$ to some point of $F_{b_2}$.\\		
		In particular,if $\gamma$ is a geodesic in $B$, then $\tilde{\gamma}$ is a geodesic. Further, if $\gamma$ is a $(k,\epsilon)$-(dotted) quasigeodesic,  for $k\geq 1$ and $\epsilon\geq 0$, then $\tilde{\gamma}$ is a $c.(k+\epsilon+1)$-qi lift. (Here, $c=1$ in the metric graph bundle case).
\end{lemma}
\begin{corollary}\cite[Corollary 3.9]{krishna-sardar}\label{path lifting remark}
	Given  any metric (graph) bundle $\pi:X\to B$ and $b_1, b_2\in B$ 
	we can define a map $\phi:F_{b_1}\to F_{b_2}$ such that 
	$d_X(x, \phi(x))\leq 3c +3cd_B(b_1, b_2)$ (resp. $d(x, \phi(x))=d_B(b_1,b_2)$) for all $x\in F_{b_1}$.
\end{corollary}

\noindent
Recall from \cite{krishna-sardar} that for $b_1,b_2\in B$, a map $f:F_{b_1}\to F_{b_2}$ satisfying $d_X(x, f(x))\leq D$ for some $D>0$ independent of $x$ is called a {\em fiber identification map}.
\begin{lemma}\cite[Lemma 3.10]{krishna-sardar}\label{fibers qi}
	Suppose $\pi:X\to B$ is an $(\eta,c)$-metric bundle or $\eta$-metric graph bundle. Let $R\geq 0$ and $b_1,b_2\in B$. Then we have the following.\\	
	$(1)$ $Hd(F_{b_1},F_{b_2})\leq 3c+ 3cd_B(b_1,b_2)$ (resp. $Hd(F_{b_1},F_{b_2})= d_B(b_1,b_2)$).\\		
	$(2)$ Suppose $\phi_{b_1b_2}:F_{b_1}\to F_{b_2}$ is a map such that $d(x,\phi_{b_1b_2}(x))\leq R$ for all $x\in F_{b_1}$. Then $\phi_{b_1b_2}$ is a $D_{\ref{fibers qi}}$-surjective, $K_{\ref{fibers qi}}=K_{\ref{fibers qi}}(R)$-quasiisometry.\\
	$(3)$ For a map $\psi_{b_1b_2}:F_{b_1}\to F_{b_2}$ with $d(x,\psi_{b_1b_2}(x))\leq R'$ for all $x\in F_{b_1}$, we have $d(\phi_{b_1b_2},\psi_{b_1b_2})\leq\eta(R+R')$.		
	In particular, the maps $\phi_{b_1b_2}$ are coarsely unique.
\end{lemma}
\begin{corollary}\cite[Corollary 3.11]{krishna-sardar}\label{fibers unif qi}
A fiber identification map $\phi_{b_1b_2}:F_{b_1}\to F_{b_2}$ satisfying $d_X(x,\phi_{b_1b_2}(x))\leq R$ for all $x\in F_{b_1}$ is a $K_{\ref{fibers unif qi}}=K_{\ref{fibers unif qi}}(R)$-quasiisometry.
\end{corollary}
\begin{corollary}\cite[Corollaries 1.14, 1.16]{pranab-mahan}\label{bdd-flaring} {\em (Bounded flaring condition)} 
	For all $k\in\RR$, $k\geq 1$ there exists a function $\mu_k:\NN\to\NN$ 
	such that the following holds:\\	
	Suppose $\pi: X\to B$ is an $(\eta,c)$-metric bundle or an $\eta$-metric graph bundle. Let $\gamma\subset B$ be a a geodesic joining $b_1,b_2\in B$, and let $\tilde{\gamma_1}$, $\tilde{\gamma_2}$ be two $k$-qi lifts of $\gamma$ in $X$. Suppose $\tilde{\gamma}_i(b_1)=x_i\in F_{b_1}$ and $\tilde{\gamma}_i(b_2)=y_i\in F_{b_2}$, $i=1,2$. Then f $d_B(b_1,b_2)\leq N$, we have $d_{b_2}(y_1,y_2)\leq \mu_k(N) \mbox{max}\{d_{b_1}(x_1,x_2),1\}$.
\end{corollary}
\begin{definition}\cite[Definition 3.13]{krishna-sardar}
	$(1)$ Suppose $(X_i, B_i,\pi_i)$, $i=1,2$ are metric bundles. A {\bf metric bundle morphism} from $(X_1,B_1,\pi_1)$ to  $(X_2,B_2,\pi_2)$ (or simply
	from $X_1$ to $X_2$) consists of a pair of coarsely $L$-Lipschitz maps $f:X_1\to X_2$ and $g:B_1\to B_2$ for some $L\geq 0$ such that $\pi_2\circ f=g\circ \pi_1$.
%
%
	
	\noindent
	$(2)$ Suppose $(X_i, B_i,\pi_i)$, $i=1,2$ are metric graph bundles. A {\bf metric graph bundle morphism} from $(X_1,B_1,\pi_1)$ to  $(X_2,B_2,\pi_2)$ (or simply
	from $X_1$ to $X_2$) consists of a pair of coarsely $L$-Lipschitz maps $f:V(X_1)\to V(X_2)$ and $g:V(B_1)\to V(B_2)$ for some $L\geq 0$ such that 
	$\pi_2\circ f=g\circ \pi_1$.\\	
	$(3)$ A morphism $(f,g)$ from a metric (graph) bundle $(X_1,B_1,\pi_1)$ to a metric (graph) bundle $(X_2, B_2, \pi_2)$ is an {\bf isomorphism} if there is a metric (graph) bundle morphism $(f',g')$ from $(X_2, B_2,\pi_2)$ to $(X_1,B_1,\pi_1)$ such that $f'$
	is a coarse inverse of $f$ and $g'$ is a coarse inverse of $g$.
\end{definition}
\noindent
For each $b\in B_1$, restriction of $f$ to $\pi^{-1}_1(b)$ is denoted by $f_{b}:\pi^{-1}_1(b)\to\pi^{-1}_2(g(b))$. 
\begin{definition}
	Suppose $(X_i, B,\pi_i)$, $i=1,2$ are metric (graph) bundles with the same base space $B$. Then $(X_1, B,\pi_1)$ is a {\bf subbundle} of $(X_2, B,\pi_2)$ if there is a metric (graph) bundle morphism $(f,g)$ from $(X_1, B,\pi_1)$ to $(X_2, B,\pi_2)$ such
	that all the fiber maps $f_b$, $b\in B$, are uniform qi embeddings and $g$ is the identity map on $B$ (resp. on $V(B)$).
\end{definition}
\begin{lemma}{\em (Restriction bundle)}\label{restriction bundle}
	Suppose $\pi:X\to B$ is a metric (graph) bundle and $A\subset B$ is a connected subset (resp. $A$ is a connected subgraph) such that
	any pair of points in $A$ can be joined by a path of finite length in $A$.
	Then the restriction of $\pi$ to $Y=\pi^{-1}(A)$ gives a metric (graph) bundle with the same
	parameters as that of $\pi: X\to B$ where $A$ and $Y$ are
	given the induced metrics from $B$ and $X$ respectively.\\
		Moreover, if $f:Y\to X$ and $g:A\to B$ are the inclusion maps then $(f,g):(Y,A)\to (X,B)$ is a metric
	(graph) bundle morphism.
\end{lemma}
\begin{definition}\label{pullback-defn} 
	$(1)$ {\bf (Pullback of a metric bundle)} Given a metric bundle $(X,B,\pi)$ and a coarsely Lipschitz map $g:B_1\to B$, a pullback of $(X,B,\pi)$ under $g$ is a metric bundle $(X_1,B_1,\pi_1)$ together with a morphism $(f:X_1\to X, g:B_1\to B)$ such that the following universal property holds: Suppose $\pi_2:Y\to B_1$ is 
	another metric bundle and $(f^Y, g)$ is a morphism from $Y$ to $X$. Then there is a coarsely unique morphism $(f', Id_{B_1})$ 
	from $Y$ to $X_1$ making the following diagram (Figure 1) commutative.

	\begin{figure}[ht]
		\centering
		\begin{tikzpicture}[node distance=1.5cm,auto]
			\node (P) {$X_1$};
			\node (B) [right of=P] {$X$};
			\node (A) [below of=P] {$B_1$};
			\node (C) [below of=B] {$B$};
			\node (P1) [node distance=1cm, left of=P, above of=P] {$Y$};
			\draw[->] (P) to node {$f$} (B);
			\draw[->] (P) to node [swap] {$\pi_1$} (A);
			\draw[->] (A) to node [swap] {$g$} (C);
			\draw[->] (B) to node {$\pi$} (C);
			\draw[->, bend right] (P1) to node [swap] {$\pi_2$} (A);
			\draw[->, bend left] (P1) to node {$f^Y$} (B);
			\draw[->, dashed] (P1) to node {$f'$} (P);
		\end{tikzpicture}
		\caption{}\label{pullback defn figure}
	\end{figure}
	\noindent
	$(2)$ {\bf (Pullback of a metric graph bundle)}
	In the case of a metric graph bundle, the diagram is replaced by one where we have the vertex sets instead of the whole spaces.
\end{definition}

\begin{lemma}\cite[Lemma 3.20]{krishna-sardar}\label{pullback lemma}
	Given $L\geq 0$ and $\phi_1,\phi_2:[0,\infty)\to[0,\infty)$ there exists $\phi:[0,\infty)\to[0,\infty)$ such that the following hold:\\	
	Suppose we have the following commutative diagram (Figure 2) of maps between metric spaces satisfying the properties (1)-(3) below.
	\begin{figure}[ht]
		\centering
		\begin{tikzpicture}[node distance=1.5cm,auto]
			\node (P) {$X_1$};
			\node (B) [right of=P] {$X$};
			\node (A) [below of=P] {$B_1$};
			\node (P1) [node distance=1cm, left of=P, above of=P] {$Y$};
			\draw[->] (P) to node {$f$} (B);
			\draw[->] (P) to node [swap] {$\pi_1$} (A);
			\draw[->, bend right] (P1) to node [swap] {$\pi_2$} (A);
			\draw[->, bend left] (P1) to node {$f^Y$} (B);
			\draw[->, dashed] (P1) to node {$f'$} (P);
		\end{tikzpicture}\caption{}\label{pic3.1}
	\end{figure}

    \noindent
	$(1)$ All the maps (except possibly $f'$) are coarsely $L$-Lipschitz.\\	
	$(2)$ If $d_{B_1}(b,b')\leq N$, then $Hd(\pi^{-1}_1(b), \pi^{-1}_1(b'))\leq\phi_1(N)$ for $b,b'\in B_1$ and $N\geq 0$.\\	
	$(3)$ The restriction of $f$ on the fibers of $\pi_1$ are uniformly properly embedded as measured by $\phi_2$.\\	
	Then $d_Y(y,y')\leq R$ implies $d_{X_1}(f'(y),f'(y'))\leq\phi(R)$ for all $y,y'\in Y$ and $R\in[0,\infty)$. In particular, if $Y$ is a geodesic space or the vertex set of a connected metric graph with restricted metric, then $f'$ is coarsely $\phi(1)$-Lipschitz.\\	
	Moreover, $f'$ is coarsely unique, i.e. there is a constant $D>0$ such that if another map $f'':Y\to X_1$ makes the above diagram commutative then $d(f',f'')\leq D$.
\end{lemma}
\begin{lemma}{\em (Pullbacks of metric bundles)}\cite[Proposition 3.19, Proposition 3.21]{krishna-sardar}\label{pullback bundle prop}
	Suppose $(X,B,\pi)$ is a metric bundle and $g:B_1\to B$ is a Lipschitz map. Then there is a pullback. Moreover, it is coarsely unique up to a metric bundle isomorphism. 
\end{lemma}
\begin{remark}
In the case of metric bundles, the pullback is the usual set-theoretic pullback.
\end{remark}
\begin{lemma}{\em (Pullbacks for metric graph bundles)}\cite[Proposition 3.23]{krishna-sardar}\label{pullback graph prop}
	Suppose $(X,B,\pi)$ is an $\eta$-metric graph bundle, $B_1$ is a metric graph and $g:V(B_1)\to V(B)$ 
	is a coarsely $L$-Lipschitz map for some constant $L\geq 1$. Then there exists a pullback $\pi_1:X_1\to B_1$ of $g$
	such that all the fiber maps $f_b:\pi^{-1}_1(b)\to \pi^{-1}(g(b))$, $b\in V(B_1)$ are isometries with respect to induced length metrics from $X_1$ and $X$ respectively and it is coarsely unique up to a metric graph bundle isomorphism.
\end{lemma}
\begin{corollary}\cite[Corollary 3.24]{krishna-sardar}\label{cor: restriction bundle}
	Suppose $\pi:X\to B$ is a metric graph bundle. 
	Suppose $A$ is a connected subgraph of $B$. Let $g:A\to B$ denote the inclusion map. Let $X_A=\pi^{-1}(A)$,
	$\pi_A$ be the restriction of $\pi$ and let $f:X_A\to X$ denote the inclusion map.
	Then $X_A$ is the pullback of $X$ under $g$.
\end{corollary}
\noindent
In \cite{pranab-mahan}, the authors proved that given a metric bundle, there exists is a natural metric graph bundle corresponding to it (see Lemma 1.17 through Lemma 1.21 in \cite{pranab-mahan}, \cite[Proposition 4.1]{krishna-sardar}), with the total space and base space of this metric graph bundle quasiisometric to the total space and the base space, respectively, of the original metric bundle. Further, by \cite[Proposition 4.2]{krishna-sardar}, there is a natural metric graph bundle morphism corresponding to a metric bundle morphism. We list these two results below.
\begin{prop}\label{bundle vs graph bundle}
	Suppose  $\pi':X'\to B'$ is an $(\eta,c)$-metric bundle. Then there is a metric graph bundle $\pi:X\to B$ along 
	with quasiisometries $\psi_B:B'\to B$ and $\psi_X: X'\to X$ such that\\ $(1)$ $\pi\circ \psi_X=\psi_B\circ \pi'$, and\\ 
	$(2)$ for $b\in B'$, $\psi_X|_{\pi'^{-1}(b)}:\pi'^{-1}(b)\to\pi^{-1}(\psi_B(b))$ is a uniform quasiisometry.\\	
	Moreover, the maps $\psi_X, \psi_B$ have coarse inverses $\phi_X$, $\phi_B$ respectively making the following diagram 
	commutative: 	
	\begin{figure}[ht]
		\centering
		\begin{tikzpicture}[node distance=2cm,auto]
			\node (A) {$X'$}; 
			\node (B) [right of=A] {$X$};
			\node (C) [below=1cm of A] {$B'$};
			\node (D) [below=1cm of B] {$B$};
			\draw [transform canvas={yshift=0.5ex},->] (A) to node [above] {$\psi_X$} (B);
			\draw [transform canvas={yshift=-0.5ex},->,dashed] (B) to node [below] [swap] {$\phi_X$} (A);
			\draw [->] (A) to node [swap] {$\pi'$} (C);
			\draw [->] (B) to node {$\pi$} (D);
			\draw [transform canvas={yshift=0.5ex},->] (C) to node [above] [swap]{$\psi_B$}(D);
			\draw [transform canvas={yshift=-0.5ex},->,dashed] (D) to node [below] {$\phi_X$} (C);
		\end{tikzpicture}
		\caption{}\label{pic4}
	\end{figure}
\end{prop}
\noindent
Let $\phi_X, \phi_B, \phi_Y$ and $\phi_A$ denote the coarse inverses of $\psi_X$, $\psi_B$, $\psi_Y$, and $\psi_A$ respectively. Suppose all these maps are $K$-quasiisometries, for $K\geq1$.
\begin{prop}\label{metric graph approx}
Let $(X',B',\pi'), (Y',A',p')$ be metric bundles and let $(f,g)$ be a metric bundle morphism from $Y'$ to $X'$. Let $(X,B,\pi), (Y,A,p)$ be the metric graph bundles corresponding to $(X',B',\pi'), (Y',A',p')$ respectively. Then\\
$(1)$ Let $\bar{f}:=\psi_X\circ f\circ\phi_Y$ and $\bar{g}:=\psi_B\circ g\circ\phi_A$. Then $(\bar{f},\bar{g})$ gives a morphism of metric graph bundles from $Y$ to $X$.\\
Moreover, if $Y'$ is the pullback of $X'$ under $g$ and $f$ is the pullback map, then $Y$ is the pullback of $X$ under $\bar{g}$ and $\bar{f}$ is the pullback map.\\
$(2)$ If $X',Y'$ are hyperbolic, then $f$ admits a CT map if and only if so does $\bar{f}$.
\end{prop}
\subsection{From metric graph bundles to coned-off metric graph bundles}
Let $\pi:X\to B$ be an $(\eta,c)$-metric bundle (resp. $\eta$-metric graph bundle) satisfying the following.\\
{\bf (0)} $B$ is $\delta_0$-hyperbolic.\\
{\bf (1)} Each fiber $F_b:=\pi^{-1}(b)$, for $b\in B$ (resp. $b\in V(B)$) is strongly hyperbolic relative to a collection of uniformly mutually cobounded subsets $\HH_b=\{H_{b,\alpha}\}_{\alpha\in\Lambda}$.\\
{\bf (2)} For every fiber, $\widehat{F}_b:=\EE(F_b,\HH_b)$ is $\delta_0$-hyperbolic with respect to the induced cone metric $\hat{d}_b$ (also denoted by $d_{\widehat{F}_b}$).\\
{\bf (3)}  There exists $D_0>2c$ such that the following holds. Suppose $b_1,b_2\in B$ (resp. $b_1,b_2\in V(B)$) and let $\gamma$ be a path of length at most $1$ in $B$ joining them (resp. let $e$ be an edge in $B$ joining them). Then for any $x\in F_{b_1}$ (resp. $x\in V(F_{b_1})$) there is a path $\tilde{\gamma}:[0,1]\to\pi^{-1}(\gamma)\subset X$ 
of length at most $c$ (resp. there is an edge $\tilde{e}$) joining $x$ to some $x'\in F_{b_2}$ (resp. $x'\in V(F_{b_2})$). If $x\in H_{b_1,\alpha}$ for some $H_{b_1,\alpha}\in\HH_{b_1}$, there exists $H_{b_2,\beta}\in\HH_{b_2}$ with $x'\in N_{D_0}(H_{b_2,\beta})$. \\
{\bf (4)} The fiber identification maps are uniformly type preserving and satisfy qi-preserving electrocution condition. Moreover for a metric bundle $(X',B',\pi')$ and the corresponding approximating metric graph bundle $(X,B,\pi)$, for each $b\in B'$, the map  $\psi_X|_{\pi'^{-1}(b)}$ (see Proposition \ref{bundle vs graph bundle}) also satisfies qi-preserving electrocution condition.\\
Such a metric (graph) bundle will be referred to as a {\bf relatively hyperbolic metric (graph) bundle}. 
\begin{corollary}{\em(Path lifting lemma-II)}\label{lifting geodesics-2}
	Suppose $\pi:X\to B$ is an $\eta$-relatively hyperbolic metric (graph) bundle. Let $b_1, b_2\in B$ (resp. $b_1,b_2\in V(B)$). For any $x\in F_{b_1}$, a geodesic $\gamma$ can be lifted to a geodesic $\tilde{\gamma}$ joining $x$ to some $y\in F_{b_2}$. In particular, if $x\in H_{b_1,\alpha}$, for some $H_{b_1,\alpha}\in\HH_{b_1}$, then there exists $D_{\ref{lifting geodesics-2}}$ depending on $d_B(b_1,b_2)$ such that for every $b\in\gamma$, $\tilde{\gamma}(b)\in N_{D_{\ref{lifting geodesics-2}}}(H_{b,\alpha})$ for some $H_{b,\alpha}\in\HH_b$. 
\end{corollary}
\begin{proof}
	We prove this only for the metric graph bundle case as the metric bundle case is similar. Let $b_1=u_0,u_1,\ldots,u_n=b_2$ be a sequence of consecutive vertices on $\gamma$ such that $d_B(u_i,u_{i+1})=1$, for $0\leq i\leq n-1$. Now  any $x_i\in F_{u_i}$ can be attached to some point in $F_{u_{i+1}}$ by an edge. Hence, given any $x=:x_0\in F_{b_1}$ we can inductively construct a sequence of points $x_i\in F_{u_i}$, $1\leq i\leq n$ and a sequence of edges $e_i$ joining $x_i$ to $x_{i+1}$ for $0\leq i\leq n-1$. We denote the concatenation of these edges by $\tilde{\gamma}$. Note that $l(\tilde{\gamma})= n$. Since $\pi$ is $1$-Lipschitz, it is a geodesic. Let $x_0$ lies in $H_1\in\HH_{b_1}$ in $F_{b_1}$. Then by the definition of relatively hyperbolic metric graph bundle, $x_1\in N_{D_0}(H_{u_1})$, for some $H_{u_1}\in\HH_{u_1}$ in $F_{u_1}$. Inductively, each $x_i$ lies in a $(iD_0+2(i-1))$-neighbourhood of a horosphere-like subset of $F_{u_i}$. Thus, $\tilde{\gamma}$ passes through an $(2(n-1)+nD_0)$-neighbourhood a union of horosphere-like subsets of $F_{u_i}$, $0\leq i\leq n$. So, we have $D_{\ref{lifting geodesics-2}}=(2(n-1)+nD_0)$.
\end{proof}
\noindent
Now we show that if $(X',B',\pi')$ is a relatively hyperbolic metric bundle, then the approximating metric graph bundle is a relatively hyperbolic metric graph bundle. We recall some parts of the construction of $X$ from \cite[Proposition 4.1]{krishna-sardar} without going into details. The space $B$ and each fiber $F_{\psi_B(b)}$ is the metric graph approximation of $B'$ and $F_b$, for $b\in B'$ respectively. Recall that $\phi_B$ is a qi inverse of $\psi_B$ and $\phi_X$ is a qi inverse of $\psi_X$ and all these maps are $K$-qi.
\begin{prop}\label{approx bundle rel hyp}
	Given $\eta':(0,\infty)\to(0,\infty)$, there exists $\eta:\NN\to\NN$ such that the following holds.
	Suppose $(X',B',\pi')$ is an $(\eta',c)$-relatively hyperbolic metric bundle. Then $(X,B,\pi)$ is an $\eta$-relatively hyperbolic metric graph bundle. 
\end{prop}
\begin{proof}
	Clearly, $B$ is a hyperbolic metric graph. The vertex set of $B$ is $B'$ and that of $X$ is $X'$. There are two types of edges in $X$. Let $x,y\in V(X)$ such that $\pi(x)=\pi(y)=b\in V(B)$. Then there exists an edge joining $x,y$ if $d_{\phi_B(b)}(x,y)\leq 1$. If $x\in\pi^{-1}(b)$ and $y\in\pi^{-1}(b')$, then there exists an edge joining $x,y$ in $X$ if and only if $d_B(b,b')=1$ and $d_{X'}(x,y)\leq c$. In this case if $x$ is an element in the horosphere-like subset of $\pi'^{-1}(b)$, then $y$ will be an element in a uniform neighbourhood of a horosphere-like subset of $\pi'^{-1}(b')$.\\	
	Now, if $\pi'^{-1}(b)$ is strongly hyperbolic relative to $\HH'_b$, then by \lemref{for mg approx}, $\pi^{-1}(b)$ is strongly hyperbolic relative to the collection of image of elements of $\HH'_b$ under $\psi_X|_{\pi'^{-1}(b)}$. Therefore, the coned-off space of $\pi^{-1}(b)$ is uniformly hyperbolic. Now we find $\eta$. Let $x,y\in\pi^{-1}(b)$ such that $d_X(x,y)\leq N$. Then $d_{X'}(x,y)\leq KN+K$. Therefore, $d_{\pi'^{-1}(b)}(x,y)\leq\eta'(KN+K)$. Then $d_{\pi^{-1}(b)}(x,y)\leq\eta'(KN+K)+1$, by \lemref{length space qi to graph}. Therefore, $\eta(N)=\eta'(KN+K)+1$.\\	
	Let $b_1,b_2\in V(B)$. Let $\phi'_{b_1b_2}:\pi'^{-1}(b_1)\to\pi'^{-1}(b_2)$ be a fiber identification map in $(X',B',\pi')$. Then $\psi_X|_{\pi'^{-1}(b_2)}\circ\phi'_{b_1b_2}\circ\phi_X|_{\pi^{-1}(b_1)}$ is a fiber identification map $\pi^{-1}(b_1)\to\pi^{-1}(b_2)$ denoted by $\phi_{b_1b_2}$. This is a composition of type preserving maps, and hence is type preserving. Similarly, if $\phi'_{b_1b_2}$ satisfies qi-preserving electrocution condition, then so does $\phi_{b_1b_2}$.
\end{proof}
\noindent
Again by \lemref{type preserve implies qi} and \lemref{for mg approx}, if $X'$ is strongly hyperbolic relative to a collection of subsets $\CC'$, then $X$ is strongly hyperbolic relative to the the collection of image of elements of $\CC'$ under $\psi_X$. So for the rest of this section, we work only with metric graph bundles.

\smallskip
\noindent
 Suppose $(X,B,\pi)$ is an $\eta$-relatively hyperbolic metric graph bundle. Let $b_1,b_2\in V(B)$ such that $d_B(b_1,b_2)=1$. Let $e$ be the edge joining $b_1,b_2$ in $B$. For each $H_{b_1,\alpha}\in\HH_{b_1}$, there exists $H_{b_2,\beta}\in\HH_{b_2}$ such that $V(H_{b_1,\alpha})$ is joined to $N_{D_0}(V(H_{b_2,\beta}))$ by edges. So, we attach $\nu(H_{b_1,\alpha})$ and $\nu(H_{b_2,\beta})$ by an edge. We denote the space obtained after attaching these new edges to $\EE(X,\HH)$ by $\XH$. Let $d_{\XH}$ denote the induced metric on it. We will show that $\XH$ is a metric graph bundle over $B$. 
\begin{prop}\label{new mbdl}
	The space $\XH$ is a $\eta_{\ref{new mbdl}}$-metric graph bundle over $B$.
\end{prop}
\begin{proof} For each $b\in V(B)$, $\widehat{F}_b:=\pi^{-1}(b)$ is a connected subgraph of $\XH$ since $F_b$ is a connected subgraph of $X$. The natural extension of $\pi$ to $\pi:\XH\to B$ is clearly $1$-Lipschitz. We need to check the following:\\
	(1)	 $\widehat{F}_b$ is uniformly properly embedded in $\XH$ with respect to $\hat{d}_b$.\\
	(2) If $b_1,b_2\in V(B)$ are adjacent vertices, then each vertex $x_1$ of $\widehat{F}_{b_1}$ is connected by an edge with a vertex in $\widehat{F}_{b_2}$.\\
	Condition (2) follows by the construction of $\XH$.
	Now we check (1). Let $x,y\in V(\widehat{F}_b)$ such that $d_{\XH}(x,y)=n$. Suppose $x,y\in V(F_b)$. Let $\alpha$ be a geodesic in $\XH$ joining $x,y$. Let $\alpha_1\ast\beta_1\ast\cdots\beta_k\ast\alpha_{k+1}$ be a decomposition of $\alpha$, where $\alpha_i$'s are the maximal portions lying outside the elements of $\HH$, and each $\beta_i$ is the maximal portion lying inside, say $H_{b_i}\in\HH_{b_i}$, with $(\alpha_1)_-=x,(\alpha_{k+1})_+=y$, $(\alpha_i)_+=(\beta_i)_-$ and $(\beta_i)_+=(\alpha_{i+1})_-$, for $1\leq i\leq k$. For $1\leq i\leq k+1$, let $p_i=(\alpha_i)_-,q_i=(\alpha_i)_+$. Then each $\beta_i$ joins $q_i,p_{i+1}$. Also note that $k\leq n$. For each $1\leq i\leq k+1$, $d_B(b,\pi(p_i))\leq n$ and $d_B(b,\pi(q_i))\leq n$. \\	
	For each $1\leq i\leq k+1$, let $\gamma_i$ be a geodesic in $B$ joining $b$ and $\pi(p_i)$ and $\gamma'_i$ be a geodesic in $B$ joining $b$ and $\pi(q_i)$. Let $\tilde{\gamma}_i$ be a lift of $\gamma_i$ starting at $p_i$ and let $\tilde{\gamma}'_i$ be be a lift of $\gamma_i$ starting at $q_i$. Then $l_X(\tilde{\gamma}_i)=l_X(\tilde{\gamma}'_i)\leq n$.
	 Suppose $\tilde{\gamma}_i(b)=x_i,\tilde{\gamma}'_i(b)=y_i$. Note that $d_{\widehat{F}}(x_i,y_i)\leq 1+D_{\ref{lifting geodesics-2}}$, where $D_{\ref{lifting geodesics-2}}$ depends on $n$.\\	 
	 Then $d_X(x,x_1)\leq l_X(\alpha_1)+l_X(\tilde{\gamma}_1)\leq 2n$. In general, for $1\leq i\leq k-1$, $d_X(y_i,x_{i+1})\leq l_X(\tilde{\gamma}'_i)+l_X(\alpha_{i+1})+l_X(\tilde{\gamma}_{i+1})\leq 3n$. Similarly, $d_X(y_k,y)\leq l_X(\tilde{\gamma}'_k)+l_X(\alpha_{k+1})\leq 2n$. All these imply that putting $x=y_0$ and $y=x_{k+1}$, we have $d_{\F_b}(y_i,x_{i+1})\leq\eta(3n)$, for $0\leq i\leq k$. Thus, $d_{\F_b}(x,y)\leq d_{\F_b}(x,x_1)+d_{\F_b}(x_1,y_1)+d_{\F_b}(y_1,x_2)+\cdots+d_{\F_b}(x_k,y_k)+d_{\F_b}(y_k,y)\leq (k+1)\eta(3n)+k(1+D_{\ref{lifting geodesics-2}})$.\\	 
	 So, we get $\eta_{\ref{new mbdl}}(n)=(n+1)\eta(3n)+n(1+D_{\ref{lifting geodesics-2}})$.
\end{proof}
\noindent
In the following lemma, we show that a metric graph bundle morphism between relatively hyperbolic metric graph bundles induces a metric graph bundle morphism between the corresponding coned-off metric graph bundle morphisms. 
\begin{lemma}\label{induced morphism}
	Let $\pi_i:X_i\to B_i$, $i=1,2$ be relatively hyperbolic metric graph bundles with $\HH_i=\bigcup_{b\in V(B_i)}\HH_b$. Let $(f,g)$ be an $R$-type preserving morphism from $(X_1,\HH_1)$ to  $(X_2,\HH_2)$, for some $R>0$. The induced map $(\hat{f}:\XH_1\to\XH_2,g)$ is a metric graph bundle morphism. 
\end{lemma}
\begin{proof}
	Note that the commutativity of $(\hat{f},g)$ easily follows from that of $(f,g)$. So it is enough to show that $\hat{f}$ is uniformly coarsely Lipschitz.\\ Let $x,y\in V(\XH_1)$ such that $d_{\XH_1}(x,y)= 1$. One possibility is that $d_{X_1}(x,y)= 1$. Then as $f$ is a metric bundle morphism, we have $d_{\XH_2}(\hat{f}(x),\hat{f}(y))\leq d_{X_2}(f(x),f(y))\leq D$. Otherwise, there exists $H\in\HH_1$ such that $x,y\in V(H)$ or one of the elements is in $H$ and the other is a cone point. Since $f$ is type preserving, $f(H)$ lies in an $R$-neighbourhood of $\HH_2$ in $X_2$ and $\hat{f}(\nu(H))$ is a cone point in $\XH_2$. So, in these two case, $d_{\XH_2}(\hat{f}(x),\hat{f}(y))\leq 2R+1$. Thus, by \lemref{proving Lipschitz}, $\hat{f}$ is uniformly coarsely Lipschitz.
\end{proof}
\begin{corollary}\label{induced isomorphism}
	Let $\pi_i:X_i\to B_i$, $i=1,2$ be relatively hyperbolic metric graph bundles. Let $(f,g)$ be a type preserving isomorphism from $(X_1,\HH_1)$ to  $(X_2,\HH_2)$. Then the induced map $(\hat{f}:\XH_1\to\XH_2,g)$ is a metric graph bundle isomorphism. 
\end{corollary}
\begin{proof}
By the definition of a metric graph bundle isomorphism, there exists a morphism $(f',g')$ from $X_2$ to $X_1$ and $D\geq 0$ such that for any $x_1\in V(X_1)$, $d_{X_1}(f'\circ f(x_1),x_1)\leq D$ and for any $x_2\in V(X_2)$, $d_{X_2}(f\circ f'(x_2),x_2)\leq D$.\\
By \lemref{induced morphism}, $(\hat{f'}:\XH_2\to\XH_1,g')$ is also metric graph bundle morphism. It is enough to check that $\hat{f}$ and $\hat{f'}$ are coarse inverse of each other.\\
Let $x\in V(X_1)$. Since $f$ is type preserving, $f(x)\in V(X_2)$. So, we have $d_{\XH_1}(\hat{f'}\circ\hat{f}(x),x)\leq d_{X_1}(f'\circ f(x),x)\leq D$.\\
Let $x=\nu(H_1)$, where $H_1\in\HH_1$. Then $d_X(f(H_1),H_2)\leq D_0$ for some $H_2\in\HH_2$. Then, $\hat{f}(x)=\nu(H_2)$ and $\hat{f'}(\nu(H_2))=\nu(H_1)$. Therefore, $d_{\XH_1}(\hat{f'}\circ\hat{f}(x),x)=0$. \\
Similarly, we can show that for all $y\in V(\XH_2)$, $d_{\XH_2}(\hat{f}\circ\hat{f'}(y),y)\leq D$.    
\end{proof}
\noindent
Now we prove the analogue of the results for pullback of relatively hyperbolic metric graph bundles.
\begin{prop}\label{pullback}
	Let $(X,B,\pi)$ be a relatively hyperbolic metric graph bundle. Let $A$ be a hyperbolic metric space and $g:A\to B$ be a coarsely $L$-Lipschitz map. Then there exists a pullback and the pullback is a relatively hyperbolic metric graph bundle. Moreover, the pullback map is type preserving.
\end{prop}
\begin{proof}
By \propref{pullback graph prop}, there exists a pullback graph bundle $p:Y\to A$ with the pullback map $f:Y\to X$. First, we recall the construction of $Y$ from the proof of \cite[Proposition 3.23]{krishna-sardar}. For each $b\in V(A)$, $F_b:=p^{-1}(b)$ is the full subgraph $F_{g(b)}:=\pi^{-1}(g(b))$.  For any pair of adjacent vertices $s,t\in V(A)$ we introduce some edges joining vertices of $F_s$ and $F_t$ in the following way. For each $s\in V(A)$, $F_s$ is an identical copy of $F_{g(s)}$. Let $f_s:F_s\to F_{g(s)}$ denote this identification.
Let $e$ be an edge joining $s,t$ in $A$ and $\alpha$ be a geodesic in $B$ joining $g(s),g(t)$. By \lemref{lifting geodesics}(1), for each $x\in V(F_s)$, $\alpha$ can be lifted to a geodesic in $X$ with $\tilde{\alpha}(g(s))=x$. We join $x$ by 
an edge to $y\in V(F_t)$ if and only if $f_t(y)=\tilde{\alpha}(g(t))$. This completes the construction of $Y$. \\
For each $b\in V(A)$, since $F_b$ is isometrically identified with $F_{g(b)}$, $F_b$ is strongly hyperbolic relative to the natural collection of subsets denoted by $\HH_b:=\{H_{b,\beta}\}_{\beta\in\Lambda}$. Let $b_1,b_2\in V(B)$ be adjacent vertices in $A$. Let $H_1$ be a horosphere-like subset in $F_{b_1}$ and let $x\in V(H_1)$. Then $x$ is joined by an edge to some $y\in V(F_{b_2})$. By the construction of $Y$, this means the following. Let $f_{b_2}:F_{b_2}\to F_{g(b_2)}$ be the map identifying the two fibers and $\alpha$ be a geodesic in $B$ joining $g(b_1)$ and $g(b_2)$ of length at most $2L$. If $\tilde{\alpha}$ is a lift of $\alpha$ in $X$ with $\tilde{\alpha}(g(b_1))=x$ then, $f_{b_2}(y)=\tilde{\alpha}(g(t))$. By \corref{lifting geodesics-2}, there exists $H_2\in\HH_{b_2}$ such that $f_{b_2}(y)\in N_{D_{\ref{lifting geodesics-2}}(2L)}(H_2)$. Identifying $\HH_{b_2}$ with the horosphere-like subsets of $F_{b_2}$, we have that $y$ lies in a uniform neighbourhood of $H_2$. This shows that the condition (3) also holds. It is easy to see that (4) is satisfied as well.\\
Thus, $(Y,A,p)$ is a relatively hyperbolic metric graph bundle. Since the identification maps collectively define the pullback map, it is type preserving. 
\end{proof}
\noindent
Let $\HH_Y=\{\HH_b\}_{b\in V(A)}$. We next show that coned-off space of $Y$ is still a pullback of the coned off space of $X$.
\begin{prop}\label{induced pullback}
	Let $(X,B,\pi)$ be a relatively hyperbolic metric graph bundle and $g:A\to B$ be a coarsely Lipschitz map. Let $(Y,A,p)$ be the pullback of $X$ under $g$ and $f$ be the pullback map. Then $\YH:=\EE(Y,\HH_Y)$ is the pullback of $\XH:=\EE(X,\HH)$ under $g$ with the induced map $\hat{f}:\YH\to\XH$ as the pullback map.
\end{prop}
\begin{proof}
	By \propref{new mbdl}, $\YH$ is a metric graph bundle. So it is enough to check the universal property. Let $\pi':X'\to A$ be a metric graph bundle with a morphism $(f',g)$ from $X'$ to $\XH$, where $f'$ is coarsely $L$-Lipschitz for some $L\geq 1$. Since $Y$ is the pullback of $X$ under $g$, we know that for each $s\in V(A)$, $p^{-1}(s)$ is $\pi^{-1}(g(s))$ and this identification is denoted by $f_s$. So we have a natural map $\widehat{f}_s:\EE(p^{-1}(s),\HH_s)\to\EE(\pi^{-1}(g(s)),\HH_{g(s)})$ which is a uniform quasiisometry. 
	\begin{figure}[ht]
		\centering
		\begin{tikzpicture}[node distance=1.5cm,auto]
			\node (P) {$\YH$};
			\node (B) [right of=P] {$\XH$};
			\node (A) [below of=P] {$A$};
			\node (C) [below of=B] {$B$};
			\node (P1) [node distance=1cm, left of=P, above of=P] {$X'$};
			\draw[->] (P) to node {$\hat{f}$} (B);
			\draw[->] (P) to node [swap] {$p$} (A);
			\draw[->] (A) to node [swap] {$g$} (C);
			\draw[->] (B) to node {$\pi$} (C);
			\draw[->, bend right] (P1) to node [swap] {$\pi'$} (A);
			\draw[->, bend left] (P1) to node {$f'$} (B);
			\draw[->, dashed] (P1) to node {$f''$} (P);
		\end{tikzpicture}
		\caption{}\label{coned-off pullback figure}\end{figure}
	
	We define $f'':X'\to\YH$ as follows. For each $s\in V(A)$, $f''_s=\widehat{f}^{-1}_s\circ{f'_s}$. Collectively, these maps define $f''$ and $f''$ makes Figure \ref{coned-off pullback figure} commutative.\\		  
	We need to show that $f''$ is a coarsely unique, coarsely Lipschitz map. This follows from \lemref{pullback lemma}. Note that condition (1) of \lemref{pullback lemma} holds, i.e., $\hat{f},f'$ are coarsely Lipschitz as they are metric graph bundle morphisms. Condition (2) follows from \lemref{fibers qi}(1), and 
	(3) follows since the fibers of metric graph bundles are uniformly properly embedded.	\end{proof}

\subsubsection{QI sections}
We recall the following result.
\begin{prop}\cite[Proposition 2.10, Proposition 2.12]{pranab-mahan}
	{\em(Global qi sections for metric (graph) bundles)}\label{existence-qi-section} For all $\delta',c\geq 0, N\geq 0$ and 
	$\eta:[0,\infty) \rightarrow [0,\infty)$ there exists $K_{\ref{existence-qi-section}}=K_{\ref{existence-qi-section}}(c,\eta,\delta',N)$ such that the following holds {\em (}$c=1$ for metric graph bundles{\em )}.\\	
	Suppose $\pi:X'\rightarrow B'$ is an $(\eta,c)$-metric bundle or an $\eta$-metric graph bundle with fibers such that for each $b\in B'$ \textup{(}$b\in V(B')$\textup{)},\\
	(i) $F_b$ is $\delta'$-hyperbolic with respect to the induced path metric.\\
	(ii) The barycenter map $\partial^3 F_b\to F_b$ is coarsely $N$-surjective.\\
	Then there is a $K_{\ref{existence-qi-section}}$-qi section over $B'$ through each point of $X'$.
\end{prop}
\noindent
As an easy corollary, we show the following.
\begin{corollary}\label{existence-qi-section-2} For all $\delta',c\geq 0, N\geq 0$ and $\eta:[0,\infty) \rightarrow [0,\infty)$ there exists $K_{\ref{existence-qi-section-2}}=K_{\ref{existence-qi-section-2}}(c,\eta,\delta',N)$ such that the following holds.\\	
	Suppose $\pi:X\rightarrow B$ is a $\eta$-relatively hyperbolic metric graph bundle such that for each $b\in V(B)$, the barycenter map $\partial^3\widehat{F}_b\to\widehat{F}_b$ is coarsely $N$-surjective. Then there is a $K_{\ref{existence-qi-section-2}}$-qi section in $\XH$ over $B$, through each point of $X$, such that it lies in $X$.
\end{corollary}
\begin{proof}
The space $\XH=\EE(X,\HH)$ is a metric graph bundle. Let $x\in V(F_b)$. By \propref{existence-qi-section}, there exists a $K_{\ref{existence-qi-section}}$-qi section	$s':B\to\XH$ with $s'(b)=x$. We define a map $s:B\to X$ as follows.\\
For $u\in V(B)$ with $s'(u)\in V(F_u)$, we put $s(u)=s'(u)$. If $s'(u)$ is a cone point in $\widehat{F}_u$, i.e., $s'(u)=\nu(H_u)$ for some $H_u\in\HH_u$, then choose $z\in V(H_u)$ and define $s(u)=z$. We need to check that $s$ is a qi section and that it is coarsely well defined.\\
(1) Clearly $s$ is a set-theoretic section. We need to check that it is a qi embedding. We know that for each $u\in V(B)$, $d_u(s(u),s'(u))\leq\frac{1}{2}$. 
Let $u,v\in V(B)$. Then $d_B(u,v)\leq d(s(u),s(v))$ and $d(s(u),s(v))\leq d(s'(u),s'(v))+1\leq K_{\ref{existence-qi-section}}d_B(u,v)+K_{\ref{existence-qi-section}}+1$. So, $s$ is a $(K_{\ref{existence-qi-section}}+1)$-qi embedding.\\
(2) Coarse well definedness of $s$ follows easily. Let $s_1$ be another qi section defined in way similar to that of $s$, i.e., for each $u\in B$, $d(s'(u),s_1(u))\leq d_u(s'(u),s_1(u))\leq\frac{1}{2}$. This implies $d(s,s_1)\leq 1$. Thus, for $K_{\ref{existence-qi-section-2}}=K_{\ref{existence-qi-section}}+1$, we have the proof.
\end{proof}	
\noindent
So through every point in $V(X)$, there exists a qi section $s:B\to\XH$ such that $s(B)\subset X$. 
Now we define a notion of cone locus in the metric graph bundle setting similar to that in \cite{mahan-reeves}.
\begin{definition}
	The {\bf cone locus} of $\XH$ is a subgraph $\BB$ of $\XH$ with vertex set $\{\nu(H_{b,\alpha})\mid b\in V(B),\alpha\in \Lambda\}$ and any two vertices $\nu(H_{b,\alpha})$ and $\nu(H_{b',\beta})$ are connected by an edge in $\BB$ if they are connected by an edge in $\XH$.
\end{definition}
\begin{remark}\label{about cone locus}
Each connected component of $\BB$ is an isometric copy of $B$. 
\end{remark}
\begin{convention}\label{qi section convention}
$(1)$ For adjacent vertices $b_1,b_2$ in $B$ and the fiber identification map $\phi_{b_1b_2}:F_{b_1}\to F_{b_2}$, by type preserving condition, $\phi_{b_1b_2}(H_{b_1,\alpha})$ lies uniformly close to an element of $\HH_{b_2}$. We denote this element by $H_{b_2,\alpha}$.\\
$(2)$ For each $\alpha\in\Lambda$, the connected component of $\BB$ consisting of vertices indexed by $\alpha$ will be denoted by $\BB_\alpha$. So, $\displaystyle{\BB=\bigsqcup_{b\in V(B)}\BB_\alpha}$.\\
\end{convention}
\noindent
Now we recall the definition of some terminology appearing in Theorem \ref{main theorem}. 
\begin{definition}\label{flare}\cite[Definition 1.12]{pranab-mahan}
	Suppose $\pi:X\rightarrow B$ is a metric (graph) bundle.
	$X$ satisfies a {\bf flaring condition} if for all   $k \geq 1$, there exist $\lambda_k>1$ and $n_k,M_k\in \NN$ such that
	the following holds:\\
	Let $\gamma:[-n_k,n_k]\rightarrow B$ be a geodesic and let
	$\tilde{\gamma}_1$ and $\tilde{\gamma}_2$ be two
	$k$-qi lifts of $\gamma$ in $X$.
	If $d_{\gamma(0)}(\tilde{\gamma}_1(0),\tilde{\gamma}_2(0))\geq M_k$,
	then we have
	\[\mbox{
		{\small $\lambda_k.d_{\gamma(0)}(\tilde{\gamma}_1(0),\tilde{\gamma}_2(0))\leq \mbox{max}\{d_{\gamma(n_k)}(\tilde{\gamma}_1(n_k),\tilde{\gamma}_2(n_k)),d_{\gamma(-n_k)}(\tilde{\gamma}_1(-n_k),\tilde{\gamma}_2(-n_k))\}$}}.
	\]
\end{definition}
\begin{definition}\label{cone flare} A coned-off metric (graph) bundle $\pi:\XH\rightarrow B$ satisfies a {\bf cone-bounded strictly flare condition} if for all $k \geq 1$, there exist $\lambda_k>1$ and $n_k\in \NN$ such that the following holds:\\
	Let $\gamma:[-n_k,n_k]\rightarrow B$ be a geodesic and let
	$\tilde{\gamma}_1$ and $\tilde{\gamma}_2$ be two
	lifts of $\gamma$ lying in the cone locus of $X$. Then 
	\[\mbox{
		{\small $\lambda_k.d_{\gamma(0)}(\tilde{\gamma}_1(0),\tilde{\gamma}_2(0))\leq \mbox{max}\{d_{\gamma(n_k)}(\tilde{\gamma}_1(n_k),\tilde{\gamma}_2(n_k)),d_{\gamma(-n_k)}(\tilde{\gamma}_1(-n_k),\tilde{\gamma}_2(-n_k))\}$}}.
	\]
\end{definition}

\noindent
For each $\alpha\in\Lambda$, $\BB_\alpha$ gives rise to a subgraph $\CC_\alpha$ of $X$, with $V(\CC_\alpha)=\bigsqcup_{b\in V(B)}V(H_{b,\alpha})$. Now, any pair of points $x,y\in V(\CC_\alpha)$ are joined by an edge in $\CC_\alpha$ if they are joined by an edge in $X$. We denote by $\CC=\{\CC_\alpha\mid\alpha\in\Lambda\}$. For each $\alpha\in\Lambda$, let $g_\alpha:\CC_\alpha\to\BB_\alpha$ be the natural map that collapses $\CC_\alpha$ to $\BB_\alpha$. Let $\GG=\{g_\alpha\mid\alpha\in\Lambda\}$. Note that each $\CC_\alpha$ is a metric graph bundle over $\BB$ with fibers $H_{b,\alpha}$.
Let $X_{pel}:=\mathcal{PE}(X,\CC,\GG,\BB)$. 
\begin{remark}\label{identifying spaces}
The space $\XH$ can be identified with $X_{pel}$.
\end{remark}

\subsection{Ladders}
Given a pair of qi sections, we can define a quasiconvex subset called a ladder. We recall its definition and properties here.
\begin{definition}\label{defn-ladder}\cite[Definition 2.13]{pranab-mahan}
Suppose $\Sigma_1$ and $\Sigma_2$ are two $K$-qi sections of the metric (graph) bundle $X$.
	For each $b\in V(B)$ we join the points $\SSS_1\cap F_b$, 
	$\SSS_2\cap F_b$ by a geodesic in $F_b$. We denote the union of these 
	geodesics by $\LL(\SSS_1,\SSS_2)$, and call it a $K$-{\bf ladder} {\em (}formed by the sections 
	$\SSS_1$ and $\SSS_2${\em )}.\end{definition}
\noindent
Note that in our case, a ladder is a disjoint union of electric fiber geodesics in $\widehat{F}_b$, $b\in V(B)$, and thus it is a subset of $\XH$.
\begin{definition}{\em(Retraction map)}
Let $\LL=\LL(\SSS_1,\SSS_2)$. We define $\Pi:\XH\to\LL$ as follows. For $x\in\XH$ with $\pi(x)=b$, $\Pi(x)=\hat{\pi}_{\hat{\lambda}_b}(x)$, where $\hat{\lambda}_b=\LL\cap\widehat{F}_b$ and $\hat{\pi}_{\hat{\lambda}_b}:\widehat{F}_b\to\widehat{\lambda}_b$ is an electric projection map.
\end{definition}
\noindent
We skip the proof of the following result as it is similar to that of \cite[Proposition 4.6]{krishna-sardar}. (Also see \cite{pranab-mahan}).
\begin{prop}\label{ladders are qi embedded}
	Given $K\geq 0$, $\delta\geq 0$ there exists $C=C_{\ref{ladders are qi embedded}}(K)\geq 0$,
	$R= R_{\ref{ladders are qi embedded}}(K)\geq 0$ and $K_{\ref{ladders are qi embedded}}(\delta, K)\geq 0$ 
	such that the following holds. 	
	Let $\pi:X\to B$ be a $\eta$-metric graph bundle and let $\SSS_1, \SSS_2$ be two $K$-qi sections in $X$ and $\LL=\LL(\SSS_1, \SSS_2)$. Then the following hold.\\	
	$(1)$ The map $\Pi:\XH\to\LL$ is a coarse $C$-Lipschitz retraction.\\	
	$(2)$ Given a $k$-qi section $\gamma$ in $\XH$ over a geodesic in $B$, $\Pi(\gamma)$ is a $(C+2kC)$-qi section in 
	$\XH$ contained in $\LL$ over the same geodesic in $B$.\\	
	$(3)$ If the barycenter map of the fibers of $\XH$ are uniformly coarsely surjective, then through any point of $\LL$ there is $(1+2K)C$-qi section contained in $\LL$.\\	
	$(4)$ The $R$-neighbourhood of $\LL$ is connected and uniformly qi embedded in $\XH$. In particular, if $\XH$ is $\delta$-hyperbolic then $\LL$ is $K_{\ref{ladders are qi embedded}}(\delta, K)$-quasiconvex in $\XH$.
\end{prop}
\begin{remark}
Let $\LL$ be a $K$-ladder as in \propref{ladders are qi embedded}. Through any $x\in\LL$, by Proposition \ref{ladders are qi embedded}(3), there exists a $(1+2K)C$-qi section $\SSS$ in $\XH$, through $x$, lying in $\LL$. This follows from \propref{existence-qi-section} and retraction map. So it is easy to see that if $x\in\LL\cap X$, then $\SSS\subset\LL\cap X$.
\end{remark}
\begin{convention}\label{qi section in ladder}
By \propref{ladders are qi embedded}, for any $K\geq 1$, and a $K$-ladder $\LL$, through any $x\in\LL$, there exists a $(C+2KC)$-qi section in $\XH$, contained in $\LL$. Let $K_0=K_{\ref{existence-qi-section-2}}$ and for $n\geq 0$, let $K_{n+1}=C+2K_nC$.  
\end{convention}

\begin{definition}\cite[Definitions 2.15, 2.16]{pranab-mahan}
Let $\pi:\XH\to B$ be a metric (graph) bundle and $\SSS_1,\SSS_2$ be two qi sections.\\	
	$(1)$ Suppose $R\geq 0$. The set $U_R(\SSS_1,\SSS_2)=\{b\in B:\,d_b(\SSS_1\cap F_b,\SSS_2\cap F_b)\leq R\}$ is called the $R$-{\bf neck} of the ladder $\LL(\SSS_1,\SSS_2)$.\\	
	$(2)$ The {\bf girth} of the
	ladder $\LL(\SSS_1,\SSS_2)$ is $\min\{d_b(\SSS_1\cap F_b, \SSS_2\cap F_b): b\in B\}$ and it is denoted by $d_h(\SSS_1,\SSS_2)$.
\end{definition}
\begin{remark}\cite[Lemma 2.18]{pranab-mahan}\label{remark about neck} Let $\LL(\SSS_1,\SSS_2)$ be a a $K$-ladder.
$(1)$ If $R\geq\max\{ M_K,d_h(\SSS_1,\SSS_2)\}$, then $U_R(\SSS_1,\SSS_2)$ is a quasiconvex subset, where the quasiconvexity constant depends only on $K$.\\
$(2)$ If $d_h(\SSS_1,\SSS_2)\geq M_K$, then the diameter of $U_R(\SSS_1,\SSS_2)$ is uniformly bounded by a constant depending on $K$ and $R$.\\
$(3)$ For any $L>0$ and a geodesic $\gamma:[t_0,t_1]\to B$ with $d_{\gamma(t_0)}(\SSS_1\cap F_{\gamma(t_0)},\SSS_2\cap F_{\gamma(t_0)})=LR$, $\gamma(t_1)\in U_R(\SSS_1,\SSS_2)$ and for all $t\in [t_0,t_1)\cap\ZZ$, $\gamma(t)\not \in U_R(\SSS_1,\SSS_2)$, we have $l_{\XH}(\gamma)$ is uniformly bounded by a constant depending on $K$ and $L$.
\end{remark} 
\begin{definition}
	{\bf (Small girth ladders)} A $K$-ladder $\LL(\SSS_1, \SSS_2)$ is called a {\bf small girth ladder} if $U_{M_K}(\SSS_1,\SSS_2)\neq\emptyset$.
\end{definition}
\noindent

\section{The combination theorem}\label{section 4}
\noindent
In this section, we prove a combination theorem for metric (graph) bundles with relatively hyperbolic fibers. Here is the set-up. 
\begin{enumerate}[leftmargin=\parindent,align=left,labelwidth=\parindent,labelsep=4pt]
	\item  $\pi:X\rightarrow B$ be an $\eta$-metric (graph) bundle with $B$, $\delta_0$-hyperbolic.
	\item For each $b\in V(B)$, $\HH_b=\{H_{b,\alpha}\}_{\alpha\in\Lambda}$ is a collection of uniformly mutually cobounded subsets and the fiber $F_b$ strongly hyperbolic relative to $\HH_b$. Moreover, $\widehat{F}_b=\mathcal{E}(F_b,\HH_b)$ is $\delta_0$-hyperbolic with respect to $\hat{d}_b$. 
	\item For each $b\in V(B)$, the barycentre map $\partial^3\widehat{F}_b\to \widehat{F}_b$ is coarsely $N$-surjective.	
	\item The natural fiber maps are uniformly type preserving and satisfy qi-preserving electrocution condition.
	\item The induced coned-off metric graph bundle $\XH$ satisfies the $(\nu_k, M_k,n_k)$-flaring condition, for each $k\geq 1$, and also a cone-bounded strictly flare condition.
\end{enumerate} 
\noindent
We recall the following theorem from \cite{pranab-mahan}.
\begin{theorem}{\em (\cite[Theorem 4.3 and Proposition 5.8]{pranab-mahan})}\label{mbdl thm}
Suppose $\pi:X\to B$ is a metric (graph) bundle such that\begin{enumerate}[leftmargin=\parindent,align=left,labelwidth=\parindent,labelsep=4pt]
	\item $B$ is a hyperbolic metric space.
	\item Each fiber $F_b$, $b\in B$ (resp. $b\in V(B)$) is a uniformly hyperbolic metric space.
	\item For each fiber, the barycentre map $\partial^3 F_b\to F_b$ is uniformly coarsely surjective.
\end{enumerate}
Then $X$ is a hyperbolic metric space if and only if a flaring condition is satisfied.
\end{theorem}
\noindent
A combination theorem for a relatively hyperbolic extension of surface groups was already given by Mj and Sardar in \cite{pranab-mahan}.
\begin{theorem}\cite[Proposition 5.17]{pranab-mahan}
	Let $K=\pi_1(S^h)$ be the fundamental group of a surface with finitely many punctures and let $K_1,\ldots,K_n$ be its peripheral subgroups. Let $1\to K\to G\stackrel{\pi}{\to} Q\to1$ and $1\to K_i\to N_G(K_i)\stackrel{\pi}{\to}Q_i\to 1$ be the induced short exact sequences of groups. Then $G$ is strongly hyperbolic relative to $\{N_G(K_i)\}^n_{i=1}$ if and only if $Q$ is a convex cocompact subgroup of the pure mapping class group of $S^h$. 
\end{theorem}
\noindent
We first prove the following.

\subsection{Step 1: Weak relative hyperbolicity}
\begin{theorem}\label{weak comb}
Let $\pi:X\rightarrow B$ be an $\eta$-metric (graph) bundle such that \begin{enumerate}
	\item[(C1)] $B$ is a $\delta_0$-hyperbolic metric space.
	\item[(C2)] For each $b\in V(B)$, the fiber $F_b$ strongly hyperbolic relative to a collection of uniformly mutually cobounded subsets $\HH_b$.
	\item[(C3)] For each $b\in V(B)$, the barycentre map $\partial^3\widehat{F}_b\to \widehat{F}_b$ is coarsely $N$-surjective.	
	\item[(C4)] The induced coned-off metric (graph) bundle satisfies the $(\nu_k, M_k,n_k)$-flaring condition, for each $k\geq 1$.
\end{enumerate}
Then $X$ is weakly hyperbolic relative to $\CC$.	
\end{theorem}
\begin{proof}
Recall that $\XH$ is a metric graph bundle over the hyperbolic base space $B$ with $\delta_0$-hyperbolic fibers $\widehat{F}_b$ for each $b\in V(B)$. Then along with conditions (C3) and (C4), $\XH$ satisfies all the hypotheses of \thmref{mbdl thm}. Therefore, $\XH$ is a hyperbolic metric space. \\
Now, $\BB$ is a collection of uniformly separated, uniformly quasiconvex subsets of $\XH$. Therefore, by \lemref{coned-off hyp}, $\EE(\XH,\BB)$ is a hyperbolic metric space. Therefore, $(\XH,\BB)$ is weakly relatively hyperbolic.  By Remark \ref{coned off-pel qi}, $\EE(\XH,\BB)$ is is quasiisometric to $\EE(X,\CC)$ and thus, $\EE(X,\CC)$ is a hyperbolic metric space. Therefore, $X$ is weakly hyperbolic relative to $\CC$. 
\end{proof}
\subsection{Step 2: Strong relative hyperbolicity}
For the strong relative hyperbolicity of $(X,\CC)$, the pair needs to satisfy the bounded penetration property. We use the following lemma to prove strong relative hyperbolicity. The proof is similar to that of \cite[Proposition 4.4]{mahan-reeves}.
\begin{lemma}\label{cone flare implies bcp}
	Suppose the metric (graph) bundle $\XH$ is hyperbolic and it satisfies a cone-bounded strictly flare condition. Then there exists $D>0$ such that $\BB$ is a mutually $D$-cobounded collection of subsets.
\end{lemma}
\begin{proof}
Suppose not. Then for every $D\geq D_{\ref{big npp small nbd}}$, there exists $\BB_\alpha,\BB_\beta\in\BB$, $x,y\in \BB_{\alpha}$ and their nearest point projections $x_1,y_1\in \BB_\beta$ respectively such that $d(x_1,y_1)\geq D$. Then by \lemref{big npp small nbd}, $[x_1,y_1]\subset N^{\XH}_{M_{\ref{big npp small nbd}}}([x,y])$. Let $x_2,y_2\in[x,y]$ such that $d_{\XH}(x_1,x_2)\leq M_{\ref{big npp small nbd}},d_{\XH}(y_1,y_2)\leq M_{\ref{big npp small nbd}}$. Clearly, $d_{\XH}(x_2,y_2)\geq d_{\XH}(x_1,y_1)-d_{\XH}(x_1,x_2)-d_{\XH}(y_1,y_2)\geq D-2M_{\ref{big npp small nbd}}$. For $i=1,2$, let $\gamma_i=[x_i,y_i]$.

Consider $\pi(\gamma_i)$ in $B$ for $i=1,2$. Since each $\BB_\alpha\in\BB$ is an isometric embedding of $B$ in $\XH$, $d_B(\pi(x_i),\pi(y_i))\geq D-2M_{\ref{big npp small nbd}}$ and moreover, $d_B(\pi(x_1),\pi(x_2))\leq M_{\ref{big npp small nbd}}, d_B(\pi(y_1),\pi(y_2))\leq M_{\ref{big npp small nbd}}$. Let $x_0,y_0$ denote the lifts of $\pi(x_1),\pi(y_1)$ in $\BB_\alpha$. Then $d_{\XH}(x_1,x_0)\leq d_{\XH}(x_1,x_2)+d_{\XH}(x_2,x_0)=d_{\XH}(x_1,x_2)+d_B(\pi(x_2),\pi(x_1))\leq 2M_{\ref{big npp small nbd}}$. Similarly, $d_{\XH}(y_1,y_0)\leq 2M_{\ref{big npp small nbd}}$. Note that the geodesic $\gamma_0:=[x_0,y_0]$ in $\BB_\alpha$ has length at least $D-2M_{\ref{big npp small nbd}}$.\\
So, we have geodesics $\gamma_0\subset\BB_{\alpha},\gamma_1\subset\BB_\beta$ of length at least $D-2M_{\ref{big npp small nbd}}$ such that $\pi(x_1)=\pi(x_0)$ and $\pi(y_1)=\pi(y_0)$. $Hd_{\XH}(\gamma_0,\gamma_1)\leq 2M_{\ref{big npp small nbd}}+2\delta$, where $\delta$ is the hyperbolicity constant of $\XH$. So, for each $x'\in\gamma_1$, there exists $y'\in\gamma_0$ such that $d_{\XH}(x',y')\leq 2M_{\ref{big npp small nbd}}+2\delta$. Since $d_{B}(\pi(x'),\pi(y'))\leq 2M_{\ref{big npp small nbd}}+2\delta$, there exists $y''\in\gamma_0$ with $\pi(y'')=\pi(y')$ such that $d_{\XH}(x',y'')\leq 4(M_{\ref{big npp small nbd}}+\delta)$ and thus, $d_{\F_{\pi(x)}}(x',y'')\leq\eta_{\ref{new mbdl}}(4(M_{\ref{big npp small nbd}}+\delta))=:M_1$. So, there exists a cone-bounded hallway $f:[-m,m]\times [0,1]\to\XH$ of length $2m\geq D-4(M_{\ref{big npp small nbd}}+\delta)$ with $f([-m,m]\times\{0\})=\gamma_1$, $f([-m,m]\times\{1\})=\gamma_0$ satisfying 
$$l(f(\{i\}\times[0,1]))\leq M_1, \text{ for } -m\leq i\leq m.$$
As we increase $D$, $m$ also increases, but $M_1$ is fixed. Thus, for every $\lambda>1$, there exists a cone-bounded hallway which does not flare, violating the cone-bounded hallway strictly flare condition. \end{proof}
\begin{theorem}\label{strong comb}
	Let $\pi:X\rightarrow B$ be an $\eta$-metric (graph) bundle such that \begin{enumerate}[leftmargin=\parindent,align=left,labelwidth=\parindent,labelsep=0pt]
		\item[(C1) ] $B$ is a $\delta_0$-hyperbolic metric space.
		\item[(C2) ] Each fiber $F_b$ is strongly hyperbolic relative to a collection of uniformly mutually cobounded subsets $\HH_b$.
		\item[(C3) ] For each fiber, the barycentre map $\partial^3\widehat{F}_b\to \widehat{F}_b$ is coarsely $N$-surjective.	
		\item[(C4) ] The induced coned-off metric graph bundle satisfies the $(\nu_k, M_k,n_k)$-flaring condition, for each $k\geq 1$.
		\item[(C5) ] The cone bounded strictly flare condition is satisfied.
	\end{enumerate}
	Then $X$ is strongly hyperbolic relative to $\CC$.	
\end{theorem}
\begin{proof}
By \thmref{weak comb} and \lemref{cone flare implies bcp}, we have that $\XH$ is a hyperbolic metric space and $\BB$ are mutually $D$-cobounded subsets of $\XH$. Then, by \lemref{coned-off strong hyp}, $(\XH,\BB)$ is strongly relatively hyperbolic. Equivalently, $(X,\CC)$ is also strongly relatively hyperbolic.
\end{proof}
\noindent
We end this section with the converse to the strong combination theorem. The proof is similar to that of \cite[Theorem 4.7]{mahan-reeves}. 
\begin{theorem}\label{Combination converse}
Let $\pi:X\rightarrow B$ be an $\eta$-metric (graph) bundle such that \begin{enumerate}[leftmargin=\parindent,align=left,labelwidth=\parindent,labelsep=4pt]
	\item[(C1)] $B$ is a $\delta_0$-hyperbolic metric space.
	\item[(C2)] For each $b\in V(B)$, the fiber $F_b$ strongly hyperbolic relative to a collection of uniformly separated closed subsets $\HH_b$.
	\item[(C3)] For each $b\in V(B)$, the barycentre map $\partial^3\widehat{F}_b\to \widehat{F}_b$ is coarsely $N$-surjective.	
\end{enumerate}
If $(X,\CC)$ is strongly relatively hyperbolic, then the induced coned-off metric (graph) bundle satisfies a flaring condition and a cone bounded strictly flare condition.
\end{theorem}
\begin{proof}
Since $(X,\CC)$ is strongly relatively hyperbolic, by \lemref{brahma-amalgeo}, $\PP\EE(X,\CC,\BB,\GG)=\XH$ is a hyperbolic metric space with $\BB$, is a collection of quasiconvex subsets of $\PP\EE(X,\CC,\BB,\GG)$. This, along with the conditions (C1),(C2),(C3) and \thmref{mbdl thm}, we have that $\XH$ satisfies a flaring condition. \\
Now we will show that cone-bounded flaring condition is satisfied. Suppose not. Then there exists some $k\geq 1$ and $M\geq 1$ such that for all $n>0$, we have the following. Let $\gamma:[-n,n]\rightarrow B$ be a geodesic with $k$-qi lifts
$\tilde{\gamma}_1$ and $\tilde{\gamma}_2$ lying in the cone locus $\BB_{\alpha}$ and $\BB_{\beta}$ respectively satisfying $$d_{\gamma(i)}(\tilde{\gamma}_1(i),\tilde{\gamma}_2(i))\leq M \text{ for all } i\in[-n,n].$$ Let $\mu$ be an electric geodesic in $\F_{\gamma(-n)}$ joining $\tilde{\gamma}_1(-n)$ and $\tilde{\gamma}_2(-n)$. Similarly, let $\mu'$ be an electric geodesic in $\F_{\gamma(n)}$ joining $\tilde{\gamma}_1(n)$ and $\tilde{\gamma}_2(n)$. Then length of $\mu$ in $\F_{\gamma(-n)}$ and that of $\mu'$ in $\F_{\gamma(n)}$ is at most $M$. Let $a=\mu\cap\CC_{\alpha}$ and $a'=\mu'\cap\CC_{\alpha}$.

Let $\sigma_1$ be a concatenation of the edge $[a,\tilde{\gamma}_1(-n)]$, $\tilde{\gamma}_1$ and the edge $[\tilde{\gamma}_1(n),a']$. Let $\sigma_2$ be the concatenation of $\mu|_{[a,\tilde{\gamma}_2(-n)]}$, $\tilde{\gamma}_2$ and $\mu'|_{[\tilde{\gamma}_2(n),a']}$. Since $M$ is fixed, by Lemma \ref{extend quasigeodesic}, both $\sigma_1$ and $\sigma_2$ are partially electrocuted quasigeodesics in $\XH$. The paths $\sigma_1$, $\sigma_2$ pass through different elements of $\BB$ for $2n$ distance. As $n$ can be arbitrarily large, we get a contradiction to the similar intersection pattern. 
\end{proof}

\section{Cannon-Thurston maps for pullback bundles}\label{section 5}
\noindent
In this section, we complete the proof of the main result (Theorem \ref{main theorem}) of the paper. We have the following assumptions. \begin{itemize}
	[leftmargin=\parindent,align=left,labelwidth=\parindent,labelsep=4pt]
	\item $(X,B,\pi)$ is a metric graph bundle or an approximating metric graph bundle such that cone-bounded strictly flare condition is satisified and the induced coned-off metric graph bundle satisfies the $(\lambda_k, M_k,n_k)$-flaring condition, for each $k\geq 1$. 
		\item $B$ is a $\delta_0$-hyperbolic metric space.
		\item For each $b\in V(B)$, $\HH_b=\{H_{b,\alpha}\}_{\alpha\in\Lambda}$ is a collection of uniformly mutually cobounded subsets and the fiber $F_b$ strongly hyperbolic relative to $\HH_b$. Also, 
		$\widehat{F}_b=\mathcal{E}(F_b,\HH_b)$ is $\delta_0$-hyperbolic.
		\item For each $b\in V(B)$, the barycentre map $\partial^3\widehat{F}_b\to \widehat{F}_b$ is coarsely $N$-surjective.	
		\item The natural fiber maps are uniformly type preserving and satisfy qi-preserving electrocution condition.
		\end{itemize} 
\noindent
We recall the main theorem of \cite{krishna-sardar} here. 
\begin{theorem}\label{CT for mgbdl}\cite[Theorem 5.2]{krishna-sardar}
Suppose $X$ is a metric (graph) bundle over a hyperbolic metric space $B$ such that $X$ is hyperbolic and all the fibers are nonelementary uniformly hyperbolic, and the barycenter map is uniformly coarsely surjective. Suppose $i:A\map B$ is a Lipschitz, qi embedding 
and $Y$ is the pullback of $X$ under $i$. Then $Y$ is hyperbolic and $f:Y\to X$ admits the CT map.
\end{theorem} 	
\noindent
By \thmref{strong comb} in \secref{section 4}, $X$ is hyperbolic relative to the family $\mathcal{C}$ of maximal cone-subbundles of horosphere-like spaces. We will show that it is enough to prove the existence of CT maps in the case of metric graph bundles. In the case of relatively hyperbolic spaces, quasiisometry need not induce a homeomorphism between the Bowditch boundaries. However, in our case, this issue does not arise. Suppose $(X',B',\pi')$ be a metric bundle strongly hyperbolic relative to a collection of horosphere-like subsets $\CC'$, satisfying the hypotheses of \thmref{main theorem}, and let $(X,B,\pi)$ be the approximating metric graph bundle along with quasiisometries $\psi_B:B'\to B$ and $\psi_X:X'\to X$ with  coarse inverse $\phi_B$ and $\phi_X$ respectively, as in \propref{bundle vs graph bundle}. Then we have the following.
\begin{prop}\label{type preserving bundle morphism}
The metric graph bundle $(X,B,\pi)$ is strongly hyperbolic relative to $\CC=\{\psi_X(\CC'_{\alpha})\mid\CC'_{\alpha}\in\CC'\}$.  
\end{prop}
\begin{proof}
We can easily check that $(X,B,\pi)$ satisfies the hypotheses of \thmref{strong comb}. (C1) is obvious and (C2) holds by \propref{approx bundle rel hyp} and \lemref{for mg approx}. Similarly, (C3) also holds. Finally (C4) and (C5) follow from the assumptions. Thus, $X$ is strongly hyperbolic relative to $\CC$.
\end{proof}
\begin{prop}\label{induced qi bundle}
	The map $\psi^h:X'^h\to X^h$ induced by $\psi$ is a  quasiisometry. 
\end{prop}
\noindent
This follows from \lemref{type preserve implies qi}. 
Thus, existence of CT map for the pullback of the approximating metric graph bundle $(X,B,\pi)$ implies the existence of the CT map for the pullback of the metric bundle $(X',B',\pi')$.

\noindent
Suppose $i:A\to B$ is a $k$-qi embedding. As in the case of \cite[Theorem 5.2]{krishna-sardar}, we take $Y$ to be a subbundle of $X$. 
\begin{remark}\label{about coned-off Y}
$(1)$ Note that the pullback map $f$ is type preserving by \propref{pullback}. So for each $b\in V(A)$, the collection of horosphere-like subsets of $F_b$ is identified with $\HH_{i(b)}\subset F_{i(b)}$. Let $\HH_Y=\{\HH_{i(b)}\}_{b\in V(A)}$.\\
$(2)$ $\CC_Y$ is a collection of subgraphs $\CC_{Y,\alpha}$ of $Y$ where $V(\CC_{Y,\alpha})$ is the disjoint union of vertex sets of the horosphere-like subsets of the fibers of $Y$. Moreover, any two such vertices are joined by an edge in $\CC_Y$ if they are joined by an edge in $Y$. So each $\CC_{Y,\alpha}$ can be identified with $\CC_\alpha\cap Y$.   
\end{remark}

\begin{lemma}\label{rel hyp of Y}
The metric graph bundle $(Y,A,p)$ is strongly hyperbolic relative to a collection of subsets $\CC_Y$.
\end{lemma}
\begin{proof}
By \propref{pullback}, $(Y,A,p)$ is a relatively hyperbolic metric graph bundle. Note that $A$ is a hyperbolic metric space. The space $\YH$ is obtained as in \propref{new mbdl}. Hyperbolicity of $\YH$ follows from \thmref{CT for mgbdl}. 
Since each fiber $F_b$ in $Y$, $b\in V(A)$, is isometric to the fiber $F_{i(b)}$ in $X$, $F_b$ is strongly hyperbolic relative to a collection of subsets $\HH_b$. Moreover, the barycenter map $\partial^3\widehat{F}_b\to F_b$ is uniformly coarsely surjective since the barycenter map of $F_{i(b)}$ is coarsely surjective. By \cite[Remark 4.4]{pranab-mahan} and bounded flaring condition, $\EE(Y,\HH_Y)$ inherits the flaring condition of $\EE(X,\HH)$ and similarly, $(Y,\HH_Y)$ inherits the cone-bounded strictly flare condition of $(X,\HH)$. By \thmref{mbdl thm}, this completes the proof.
\end{proof}

\noindent
By Lemma \ref{CT existence-1}, to prove the existence of CT map, the main step is the construction of quasigeodesics in $\EE(Y,\HH_Y)$ and $\EE(X,\HH)$ to show that Mitra's criterion holds. The result below follows from \propref{induced pullback}.
\begin{prop}\label{coned-off pullback}
	Suppose $\pi:X\map B$ is a metric (graph) bundle satisfying the hypotheses of Theorem \ref{main theorem}. Let $i: A\map B$ is a Lipschitz $k$-qi embedding and $p:Y\map A$ be the pullback bundle with the pullback map $f:Y\map X$. Then $\YH$ is a pullback of $\XH$. 
\end{prop}
\noindent
Let $\displaystyle{\AA=\sqcup_{\alpha}\AA_\alpha}$ be the cone locus of $Y$ and let $\GG_Y$ be the collection of natural maps $g_{Y,\alpha}$ collapsing $C_{Y,\alpha}$ to $\AA_\alpha$, $\alpha\in\Lambda$. Let $Y_{pel}=\mathcal{PE}(Y,\CC_Y,\GG_Y,\AA)$. We identify $Y_{pel}$ and $\YH$, in fact, we identify $\YH$ with a subset of $\XH$. By the result below, $\YH$ is properly embedded in $\XH$.

\begin{lemma}\label{prop-embed}\cite[Lemma 5.20]{krishna-sardar}
	The pullback of a metric (graph) bundle is metrically properly embedded in it. 
\end{lemma}
\noindent
By \cite{krishna-sardar}, we have a construction of quasigeodesics in $\XH$ and $\YH$. Let $y,y'\in Y$. Let $\SSS,\SSS'$ be the $K_0$-qi sections in $\XH$ containing $y,y'$ respectively. The existence of these qi sections is due to \propref{existence-qi-section-2}. We denote the ladder formed by these qi sections by $\LL=\LL(\SSS,\SSS')$. By the main construction in \cite{pranab-mahan}, we have a uniform quasigeodesic $c(y,y')$  joining $y,y'$ in $\XH$ such that $c(y,y')\subset\LL$. Note that $\LL\cap\YH$ is a ladder in $\YH$. In \cite{krishna-sardar}, the authors `modify' $c(y,y')$ to get a uniform quasigeodesic in the pullback. So, in our case, for $y,y'\in Y$, the path $c(y,y')$ is modified to get a path $\bar{c}(y,y')$ in $\YH$ such that $\bar{c}(y,y')\subset\LL\cap\YH$. We make a note of the steps involved, without going into details. The reader is referred to \cite[Section 5]{krishna-sardar} for the detailed construction of these paths. 
We start with the ladder decomposition. We will retain the notations and conventions from \cite{krishna-sardar}. The ladder $\LL$ is decomposed into small girth subladders which facilitate the construction of quasigeodesics.
The following result encapsulates the ladder decomposition.
\begin{prop}\label{ladder subdivision}\cite[Corollary 5.10]{krishna-sardar} There exists constants $R_0>0, R_1, R'_1, D_1$ and there is a partition $0=t_0<t_1<\cdots <t_n=l$ of $[0,l]$ and $K_1$-qi sections $\SSS_i$ passing through $\alpha(t_i)$, $0\leq i\leq n$ inside $\LL(\SSS, \SSS')$ such that the following holds.\\
	$(1)$ $\SSS_0=\SSS, \SSS_n=\SSS'$.\\
	$(2)$ For $0\leq i\leq n-2$, $\SSS_{i+1}\subset \LL(\SSS_i, \SSS')$.\\
	$(3)$ For $0\leq i\leq n-2$, there are two possibilities. $(\textup{I})$ $d_h(\SSS_i, \SSS_{i+1})=R_0$,\\
	$(\textup{II})$ $d_h(\SSS_i, \SSS_{i+1})>R_0$ and there is a $K_2$-qi section $\SSS'_i$ through $\alpha(t_{i+1}-1)$ inside $\LL(\SSS_i, \SSS_{i+1})$ such that $d_h(\SSS_i, \SSS'_i)<R_1$.		
		In either case $d_X(\SSS_i, \SSS_{i+1})>R'_1$ and $\SSS_i, \SSS_{i+1}$ are $D_1$-cobounded in $\XH$.\\
	$(4)$ $d_h(\SSS_{n-1}, \SSS_n)\leq R_0$.
\end{prop}
\noindent
Now, we recall the construction of quasigeodesics in $\XH$. Put $y_0=y$ and $y_{n+1}=y'$. For $1\leq i\leq n$, let $y_i$ be a uniform approximate nearest point projection of $y_{i-1}$ on $\SSS_i$ in $\XH$. The points $y_i,y_{i+1}$ are joined in $\XH$ by a path $\gamma_i$ such that the concatenation of such $\gamma_i$'s is a uniform quasigeodesic in $X$ joining $y,y'$.

\smallskip
\noindent
The path construction is by induction. Suppose we have $y_0,\ldots,y_i$ and $\gamma_0,\ldots,\gamma_{i-1}$, $0\leq i\leq n-2$.\\
Suppose $\LL_i= \LL(\SSS_{i},\SSS_{i+1})$ is of type (I) or $i=n-1$. We know that in this case, $U_{R_0}(\SSS_{i},\SSS_{i+1})\neq\emptyset$. Let $u_i$ be a nearest point projection of $\pi(y_i)$ on $U_{R_0}(\SSS_{i},\SSS_{i+1})$. Then we put $y_{i+1}=\SSS_{i+1}\cap \widehat{F}_{u_i}$. Let $\alpha_i$ denote the lift of $[\pi(y_i),u_i]$ in $\SSS_i$, and $\sigma_i$ denote the subsegment of $\widehat{F}_{u_i}\cap \LL_i$ joining $\alpha_i(u_i)$ and $y_{i+1}$. The path $\gamma_i$ is the concatenation of $\alpha_i$ and $\sigma_i$. \\
Now suppose $\LL_i= \LL(\SSS_{i},\SSS_{i+1})$ is of type (II), i.e. $d_h(\SSS_{i},\SSS_{i+1})>R_0$. In this case, there is a $K_2$-qi section $\SSS'_i$ in $\LL_i = \LL(\SSS_i,\SSS_{i+1})$ passing through $\alpha(t_{i+1}-1)$ with $d_h(\SSS_i, \SSS'_i)\leq R_1$. Here the earlier case is repeated twice. Let $v_i$ be a nearest point projection of $\pi(y_i)$ on $U_{R_1}(\SSS_i,\SSS'_i)$ and $w_i$ be a nearest point projection $v_i$ on $U_{R_1}(\SSS'_i, \SSS_{i+1})$. Let $y'_i=\SSS'_i\cap \widehat{F}_{v_i}$ and $y_{i+1}=\SSS_{i+1}\cap \widehat{F}_{w_i}$. Let $\alpha_i$ be the lift of $[\pi(y_i),v_i]$ in $\SSS_i$ and let $\beta_i$ denote the lift of $[v_i,w_i]$ in $\SSS'_i$.
Then $\gamma_i$ is the concatenation of $\alpha_i$, 
$[\SSS_i\cap\widehat{F}_{v_i},\SSS'_i\cap\widehat{F}_{v_i}]_{v_i}$, $\beta_i$ and 
$[\SSS'_i\cap\widehat{F}_{w_i},\SSS_{i+1}\cap\widehat{F}_{w_i}]_{w_i}$.\\
Finally, for $i=n$, $\gamma_n$ is the lift of $[\pi(y_n),\pi(y')]$ in $\SSS'$. We denote the concatenation of ${\gamma}_i$'s by $c(y,y')$.
%

\smallskip
\noindent
The `modification' of $c(y,y')$ to get a path in $\YH$ is as follows. \\
For $0\leq i\leq n$, let $b_i$ be a nearest point projection of $\pi(y_i)$ on $A$ and let $\bar{y}_i=\widehat{F}_{b_i}\cap \SSS_i$.
We define a path $\bar{\gamma}_i \subset\YH$ joining the points $\bar{y}_i, \bar{y}_{i+1}$ for $0\leq i\leq n$. 
The path $\bar{\gamma}_n$ is the lift of $[\pi(y_{n+1}), \pi(\bar{y}_n)]_A$ in $\SSS'\cap\YH$. We define the rest of the $\bar{\gamma}_i$'s here.\\
Suppose $\LL_i$ is of type (I) or $i=n-1$. 
Let $\bar{\alpha}_i$ denote the lift of $[b_{i},b_{i+1}]_A$ in $\SSS_{i}$ starting at $\bar{y}_{i}$. The path $\bar{\gamma}_i$ is defined to be the concatenation of $\bar{\alpha}_i$ and the fiber geodesic $\widehat{F}_{b_{i+1}}\cap \LL(\SSS_{i},\SSS_{i+1})$.\\
Suppose $\LL_i$ is of type (II). Let $b'_i\in A$ be a nearest point projection $\pi(y'_i)$ on $A$ and $\bar{y}'_i=\widehat{F}_{{b}'_i}\cap\SSS'_i$.
Let $\bar{\alpha}_i$ and $\bar{\beta}_i$ be the lift of $[b_i,b'_i]_A$ in $\SSS_i\cap\YH$ and $[b'_i,b_{i+1}]_A$ in $\SSS'_i\cap\YH$ respectively. 
The concatenation of the paths $\bar{\alpha}_i$, 
$[\SSS_i\cap\widehat{F}_{b'_i}, \SSS'_i\cap\widehat{F}_{b'_i}]_{b'_i}$, $\bar{\beta}_i$ and 
$[\SSS'_i\cap\widehat{F}_{b_{i+1}}, \SSS_{i+1}\cap\widehat{F}_{b_{i+1}}]_{b_{i+1}}$ is 
defined to be $\bar{\gamma}_i$. The concatenation of $\bar{\gamma}_i$'s is denoted by $\bar{c}(y,y')$
\begin{lemma}\cite[Lemmas 5.12,5.18, Proposition 2.38]{krishna-sardar}\label{main quasigeodesics}
	There exists $K_{\ref{main quasigeodesics}}\geq1$ such that the path $c(y,y')$ (resp. $\bar{c}(y,y')$) is a $K_{\ref{main quasigeodesics}}$-quasigeodesic in $\XH$ (resp. $\YH$).
\end{lemma}
\noindent	
We will use these paths to show Mitra's criterion. To that end, we borrow the notion of vertical quasigeodesic rays from \cite{mahan-pal}.
\subsection{Vertical quasigeodesic ray}\label{vertical qg ray}
Suppose we have a metric graph bundle as in \thmref{main theorem}. Let $\LL=\LL(\SSS,\SSS')$ be a $K_0$-ladder in $\XH$ and $\LL$ has a decomposition as in \propref{ladder subdivision}. Fix $b_0\in V(B)$. Let $b\in B$ such that $d_B(b_0,b)=n$. Let $b_0=u_0,u_1,\ldots,u_n=b$ be the consecutive vertices in the geodesic $[b_0,b]$. We first fix some notations. For each $0\leq j\leq n$, we have the following.\\
(i) $F_j:=F_{u_j}$ and $\widehat{F}_j:=\widehat{F}_{u_j}$\\
(ii) $\widehat{\lambda}_j=\LL\cap\widehat{F}_j$ is an electric geodesic in $\widehat{F}_j$ and $\lambda^b_j=\widehat{\lambda}_j\cap F_j$.\\
(iii) $\phi_j:F_j\to F_{j-1}$ is a fiber identification map with the induced $L$-quasiisometry $\widehat{\phi}_j:\widehat{F}_j\to\widehat{F}_{j-1}$, where $L\geq 1$. Recall that $\widehat{\phi}_j(z)=\phi_j(z)$ for all $z\in F_j$.\\
(iv) A geodesic joining $x_1,x_2\in\widehat{F}_i$, for any $0\leq i\leq n$ is denoted by $[x_1,x_2]$.\\
Let $x\in\widehat{\lambda}_n\cap F_n$. We consider the following cases.\\
{\bf Case 1:} $x\in\QQ$ for $\QQ=\SSS_j$ or $\SSS'_j$ from the ladder decomposition. Then $\QQ$ is a $K_2$-qi section. Suppose $y_i=\QQ\cap F_{i}$ for $1\leq i\leq n-1$. Let $\mu_n$ be the maximal connected component of $\lambda^b_n$ containing $x$. Suppose $\QQ\neq\SSS$. We consider two possibilities. Suppose $\mu_n$ is a non-trivial segment. Let $a_n\in\mu_n$ such that $d_{F_n}(x,a_n)\leq 1$. Without loss of generality, assume $a_n\in\LL(\SSS,\QQ)$. Let $\mu'_n$ be the subsegment of $\mu_n$ joining $a_n,x$. Now, suppose $\mu_n$ is the single point $x$, i.e, there exists horosphere-like subsets $H_n,H'_n\in\HH_{b_n}$ such that $x\in H_n\cap H'_n$. In this case, take $a_n=\nu(H_n)$ and $\mu'_n$ to be the edge joining $x$ and $\nu(H_n)$, where $\nu(H_n)\in\widehat{\lambda}_n\cap\LL(\SSS,\QQ)$. Also, let $\QQ_n$ be a $K_3$-qi section in $\LL(\SSS,\QQ)$ through $a_n$ (here, $\QQ_n$ need not lie completely in $X$.)

 Now, $\widehat{\phi_n}(\mu'_n)$ is an $L$-quasigeodesic in $\widehat{F}_{n-1}$, and $d_{\widehat{F}_{n-1}}(\widehat{F}_{n-1}\cap\QQ_n,\widehat{\phi}(a_n))\leq\eta_{\ref{new mbdl}}(2K_3+1)$. Similarly, $d_{\widehat{F}_{n-1}}(\widehat{F}_{n-1}\cap\QQ,\widehat{\phi}(x))\leq\eta_{\ref{new mbdl}}(2K_2+1)$. The subsegment of $\lambda^b_{n-1}$ joining  $\widehat{F}_{n-1}\cap\QQ_n,\widehat{F}_{n-1}\cap\QQ$ and the quasigeodesic $\widehat{\phi}(\mu_n)$ track each other outside horosphere-like sets by BPP and Lemma \ref{extend quasigeodesic}. So, there exists $y'_{n-1}\in\mu'_n\cap X$ such that $d_X(y_{n-1},\widehat{\phi}_n(y'_{n-1}))\leq C$. But $d_{\widehat{F}_n}(y'_{n-1},x)\leq 1$. In the second possibility, we in fact have that $\mu'_n\cap X=x$, therefore, $y'_{n-1}=x$. Therefore, $d_X(x,y_{n-1})\leq C+2$. (If $\QQ=\SSS$, we replace $\LL(\SSS,\QQ)$ by $\LL(\QQ,\SSS')$).\\
We proceed inductively. Suppose we have $y_k\in F_k\cap\QQ$ such that $d_X(y_k,x)\leq (n-k+1)(C+2)$ for some $0< k\leq n-1$.

As before, let $a_k\in\lambda_k\cap\LL(\SSS,\QQ)$ such that $d_{F_k}(y_k,a_k)\leq 1$. Let $\mu'_k$ be the subsegment of $\lambda_k$ joining $a_k,y_k$. Note that if $a_k$ is a cone point, then $\mu'_k$ is an edge of length $\frac{1}{2}$ joining $y_k$ and $a_k$, and if not $\mu'_k\subset\lambda^b_k$, Also, let $\QQ_k$ be $K_3$-qi sections in $\LL(\SSS,\QQ)$ through $a_k$. Now, $\widehat{\phi}_k(\mu'_k)$ is an $L$-quasigeodesic in $\widehat{F}_{k-1}$, and $d_{\widehat{F}_{k-1}}(\widehat{F}_{k-1}\cap\QQ_n,\widehat{\phi}(a_k))\leq\eta_{\ref{new mbdl}}(2K_3+1)$ and $d_{\widehat{F}_{k-1}}(\widehat{F}_{k-1}\cap\QQ,\widehat{\phi}(y_k))\leq\eta_{\ref{new mbdl}}(2K_2+1)$. Applying Lemma \ref{extend quasigeodesic}, the subsegment of $\lambda^b_{k-1}$ joining  $\widehat{F}_{k-1}\cap\QQ_n,\widehat{F}_{k-1}\cap\QQ$ and $\widehat{\phi}(\mu_k)$ track each other outside horosphere-like sets and thus, there exists $y'_{k-1}\in\mu'_k\cap X$ such that $d_X(y_{k-1},\widehat{\phi}(y'_{k-1}))\leq C$. Thus, $d_X(y_{k-1},y_k)\leq C+2$, which further implies that $d_X(y_{k-1},x)\leq (n-k)(C+2)$. This gives us $y=y_0\in F_0\cap\QQ$ such that $d_X(y,x)\leq d_B(b_0,b)(C+2)$.

\medskip
\noindent
{\bf Case 2:} $x\notin\QQ$ for any qi section $\QQ$ in the ladder decomposition. Suppose $x\in\LL(\QQ',\QQ'')$ where $\QQ',\QQ''$ are the consecutive qi sections in the ladder decomposition, i.e., if $\QQ'=\SSS_j$ (resp. $\SSS'_j)$, then $\QQ''=\SSS'_j$ or $\SSS_{j+1}$ depending on the type of the ladder (resp. $\SSS_{j+1}$). For all $0\leq i\leq n$, let $p_i=\QQ'\cap F_i$ and $q_i=\QQ''\cap F_i$. Also let $\widehat{\mu}_i$ be the subsegment of $\widehat{\lambda}_i$ joining $p_i,q_i$. We have $d_{\widehat{F}_{n-1}}(\widehat{\phi}_n(p_n),p_{n-1})\leq \eta_{\ref{new mbdl}}(2K_2+1)$ and $d_{\widehat{F}_{n-1}}(\widehat{\phi}_n(q_n),q_{n-1})\leq \eta_{\ref{new mbdl}}(2K_2+1)$. Applying \lemref{extend quasigeodesic} and BPP, the paths $\widehat{\mu}_{n-1}$ and $\widehat{\phi}_n(\widehat{\mu}_n)$ track each other  outside horosphere-like sets. Therefore, there exists $C>0$, depending on $L$, and $x_{n-1}\in\lambda^b_{n-1}$ such that $d_{F_{n-1}}(\widehat{\phi}(x),x_{n-1})\leq C$. Thus, $d_X(x,x_{n-1})\leq C+1$.  \\
As in the above case, by induction we can show that for $0\leq k\leq n-1$, there exists $x_k\in\lambda^b_k\cap\LL(\QQ',\QQ'')$ such that $d_X(x,x_k)\leq (n-k)(C+1)$. In particular, there exists $y\in\lambda^b\cap\LL(\QQ',\QQ'')$ such that $d_X(x,x_k)\leq d_B(b,b_0)(C+1)$. We have shown the following result.
\begin{lemma}\label{vertical qg}
There exists $K_{\ref{vertical qg}}\geq 1$ such that for any $x\in\LL\cap X$, there exists a vertical quasigeodesic ray $r_x:V([\pi(x),b_0])\to X$ such that $r_x(\pi(x))=x$ and $d_B(\pi(x),v)\leq d_X(r_x(\pi(x)),r_x(v))\leq K_{\ref{vertical qg}}d_B(\pi(x),v)$.
\end{lemma}
\noindent
Now, we can prove Mitra's criteria. Let $y_0\in F_{b_0}$. Let us denote $c(y,y')$ by $c$ and $\bar{c}(y,y')$ by $\bar{c}$. We have the following from \cite{krishna-sardar}.
\begin{lemma}\label{criterion satisfied}\cite[Lemma 5.21]{krishna-sardar}
	Given $D>0$, there is $D_{\ref{criterion satisfied}}>0$ such that the following holds.		
	If $d_{\XH}(y_0,c)\leq D$ then $d_{\YH}(y_0,\bar{c})\leq D_{\ref{criterion satisfied}}$.
\end{lemma}

\noindent
Let $k_0$ be the quasiconvexity constant of $A$ in $B$. The proof is very similar to that of \cite[Lemma 5.21]{krishna-sardar}.
\begin{lemma}\label{criterion satisfied-2}
	Given $D>0$, there is $D_{\ref{criterion satisfied-2}}=D_{\ref{criterion satisfied-2}}(D)>0$ such that the following holds.	
	If $d_{X}(y_0,c\cap X)\leq D$ then $d_{Y}(y_0,\bar{c}\cap Y)\leq D_{\ref{criterion satisfied-2}}$.
\end{lemma}
\begin{proof}
	Let $x\in c\cap X$ such that $d_X(y_0,x)\leq D$. Then we have, $d_{\XH}(y_0,x)\leq D$. 
	
	\noindent	
	Recall that $c$ is a concatenation of lift of geodesic segments of $B$ in  $K_2$-qi section $\SSS_k$ or $\SSS'_k$ and fiber geodesics. Let $\QQ,\QQ'$ be consecutive $K_2$-qi sections with $d_h(\QQ,\QQ')\leq R_0$ and $x\in\LL(\QQ,\QQ')$. Then $c\cap\LL(\QQ,\QQ')$ is either a lift of a geodesic segment of $B$ in $\QQ$, joining the points $z$ to $w$ or it is the concatenation of this lift and a fiber geodesic $\sigma$ in $\widehat{F}_{\pi(w)}$ joining $\QQ$ to $\QQ'$. Let $z'=\QQ'\cap\sigma$ and $b'=\pi(z')$. Note that $z$ and $z'$ are one of the $y_k$'s or $y'_k$'s. 
	
	So we have two cases. First, suppose $x\in\sigma\cap X$. Let $b$ be a nearest point projection of $\pi(x)$ on $A$. Then, $d_B(\pi(x),b)\leq d_B(\pi(x),b_0)\leq D$. By the definition of $\bar{c}$, $\LL(\QQ,\QQ')\cap F_{b}\subset\bar{c}$. By Lemma \ref{vertical qg}, there exists $y\in\LL(\QQ,\QQ')\cap F_b$ such that $d_X(x,y)\leq K_{\ref{vertical qg}}d_B(\pi(x),\pi(y))\leq K_{\ref{vertical qg}}D$. 
	
	Now suppose $x$ lies in the lift of the geodesic segment $[\pi(z),\pi(w)]_B$ in $\QQ$. So, $\pi(x)\in[\pi(z),\pi(w)]_B$. Let $\overline{\pi(z)},\overline{\pi(w)}$ denote the nearest point projections of $\pi(z),\pi(w)$ respectively on $A$. By \lemref{trivial lemma} and stability of quasigeodesics, 
	$d_B(\pi(x), [\overline{\pi(z)},\overline{\pi(w)}]_A)\leq D_{\ref{trivial lemma}}(D,k_0,\delta)+D_{\ref{stab-qg}}(\delta,k)$. Let $b\in[\overline{\pi(z)},\overline{\pi(w)}]_A$ such that $d_B(\pi(x),b)\leq D_{\ref{trivial lemma}}(D,k_0,\delta)+D_{\ref{stab-qg}}(\delta,k)$. 	
	By the definition of $\bar{c}$, the lift of the geodesic segment $[\overline{\pi(z)},\overline{\pi(w)}]_A$ in $\QQ$ is a subsegment of $\bar{c}$. Let $y=\QQ\cap F_b$. Again by Lemma \ref{vertical qg}, $d_X(x,y)\leq K_{\ref{vertical qg}}d_B(\pi(x),b)\leq K_{\ref{vertical qg}}d_B(\pi(x),b)\leq K_{\ref{vertical qg}}(D_{\ref{trivial lemma}}(D,k_0,\delta)+D_{\ref{stab-qg}}(\delta,k))$.
	
	\noindent
	Let $D_1=\max\{K_{\ref{vertical qg}}D,K_{\ref{vertical qg}}(D_{\ref{trivial lemma}}(D,k_0,\delta)+D_{\ref{stab-qg}}(\delta,k))\}$. Then, $d_X(y_0,y)\leq d_X(y_0,x)+d_X(x,y)\leq D+D_1$. 	
	By \lemref{prop-embed}, there exists a function $\eta_0:\NN\to\NN$ such that, $d_Y(y_0,y)\leq\eta_0(D+D_1)$. We take $D_{\ref{criterion satisfied-2}}(D)=\eta_0(D+D_1)$.\end{proof}


\section{Examples and Applications}\label{section 6}
\subsection{Short exact sequence of groups}
Let $$1\to (K,K_1)\to (G,N_G(K_1))\stackrel{\pi}{\to} (Q,Q_1)\to1$$ be a short exact sequence of pairs of finitely generated groups. Fixing generating sets $S$ for $G$, $S_K$ for $K$ such that $S_K\subset S$ and $S_Q=\pi(S)\setminus\{1\}$, there is a natural metric graph bundle associated to $1\to K\to G\to Q\to 1$, with the base space as the Cayley graph $\Gamma_Q$ of $Q$ and the fibers as the copies of Cayley graph $\Gamma_K$ of $K$ with respect to $S_K$ (see \cite[Example 3.3.1]{krishna-sardar}), with the total space $\Gamma_G$, Cayley graph of $G$ with respect to $S$. 

Recall that $G$ is said to {\em preserve cusps} if for every $g\in G$, there exists $k_g\in K$ such that $gK_1g^{-1}=k_gK_1k_g^{-1}$. In \cite{Pal}, Pal showed existence of CT map for relatively hyperbolic extension of groups.
\begin{theorem}\cite[Theorem 3.11]{Pal}\label{pal ct}
Suppose we have a short exact sequence of pairs of finitely generated groups $$1\to (K,K_1)\to (G,N_G(K_1))\stackrel{\pi}{\to} (Q,Q_1)\to1$$ such that $K$ is strongly hyperbolic relative to a non-trivial proper subgroup $K_1$ and $G$ preserves cusps. Moreover, if $G$ is weakly hyperbolic relative to $K_1$ and strongly hyperbolic relative to $N_G(K_1)$, then the embedding $i:\Gamma_K\to\Gamma_G$ admits a CT map. 	
\end{theorem}

\noindent
Under the first three conditions in \thmref{pal ct}, i.e., that $K$ is a strongly hyperbolic relative to $K_1$, $G$ preserves cusps and $G$ is weakly hyperbolic relative to $K_1$, Pal (cf. \cite[Corollary 2.11]{Pal}) showed that there exists a qi section $s:Q\to G$ and that $Q$ is a hyperbolic metric space.\\
Now, for $B=\Gamma_Q$, since $(K,K_1)$ is strongly relatively hyperbolic, we can associate a relatively hyperbolic metric graph bundle $\pi:X\to B$ to this short exact sequence (cf. \cite[Subsection 5.1]{pranab-mahan}). 

\noindent
As an application of the combination theorem part of the main result, we have the following.
\begin{theorem}\label{combi for short exact seqn}
Let $1\to (K,K_1)\to (G,N_G(K_1))\to (Q,Q_1)\to 1$ be a short exact sequence of pairs of finitely generated groups such that $K$ is strongly hyperbolic relative to a proper non-trivial subgroup $K_1$ and suppose $G$ preserves cusps. Let $G$ be weakly hyperbolic relative to $K_1$ and a cone-bounded strictly flare condition be satisfied by $\EE(G,K_1)$. Then $G$ is strongly relatively hyperbolic relative to $N_G(K_1)$.	 
\end{theorem}
Let $Q'<Q$ be a finitely generated subgroup and $G_1=\pi^{-1}(Q')$. Let $S_1$ be a generating set of $G_1$ such that $S_K\subset S_1\subset S$. Let $S_{Q'}=\pi(S_1)\setminus\{1\}$. Clearly $\Gamma_{Q'}$ (with respect to $S_{Q'}$) is a subgraph of $\Gamma_Q$ and $\Gamma_{G_1}=\pi^{-1}(\Gamma_{Q'})$. By \corref{cor: restriction bundle} it follows that $\Gamma_{G_1}$ is the pullback of $\Gamma_G$ under $\Gamma_{Q'}\hookrightarrow \Gamma_Q$.

If $Q'$ is qi embedded in $Q$, then $\Gamma_{Q'}$ is also a hyperbolic metric space and we can associate a relatively hyperbolic metric graph bundle structure to $\pi|_{\Gamma_{G_1}}:\Gamma_{G_1}\to\Gamma_{Q'}$. By \propref{coned-off pullback}, $\EE(G_1,K_1)$ is a pullback of $\EE(G,K_1)$ under $\Gamma_{Q'}\hookrightarrow \Gamma_Q$.

\noindent
Further, we have the following corollary to \thmref{main theorem}.
\begin{theorem}\label{main application}
	Suppose $1\to (K,K_1)\to (G,N_G(K_1))\to (Q,Q_1)\to 1$ is a short exact sequence of relatively hyperbolic groups
	where $K$ is strongly hyperbolic relative to a proper non-trivial subgroup $K_1$ and $G$ preserves cusps. Suppose $G$ is weakly hyperbolic relative to $K_1$ and it is strongly hyperbolic relative to $N_G(K_1)$. Suppose $Q'$ is a qi embedded subgroup of $Q$ and $G_1=\pi^{-1}(Q')$.
	Then $(G_1,N_{G_1}(K_1))$ is strongly relatively hyperbolic and $(G_1,N_{G_1}(K_1))\to(G,N_{G}(K_1))$ admits CT.
\end{theorem}

\noindent
For the rest of the section, we have the following set-up. 
\begin{enumerate}[leftmargin=\parindent,align=left,labelwidth=\parindent,labelsep=4pt]
	\item  $\pi:X\rightarrow B$ be an $\eta$-metric graph bundle, with $B$, a $\delta_0$-hyperbolic metric graph.
	\item For each $b\in V(B)$, the fiber $F_b$ strongly hyperbolic relative to $\HH_b=\{H_{b,\alpha}\}_{\alpha\in \Lambda}$ and the coned-off fiber $\widehat{F}_b=\mathcal{E}(F_b,\HH_b)$ is $\delta_0$-hyperbolic.
	\item $X$ is strongly hyperbolic relative to the collection $\CC=\{\CC_{\alpha}\}_{\alpha\in\Lambda}$, and for each $\alpha\in\Lambda$, the connected component of the cone locus is $\BB_\alpha$.
	\item The spaces $\XH$, $\EE(X,\CC)$ and $X^h:=\GG(X,\CC)$ are $\delta$-hyperbolic for $\delta>0$.	
\end{enumerate} 
\noindent
We also fix some conventions and notations that we follow for the rest of this section.

\noindent
\begin{itemize}[leftmargin=\parindent,align=left,labelwidth=\parindent,labelsep=4pt]
\item Through any point of $X$ there is a $K_0$-qi section over $B$.
\item The pullback map is $i_{Y,X}:Y\to X$ and induces $i_{\YH,\XH}:\YH\to\XH$ and the corresponding CT maps are $\partial i_{Y,X}:\partial Y^h\to \partial X^h$ and $\partial i_{\YH,\XH}:\partial\YH\to\partial\XH$. 
\item For every $x,y\in X$, the path $c(x,y)\subset\XH$ from \secref{section 5} is a $K_{\ref{section 5}}$-quasigeodesic. Similarly, for $x,y\in Y$, the modified path $\tilde{c}(x,y)\subset\YH$ is also a $K_{\ref{section 5}}$-quasigeodesic.\\
For the rest of this section, we fix $b_0\in V(A)$ and $F:=F_{b_0}$ and $p\in F$. 
\item The corresponding CT maps are $\partial i:\partial F^h\to\partial X^h$, $\partial i_{\XH}:\partial\widehat{F}\to\partial\XH$, $\partial i_Y:\partial F^h\to\partial Y^h$ and $\partial i_{\YH}:\partial\widehat{F}\to\partial\YH$.
\end{itemize} 
\cite[Proposition 6.6]{krishna-sardar} gives a set-theoretic description of the boundary of a hyperbolic metric (graph) bundle. Note that $X^h$ no longer has a metric (graph) bundle structure. But we have a set-theoretic description of $\partial\XH$. Let $\xi\in\partial B$ and let $\alpha$ be a quasigeodesic ray in $B$ with $\alpha(\infty)=\xi$. Let $\partial^{\xi}\XH=\{\tilde{\alpha}(\infty)\mid \tilde{\alpha}\text{ is a lift of }\alpha\text{ in }\XH\}$. By \cite[Lemma 6.4]{krishna-sardar}, $\partial^{\xi}\XH$ is independent of the choice of $\alpha$.
\begin{prop}\label{bundle boundary}
$\displaystyle{\partial \XH=\Lambda\widehat{F}\cup\bigcup_{\xi\in\partial B}\partial^{\xi}\XH}$.  
\end{prop}
\noindent
Now we prove a variation of \cite[Lemma 6.9]{krishna-sardar}. The proof is similar.
\begin{lemma}\label{lemma: ct to base}
	Suppose $b\in B$ and $\alpha_n$ is a sequence of geodesics in $B$ joining $b$ to $b_n$. Let $\tilde{\alpha}_n$ is a uniform qi lift of $\alpha_n$ for all $n$. Suppose $\tilde{\alpha}_n$ joins $\tilde{\alpha}_n(b)=x^0_n$ to $\tilde{\alpha}_n(b_n)=x_n$ and suppose the set $\{x^0_n\}$ has finite diameter. If $x_n\to z\in \partial\XH$, then $\lim_{n\to\infty}b_n$
	exists. If $\xi=\lim_{n\to\infty}b_n$ and $\alpha:[0,\infty)\to B$ is a geodesic ray joining $b$ to $\xi$, then there is a uniform qi lift 
	$\tilde{\alpha}$ of $\alpha$ such that $\tilde{\alpha}(\infty)=z$.
\end{lemma}
\begin{proof}
Since $x_n\to z$, there exists $D>0$ such that for all $M>0$ there is $N=N(M)>0$ with $Hd_{\XH}(\tilde{\alpha}_m|_{[0,M]}, \tilde{\alpha}_n|_{[0,M]})\leq D$ for all 
$m,n\geq N$ by \lemref{lem: bdry defn}(1). Then, for all $M>0$,
$Hd_{B}(\alpha_m|_{[0,M]}, \alpha_n|_{[0,M]})\leq D$ for all $m,n\geq N$. Again by \lemref{lem: bdry defn}(1), the sequence $\{b_n\}$ converges to a point of $\xi\in\partial B$. Let $\alpha$ be a geodesic ray in $B$ joining $b$ to $\xi$.	We claim that $z\in \partial^{\xi} X$.

Given $t\geq 0$, by \lemref{lem: bdry defn}(2) there exists $N'=N'(t)$ such that
$d(\alpha(t),\alpha_n)\leq D'$ for all $n\geq N'$, for some $D'$ depending on $\delta$.
Let $N_0=\max\{N(t),N'(t)\}$. Let $t'>0$ be such that $d_X(\alpha(t),\alpha_{N_0}(t'))\leq D'$.
Define $\tilde{\alpha}(t)=\hat{\phi}_{uv}(\tilde{\alpha}_{N_0}(t'))$ where $u=\alpha_{N_0}(t'),v=\alpha(t)$
and $\phi_{uv}$ is a fiber identification map. We can easily see that this defines a qi section over $\alpha$ and $z=\tilde{\alpha}(\infty)$.
\end{proof}
\begin{lemma}\label{qi lift in X}
Let $z\in\bigsqcup_{\xi\in\partial B}\partial^{\xi}\XH$. Then there exists a qi lift $\tilde{\alpha}$ in $\XH$ such that $\tilde{\alpha}\subset X$ and $\tilde{\alpha}(\infty)=z$ in $\partial\XH$.
\end{lemma}
\begin{proof}
Suppose $z\in\partial^{\xi}\XH$ for some $\xi\in\partial B$. Let $\beta$ be geodesic ray $B$ with $\beta(\infty)=\xi$ and let $\tilde{\beta}$ be a qi lift with $\tilde{\beta}(
\infty)=z$. We define a path $\tilde{\alpha}$ as follows. For each $n\in\ZZ_{\geq0}$, $\tilde{\alpha}(n)=\tilde{\beta}(n)$ if $\tilde{\beta}(n)\in X$. If $\tilde{\beta}(n)\in\nu(H_{b,\alpha})$ for some $H_{b,\alpha}\in\HH_b$ for some $b\in V(B)$, then define $\tilde{\alpha}(n)=h$ for some $h\in H_{b,\alpha}$. This is a discrete path satisfying $d_{\XH}(\tilde{\alpha},\tilde{\beta})\leq\frac{1}{2}$. Moreover, it is coarsely unique. If for any $n\in\ZZ_{\geq0}$, $\tilde{\beta}(n)\in\nu(H_{b,\alpha})$ and $\tilde{\alpha}(n)=h$ and if we have another path $\tilde{\alpha'}$ defined as above with $\tilde{\alpha'}(n)=h'$ for some $h'\in H_{b,\alpha}$, then clearly, $d_{\XH}(\tilde{\alpha},\tilde{\alpha'})\leq 1$. Next we show that $\tilde{\alpha}$ is a quasigeodesic ray. Let $m>n\geq 0$. Then $l_{\XH}(\tilde{\alpha}([n,m]))=\sum_{i=n}^{m-1}d_{\XH}(\tilde{\alpha}(i),\tilde{\alpha}(i+1))\leq\sum_{i=n}^{m-1}d_{\XH}((\tilde{\beta}(i),\tilde{\beta}(i+1))+1)=2(m-n)\leq 2d_{\XH}(\tilde{\alpha}(n),\tilde{\alpha}(m))$. Moreover, clearly, $\pi(\tilde{\alpha})=\beta$. Therefore, $\tilde{\alpha}$ is a qi lift with $\tilde{\alpha}(\infty)=z$.
\end{proof}
\subsection{Cannon-Thurston lamination} \label{defn: ct lamination}
\noindent
\noindent
In \cite{krishna-sardar}, the authors generalize the definition of ladder to ladder with infinite girth. We recall it here. 
\begin{definition}\label{ladder defn}
	Let $\SSS$ be a $K_0$-qi section over $B$ in $X$, $z,z'\in \partial F_b$ such that $z\neq z'$. For each $s\in B$, let
 $z_s=\partial\phi^h_{bs}(z)$ and $z'_s=\partial \phi^h_{bs}(z')$, where $\phi_{bs}:F_b\to F_s$ be the fiber identification map.\\
 A {\bf semi-infinite ladder}, $\LL(\SSS,z)$, is the union of quasigeodesic rays in $\widehat{F}_s,s\in B$ joining $\SSS\cap \widehat{F}_s$ to $z_s$ and a {\bf bi-infinite ladder}, $\LL(z;z')$ is the union of quasigeodesic lines in $\widehat{F}_s,s\in B$ joining $z_s$ and $z'_s$.
 Both the ladders are contained in $\XH$. \end{definition}
\noindent
Note that the ladders are coarsely well-defined by Lemma \ref{ideal triangles are slim}.
\begin{lemma}\label{infinite ladder property}\cite[Lemma 6.14]{krishna-sardar}
	Suppose $\LL$ is an infinite girth ladder. 
	\begin{enumerate}[leftmargin=\parindent,align=left,labelwidth=\parindent,labelsep=4pt]
		\item There is a uniformly coarsely Lipschitz retraction $\Pi:\XH\map \LL$ such that for all $b\in B$ and $x\in \widehat{F}_b$, $\Pi(x)$ is an electric projection of $x$ in $\widehat{F}_b$ on $\LL\cap \widehat{F}_b$.		
		\item Infinite girth ladders are uniformly quasiconvex and there exists $C>0$ such that a $C$-neighbourhood of the ladder is qi embedded in $\XH$.
		\item Through any point of $\LL$, there exists a $\bar{K}_0$-qi section contained in $\LL$.
		\item Any qi section in $\LL$ coarsely bisects it into two subladders.
		\end{enumerate}
\end{lemma}
\begin{definition}
$(1)$ Let $Z_1\subset Z_2$ be hyperbolic geodesic metric spaces such that the inclusion map $I:Z_1\to Z_2$ admits CT map, $\partial I:\partial Z_1\to\partial Z_2$. Then the {\bf Cannon-Thurston lamination} is given by $\Lam{Z_2}(Z_1)=\{(a,b)\in \partial^{(2)}Z_1:\partial I(a)=\partial I(b)\}$.\\
$(2)$ A (quasi)geodesic line $\alpha$ in $Z_1$ with $(\alpha(-\infty),\alpha(\infty))\in\Lam{Z_2}(Z_1)$ is called a {\bf leaf} of the CT lamination.\\
$(4)$ Let $Z\subset Z_1$ be a qi embedded subset. A leaf $\alpha\subset Z_1$ of $\Lam{Z_2}(Z_1)$ is said to be {\bf carried by $Z$} if $\alpha(\pm\infty)\in\partial Z$.
\end{definition}
\noindent
With the above definition, in our case, we have $\Lam{X}(F):=\Lam{X^h}(F^h)=\{(z,z')\in \partial^{(2)}F^h:\partial i(z)=\partial i(z')\}$. 
\begin{lemma}\label{cor: asymptotic lifts}\cite[Lemma 6.19]{krishna-sardar}.
	Let $\LL=\LL(z_1;z_2)$ for $(z_1,z_2)\in\partial^{(2)}_{\XH}(\widehat{F})$.
	Let $\xi\in\partial B$ and $\beta:[0,\infty)\to B$ be a continuous, arc length parameterized $\kappa_0$-quasigeodesic in $B$. Then any two qi lifts of $\beta$ in $\LL$ are asymptotic in $\XH$.
\end{lemma}

\noindent
First we recall the part of \cite[Theorem 6.25]{krishna-sardar} relevant to us in our induced coned-off metric graph bundle setting. 
\begin{theorem}\label{CT leaf in Y-0}
	Suppose we have a metric graph bundle satisfying the hypotheses of Theorem \ref{main theorem}. Suppose $\gamma$ is a quasigeodesic line in $\YH$ 
	such that $\gamma(\infty),\gamma(-\infty)$ are identified in $\partial\XH$. Then $\gamma(\pm\infty)\in\Lambda_{\YH}(\widehat{F}_b)\setminus\bigcup_{\xi\in\partial A}\partial^{\xi}\YH$ for some (any) $b\in V(A)$. 
\end{theorem}
\noindent 
Now we prove a weak analogue of \cite[Theorem 6.25]{krishna-sardar}. 
\begin{theorem}\label{CT leaf in Y}
	Suppose we have a metric graph bundle satisfying the hypotheses of \thmref{main theorem} such that $X$ is a proper metric space. Let $\gamma$ be a geodesic line in $Y^h$ with $\partial i_{Y,X}(\gamma(\infty))=\partial i_{Y,X}(\gamma(-\infty))$ in $\partial X^h$. Let $b\in V(A)$ and $F=F_b$. Then\\
	$(1)$ $\gamma(\pm\infty)\in\partial i_{F,Y}(\partial F^h)$.\\	
	$(2)$ 
	There exists a uniform quasigeodesic line $\sigma$ in $F^h$ at a finite Hausdorff distance from $\gamma$, in $Y^h$, satisfying $\partial i_{F,Y}(\sigma(\pm\infty))=\gamma(\pm\infty)$.
\end{theorem}
\begin{proof}
Throughout the proof, by Lemmas \ref{hamenstadt} and \ref{all boundaries}, we identify each non-parabolic point of $\partial X^h$ (resp. $\partial Y^h$) with the corresponding point in $\partial\XH$ (resp. $\partial\YH$).
We first note that both $\gamma(\pm\infty)$ cannot be parabolic limit points in $\partial Y^h$, since each fiber is embedded in $Y$ and $X$ by a uniform type preserving map. \\
{\bf Case 1.} Suppose $\partial i_{Y,X}(\gamma(\pm\infty))$ is non-parabolic in $\partial X^h$. Then both $\gamma(\pm\infty)$ are non-parabolic limit points. Identifying them with the boundary points in $\partial\YH$, and applying Theorem \ref{CT leaf in Y-0}, we have that $\gamma(\pm\infty)\in\Lambda_{\YH}(\widehat{F})$, and again by Lemma \ref{hamenstadt}, $\gamma(\pm\infty)\in\Lambda_{Y^h}(F^h)$. As $F^h$ is proper, by Lemma \ref{CT-limset}, $\gamma(\pm\infty)\in\partial i_{F,Y}(\partial F^h)$.

\medskip
\noindent
We prove $(2)$ in this case. This part is very similar to that in \cite[Theorem 6.25]{krishna-sardar} Let $\sigma:\RR\to F^h$ be a geodesic line such that $\partial i_{b,Y}(\sigma(\pm\infty))=\gamma(\pm\infty)$. 
Let $\{a_n^{\pm}\}$ be a sequence of points in $\sigma((0,\pm\infty))\cap F$ such that $a_n^{\pm}\to\sigma(\pm\infty)$ in $F^h$. Suppose $a_n^{\pm}\to z_{\pm}$ in $\widehat{F}$.
Let $\LL=\LL(z_+;z_-)$ and $\Sigma_{\pm n}$ be a $\bar{K}_0$-qi section in $\LL$ passing through $a_n^{\pm}$, $n\in \NN$. Let $\alpha$ be a geodesic ray in $B$ with $\alpha(0)=b$ and $\alpha(\infty)=\xi$. Let $\tilde{\alpha}_{\pm n}$ be a qi lift of $\alpha$ in $\SSS_{\pm n}$. By \corref{cor: asymptotic lifts}, $\tilde{\alpha}_{\pm n}$'s are asymptotic. Also, these are $\bar{K}_0$-quasigeodesics. Recall that by \lemref{ideal triangles are slim}, the ideal triangle $\tilde{\alpha}_n\cup\tilde{\alpha}_{-n}\cup[a^+_n,a^-_n]_{\XH}$ is $R=D_{\ref{ideal triangles are slim} }(\delta,\bar{K}_0)$-slim. So, there exists $t_n\geq 0$ such that $\tilde{\alpha}_n(t_n)\in N_R(\tilde{\alpha}_{-n})$ and $\tilde{\alpha}_{-n}(t_n)\in N_R(\tilde{\alpha}_n)$ (as long as $\tilde{\alpha}_n(t_n)$ and
$\tilde{\alpha}_{-n}(t_n)$ are not contained in the $R$-neighborhood of $[a^+_n,a^-_n]_{\XH}$). This further implies that $d_{\XH}(\tilde{\alpha}_n(t_n),\tilde{\alpha}_{-n}(t_n))\leq R+\bar{K}_0R+\bar{K}_0$.\\
Let $R_2=\max\{R+\bar{K}_0R+\bar{K}_0,M_{\bar{K}_0}\}$.
Thus for all $n\in\NN$, $U_{n}=U_{R_2}(\SSS_n, \SSS_{-n})\neq \emptyset$. Let $\LL_n=\LL(\SSS_n,\SSS_{-n})$. Let $b_n$ be a nearest point projection of $b$ on $U_n$ and let $b'_n$ be a nearest point projection of $b_n$ on $A$. Let $\beta_n=[b,b'_n]$.
Then $\bar{c}_n=\bar{c}(a^+_n,a^-_n)$ is formed by the concatenation of the lifts $\tilde{\beta}_{\pm n}$ in $\SSS_{\pm n}$ and $\widehat{F}_{b'_n}\cap\LL_n$. This is a uniform quasigeodesic in $\YH$. Let $p\in F$. Let $\eta_n$ be hyperbolic geodesics in $Y^h$ joining $a^+_n,a^-_n$.

By \lemref{convergence explained}, there exists $M>0$ such that $d_{Y^h}(p,\tilde{\eta}_n)\leq M$ for all $n>0$. For all $n>0$, by \lemref{lemma from limint}, $d_Y(p,\tilde{\eta}_n\cap Y)\leq K_{\ref{lemma from limint}}(M)$ and by \lemref{hyp-pel-track}, $d_Y(p,\tilde{c}_n\cap Y)\leq K_{\ref{lemma from limint}}(M)+C$. Let $D=K_{\ref{lemma from limint}}(M)+C$. Then, $d_{\YH}(p,\tilde{c}_n\cap Y)\leq D$ for all $n>0$. But as $n\to\infty$, $d_{\YH}(p,a^{\pm}_n)\to\infty$ implies that $d_{\YH}(p,\tilde{\beta}_{\pm n})\to\infty$. So,  \begin{equation*}
	d_{\YH}(p,F_{b'_n}\cap\LL_n)\leq D \text{ for large enough }n.
\end{equation*}
For $n>0$, let $\sigma_n=\widehat{F}\cap\LL_n$. We first show that $\sigma_n$ is a uniform quasigeodesic in $\YH$.\\ 
Since $\sigma_n$ is a geodesic in $\widehat{F}$, it is uniformly coarsely Lipschitz and properly embedded in $\YH$. By \propref{ladders are qi embedded}(2), $Hd_{\YH}(\sigma_n,\bar{c}_n)$ is uniformly bounded. So by \lemref{quasigeod criteria}, $\sigma_n$ is a uniform quasigeodesic in $\YH$. Let $\sigma^h_n$ be a hyperbolic geodesic in $Y^h$ with endpoints $a^+_n,a^-_n$. Then $Hd_Y(\sigma_n\cap Y,\sigma^h_n\cap Y)$ is uniformly small and $\sigma_n\cap F$ lies in a uniform neighbourhood of $\sigma\cap F$ in $F$. So, $\sigma^h_n$ lies in a uniform neighbourhood of $\sigma$ in $Y^h$. As $a^{\pm}_n\to\infty$, since $Y^h$ is proper, $\sigma^h_n$ converges to a hyperbolic geodesic line in $Y^h$, which is uniformly Hausdorff close to $\sigma$. Again by \lemref{quasigeod criteria}, $\sigma$ is a quasigeodesic line in $Y^h$.

\noindent
{\bf Case 2. }Now suppose $\partial i_{Y,X}(\gamma(\pm\infty))$ is a parabolic limit point in $\partial X^h$. Let $\partial i_{Y,X}(\gamma(\pm\infty))=\Lambda(\CC_{\alpha}^h)$, where $\CC_{\alpha}\in\CC$. Denote by $\CC_{Y,\alpha}=\CC_{\alpha}\cap Y$.\\
{\bf Case 2.1. } Suppose $\gamma(\infty)$ is a non-parabolic limit point of $\partial Y^h$ and $\gamma(-\infty)$ is a parabolic limit point. Let $\gamma(\infty)\in\partial^{\xi}{\YH}$ for some $\xi\in\partial A$. Let $\beta_1:[0,\infty)\to A$ be a geodesic ray with $\beta_1(0)=b,\beta_1(\infty)=\xi$ and $\tilde{\beta}_1$ be a qi lift in $\YH$ such that $\tilde{\beta}_1(\infty)=\gamma(\infty)$ in $\partial^{\xi}{\YH}$. Let $\tilde{\beta}_2\subset\CC_{Y,\alpha}$ be a qi lift of $\beta$ in $\YH$. Then $\partial i_{\YH,\XH}(\tilde{\beta}_1(\infty)) =\partial i_{\YH,\XH}(\tilde{\beta}_2(\infty))\in\Lambda\CC_{\alpha}$ and this implies that $Hd_{\XH}(\tilde{\beta}_1,\tilde{\beta}_2)<\infty$ and as $\YH$ is properly embedded in $\XH$, $Hd_{\YH}(\tilde{\beta}_1,\tilde{\beta}_2)<\infty$. This gives a contradiction. Therefore, $\gamma(\infty)\in\partial i_{F,Y}({\partial{F^h}})$.\\
Let $z_1,z_2\in\partial F^h$ such that $\partial i_{F,Y}(z_1)=\gamma(\infty)$ and $\partial i_{\widehat{F},\YH}(z_2)=\gamma(-\infty)$. We know that for $z_2=\Lambda H_{b,\alpha}^h$, $\partial i(z_2)=\Lambda(\CC_{\alpha}^h)$. We consider this case first. \\
{\bf Case 2.1.1.} $z_2=\Lambda H_{b,\alpha}^h$. Let $\sigma:\RR\to F^h$ be a geodesic line with $\sigma(\infty)=z_1$ and $\sigma(-\infty)=z_2$. Let $a\in H_{b,\alpha}$ be the point where $\sigma$ enters $H_{b,\alpha}^h$. Let $\Sigma$ be a qi section in $\XH$, through $a$, lying in $\CC_{\alpha}$. Without loss of generality, let $\{a_n\}\subset\sigma\cap F$ such that $a_n\to z_1$ in $F^h$. Identifying $z_1\in\partial\widehat{F}$, we know that $\partial i_{\widehat{F},\XH}(z_1)\in\Lambda(\BB_{\alpha})$ in $\partial \XH$ (see Lemma \ref{hamenstadt}, Remark \ref{remark about bdry} and Lemma \ref{all boundaries}). Since $\BB_{\alpha}$ is an isometric embedding, there exists $\xi\in\partial B$ such that $i_{\widehat{F},\XH}(z_1)=\Lambda(\BB_{\alpha})\cap\partial^{\xi}\XH$. Consider $\LL=\LL(\Sigma,z_1)$. Let $\SSS_n$ be a $\bar{K}_0$-qi section in $\LL$ passing through $a_n$. Let $\beta$ be a geodesic ray in $B$ with $\beta(0)=b$ and $\beta(\infty)=\xi$. Let $\tilde{\beta},\tilde{\beta}_n$ be qi lift of $\beta$ in $\SSS,\SSS_{n}$ for each $n$, respectively. By \corref{cor: asymptotic lifts}, $\tilde{\beta}_n$ and $\tilde{\beta}$ are asymptotic and are $\bar{K}_0$-quasigeodesics. Putting $a^-_n=a$ for all $n>0$ and $a^+_n=a_n$ for all $n$, and taking $\tilde{\alpha}_{-n}$ (resp. $\SSS_{-n}$) to be $\tilde{\beta}$ (resp. $\SSS$) for all $n>0$, $\tilde{\alpha}_{+n}$ (resp. $\SSS_{+n}$) to be $\tilde{\beta}_n$ (resp. $\SSS_n$)  we repeat the steps of $(2)$ of Case 1.  
Then, the subsegment of $\sigma$ starting at $a$ and limiting to $z_1$ is a hyperbolic quasigeodesic ray in $Y^h$ and further, by Lemma \ref{new quasigeodesic} and Lemma \ref{quasigeod criteria}, $\sigma$ itself is a quasigeodesic line in $Y^h$. \\
{\bf Case 2.1.2.} $z_2$ is a non-parabolic limit point of $F^h$. Let $z_0=\Lambda H_{b,\alpha}^h$. For $i=1,2$, let $\sigma_i$ be a geodesic line in $F^h$ with $\sigma_i(\infty)=z_i$ and $\sigma_i(-\infty)=z_0$. Applying Case 2.1.1 for the pairs $(z_0,z_1)$ and $(z_0,z_2)$, we have that $\sigma_1$ and $\sigma_2$ are quasigeodesic lines in $Y^h$. Since ideal (quasi)geodesic triangles of $F^h$ (resp. $Y^h$) are slim and by Lemma \ref{quasigeod criteria}, a fiber geodesic line $\sigma$ of $F^h$, with endpoints $z_1,z_2$, is a quasigeodesic line in $Y^h$.\\
{\bf Case 2.2.} Suppose $\gamma(\pm\infty)$ are both non-parabolic limit points of $\partial Y^h$. Let $x_0=\Lambda(\CC_{Y,\alpha}^h)$. Applying Case 2.1 to the pairs $(\gamma(\infty),x_0)$ and $(\gamma(-\infty),x_0)$, we have that $(\gamma(\pm\infty))\in\partial i_{F,Y}(\partial F^h)$. Let $z_1,z_2\in\partial F^h$ such that $\partial i_{F,Y}(z_1)=\gamma(\infty)$ and $\partial i_{F,Y}(z_2)=\gamma(-\infty)$. Then $z_1,z_2$ are non-parabolic limit points of $\partial F^h$ that are mapped to $\Lambda(\CC_{\alpha}^h)$ in $\partial X^h$. Then applying Case 2.1.2, we have that a geodesic line in $F^h$ joining $z_1,z_2$ is a quasigeodesic line in $Y^h$.\end{proof}
\noindent
Let $\xi\in\partial B$. We denote the pullback $\pi^{-1}([b_0,\xi))$ by $Y_{\xi}$ and let $\YH_{\xi}$ denote the induced coned-off metric graph bundle. Let $F:=F_{b_0}$ as we had fixed it before and let $i_{\xi}:F\to Y_{\xi}$ denote the inclusion map. We know that by Theorem \ref{main theorem}, $Y_{\xi}$ is strongly hyperbolic relative to $\CC\cap Y_{\xi}$ and that both $\partial i_{\xi}:\partial F^h\to \partial Y^h_{\xi}$ and $\partial Y^h_{\xi}\to \partial X^h$ exist. Let $\Lam{\xi}(F):=\Lam{Y_{\xi}^h}(F^h)$. Let $(z_1,z_2)\in\Lam{X}(F)$ and $\partial i(z_1)=\partial i(z_2)=z$. We consider the following possibilities.\\
{\bf Case 1:} Both $z_1,z_2$ are non-parabolic and $z$ is non-parabolic in $\partial X^h$. Identifying $z_1,z_2$ with the elements in $\partial\widehat{F}$, we have $(z_1,z_2)\in\Lam{\XH}(\widehat{F})$. Recall from \cite[Lemma 6.17]{krishna-sardar} that
$\Lam{\XH}(\widehat{F})=\sqcup_{\xi\in\partial B}\Lam{\xi,\XH}(\widehat{F})$, where  
$\Lam{\xi,\XH}(\widehat{F})=\{(z_1,z_2)\in\Lam{\XH}(\widehat{F})\mid \partial i_{\XH}(z_1)=\partial i_{\XH}(z_2)\in\partial^{\xi}\XH\}$. Then there exists a unique $\xi\in\partial B$ such that $(z_1,z_2)\in\Lam{\xi,\XH}(\widehat{F})$. Let $\{x_n\},\{y_n\}\subset F$ such that $x_n\to z_1$ and $y_n\to z_2$ in $F^h$. Let $\tilde{\alpha}_n$ and $\alpha^h_n$ denote a partially electrocuted geodesic in $\YH_{\xi}$ and a hyperbolic geodesic in $Y^h_{\xi}$, respectively, joining $x_n,y_n$. Then as $n\to\infty$, $d_{\YH_{\xi}}(p,\tilde{\alpha}_n)\to\infty$. By \lemref{hyp-pel-track} and Lemma \ref{proper embedding in cone space}, we have that $d_{Y^h_{\xi}}(p,\alpha^h_n)\to\infty$. Therefore, $\partial i_{\xi}(z_1)=\partial i_{\xi}(z_2)\in\partial Y_{\xi}^h$. Therefore, $(z_1,z_2)\in\Lam{\xi}(F)$.\\
{\bf Case 2:} $z_1$ is a parabolic limit point and $z_2$ is non-parabolic. Let $z_2=\Lambda(H_{b_0,\alpha}^h)$. So, $z\in\Lambda(\CC_{\alpha}^h)$. Identifying $z_1\in\partial\widehat{F}$, we have $\partial i_{\XH}(z_1)\in\partial\XH$. Moreover, $\partial i_{\XH}(z_1)\in\Lambda\BB_{\alpha}$ in $\partial\XH$. Since each $\BB_{\alpha}$ is an isometric section of $B$ in $\XH$, there exists $\xi\in\partial B$ such that $\partial i_{\XH}(z_1)\in\partial^{\xi}\XH$, i.e., $\partial i_{\XH}(z_1)=\Lambda\BB_{\alpha}\cap\partial^{\xi}\XH$. Then for a qi lift $\tilde{\beta}$ of a geodesic ray $[b_0,\xi)\subset B$ in $\XH$ such that $\tilde{\beta}\subset\CC_{\alpha}$, we have that $\tilde{\beta}(\infty)=\partial i_{\XH}(z_1)$. Note that $\tilde{\beta}\subset\YH_{\xi}$. For each $n>0$, let $b_n\in[b_0,\xi)$ such that $d_B(b_0,b_n)=n$. Then clearly, the sequence $\{\tilde{\beta}(b_n)\}$ converges to $\Lambda(\CC_{Y_{\xi},\alpha}^h)$ in $Y_{\xi}^h$, where $\CC_{Y_{\xi},\alpha}=\CC_{\alpha}\cap Y_{\xi}$. Now for a sequence $\{x_n\}$ in $F$ satisfying $x_n\to z_1$ in $F^h$, as before, let $\tilde{\alpha}_n$ and $\alpha^h_n$ denote a partially electrocuted geodesic in $\YH_{\xi}$ and a hyperbolic geodesic in $Y^h_{\xi}$, respectively, joining $x_n,\tilde{\beta}(b_n)$. Then as $n\to\infty$, $d_{\YH_{\xi}}(p,\tilde{\alpha}_n)\to\infty$ which implies that $d_{Y^h_{\xi}}(p,\alpha^h_n)\to\infty$. 
As $z_2$ is a parabolic limit point, $\partial i_{\xi}(z_2)=\Lambda(\CC_{Y_{\xi},\alpha}^h)$ in $Y_{\xi}^h$ and thus, $\partial i_{\xi}(z_1)=\partial i_{\xi}(z_2)$.\\
{\bf Case 3:} Both $z_1,z_2$ are non-parabolic and $z$ is a parabolic limit point of $\partial X^h$. Let $z=\Lambda\CC^h_{\alpha}$, and let $z_0=\Lambda H^h_{b_0,\alpha}$. Then $\partial i(z_0)=z$. By Case 2, there exists $\xi_1,\xi_2\in\partial B$ such that $(z_1,z_0)\in\Lam{\xi_1}(F)$ and $(z_0,z_2)\in\Lam{\xi_2}(F)$. \\
Let $\displaystyle{\big(\cup_{\xi\in\partial B}\Lam{\xi}(F)\big)^*}$ denote the transitive closure of of the relation $x\sim y$ on $\partial F^h$ if and only if $\partial i(x)=\partial i(y)$.
Then, $\Lam{X}(F)\subset\displaystyle{\big(\cup_{\xi\in\partial B}\Lam{\xi}(F)\big)^*}$. We also know that $\cup_{\xi\in\partial B}\Lam{\xi}(F)\subset\Lam{X}(F)$. So, we have the following result.
\begin{lemma}\label{trasitive closure}
	$\Lam{X}(F)=\displaystyle{\big(\cup_{\xi\in\partial B}\Lam{\xi}(F)\big)^*}$. 
\end{lemma} 
\subsubsection{Surjectivity of the CT maps}
 This is the relatively hyperbolic analogue of \cite[Theorem 6.26]{krishna-sardar}.
\begin{theorem}\label{thm:CT surjective}
	Suppose we have the hypothesis of Theorem \ref{main theorem} such that the fibers of the bundle are proper metric spaces. 
	Let $b\in A$. Suppose the CT map 
	$\partial i:\partial F_b^h\to\partial X^h$ is surjective. Then the CT map $\partial i_Y:\partial F_b^h\to\partial Y^h$ is also surjective.	
\end{theorem}
\begin{proof}Let $\partial i_{Y,X}:\partial Y^h\to\partial X^h$ denote the CT map of the pullback.
Let $z\in \partial Y^h$. Since $\partial i$ is surjective,
there exists $z_1\in\partial F_b^h$ such that $\partial i(z_1)=\partial i_{Y,X}(z)$. Suppose $\partial i_Y(z_1)\neq z$.
Then, $\partial i_Y(z_1)$ and $z$ are elements of $\partial Y^h$ identified under $\partial i_{Y,X}$. Then by \thmref{CT leaf in Y}(1), there exists $z_2\in\partial F_b^h\setminus\{z_1\}$ such that $\partial i_Y(z_2)=z$. Therefore, we are done. \end{proof}
%
%
%
%
%

\noindent
We have the obvious corollary.
\begin{corollary}
	Suppose $1\to (K,K_1)\to (G,N_G(K_1))\to (Q,Q_1)\to 1$ is a short exact sequence as in \thmref{combi for short exact seqn}. Suppose $Q'< Q$ is qi embedded and $G_1=\pi^{-1}(Q')$. 
	Then the CT map $\partial K^h\to\partial G_1^h$ is surjective.
\end{corollary}
\subsubsection{Quasiconvexity}
Let $$1\to K\to G\stackrel{\pi}{\to} Q\to 1$$ be a short exact sequence of (resp. relatively) hyperbolic groups. In \cite{mj-rafi}, the authors show that under certain conditions, a finitely generated infinite index subgroup of $K$ is a (resp. relatively) quasiconvex subgroup of $G$. We recall the results here.
\begin{theorem}\cite[Theorem 1.1]{mj-rafi}\label{alg qc mjrafi}
Let $1\to K\to G\stackrel{\pi}{\to} Q\to 1$ be a short exact sequence of hyperbolic groups, where $K$ is a (closed) surface group or a free group and $Q$ is convex cocompact in Teichmuller space or outer space respectively. Let $H<K$ be a  finitely generated infinite index subgroup of $K$. Then $H$ is a quasiconvex subgroup in $G$.
\end{theorem}
\begin{theorem}\cite[Theorem 1.3]{mj-rafi}\label{alg rqc mjrafi}
Let $K=\pi_1(S^h)$ be the fundamental group of a surface with finitely many punctures and let $K_1,\ldots,K_n$ be its peripheral subgroups. Let $1\to K\to G\stackrel{\pi}{\to} Q\to1$ and $1\to K_i\to N_G(K_i)\stackrel{\pi}{\to}Q\to 1$ be the induced short exact sequences of groups. Suppose $Q$ is a convex cocompact subgroup of the pure mapping class group of $S^h$. Let $K'<K$ be a finitely generated infinite index subgroup of $K$. Then $K'$ is a relatively quasiconvex subgroup of $G$. 	
\end{theorem}
\noindent
Generalizing \thmref{alg qc mjrafi}, Mj and Sardar proved the following. We keep the setup of Theorem \ref{alg qc mjrafi}.
\begin{theorem}\cite[Theorem 4.3]{sardar-mj}\label{alg qc mjsardar}
Let $Q_1<Q$ be a quasiconvex subgroup of $Q$ and $G_1=\pi^{-1}(Q_1)$. Let $H_1<G_1$ be an infinite index quasiconvex subgroup of $G_1$. Then $H_1$ is quasiconvex in $G$. 	
\end{theorem}
\noindent
We prove a relatively hyperbolic analogue of this result. We first recall the following from \cite{sardar-mj}.
\begin{prop}\cite[Proposition 2.8]{sardar-mj}\label{prop from ms}
Let $G$ be a finitely presented group, $H,K<G$, where $H$ is finitely generated. Let $\Gamma=\Gamma(G,S)$ and $\Gamma_1=\Gamma(H,S_1)$ be the Cayley graphs, where $S,S_1$ are finite generating sets of $G,H$ respectively. Let $D>0$ and $A\subset N_D(K)\cap H$ be an infinite subset. Suppose there exists $r\geq 1$ such that any pair of points in $A$ can be joined by a path which lies in an $r$-neighbourhood of $A$ in $\Gamma_1$. Then there exists a finitely generated infinite subgroup $K'\subset H\cap K$ such that $A/K'$ is finite. In particular, $A$ is contained in a finite neighbourhood of $K'$ in $\Gamma$.
\end{prop}
\begin{convention}
Let $G$ be a group hyperbolic relative to a collection of subgroups and let $H$ be a relatively hyperbolic subgroup such that the inclusion admits CT. We denote the CT lamination $\Lam{G^h}H^h$ simply by $\Lam{G}H$.
\end{convention}

\begin{lemma}\label{qc iff lam}\cite[Lemma 3.5]{mj-rafi}
Let $(G,\HH)$ be a relatively hyperbolic group and $(H,\HH')$ be a relatively hyperbolic subgroup such that $(H,\HH')\to(G,\HH)$ admits CT. Then $H$ is relatively quasiconvex in $G$ if and only if $\partial^{(2)}_G(H)=\emptyset$.	
\end{lemma}
\begin{lemma}\label{lams}\cite{mitraonss}
Let $(G,\HH)$ be a relatively hyperbolic group and $(H_1,\HH_1)$ and $(H_2,\HH_2)$ be relatively hyperbolic subgroups such that $H_2<H_1$ and $(H_2,\HH_2)\stackrel{i}{\to}(H_1,\HH_1)$ and $(H_1,\HH_1)\stackrel{j}{\to}(G,\HH)$ admit CT map. Then $\partial^{(2)}_G(H_2)=\partial^{(2)}_{H_1}(H_2)\cup\partial i^{-1}(\partial^{(2)}_{G}(H_1))$.
\end{lemma}
\begin{theorem}
Let $\KK=\{K_i\}_{i=1}^n,\HH=\{N_G(K_i)\}_{i=1}^n$. Let $$1\to (K,\KK)\to (G,\HH)\stackrel{\pi}{\to} (Q,Q)\to 1$$ be a short exact sequence of relatively hyperbolic groups. Let  $Q_1<Q$ be a quasiconvex subgroup of $Q$ and let $\pi^{-1}(Q_1)=G_1$. Suppose $H$ is an infinite index relatively quasiconvex subgroup of $G_1$. If $H$ is not relatively quasiconvex in $G$, then there exists an infinite index relatively quasiconvex subgroup of $K$ which is not relatively quasiconvex in $G$.
\end{theorem}
\begin{proof}
	Since $H$ is not relatively quasiconvex in $G$, by Lemma \ref{qc iff lam}, $\partial^{(2)}_G(H)\neq\emptyset$. Let $\gamma$ be a leaf of $\partial^{(2)}_G(H)$. Since $H$ is relatively quasiconvex in $G_1$, $\gamma$ is a leaf of $\partial^{(2)}_G(G_1)$. By Theorem \ref{CT leaf in Y}, there exists some $D'>0$ such that $\gamma$ lies in a $D'$-neighbourhood of $\Gamma_K^h$, where $\Gamma_K$ is the Cayley graph of $K$. Moreover, $\gamma(\pm\infty)$ are non-parabolic limit points of $\partial K^h$ as well as $\partial H^h$. This implies that there exists some $D>0$ such that $\gamma\cap H$ is an infinite subset lying in $N_D(K)\cap H$. By \propref{prop from ms}, there exists a finitely generated infinite subgroup $K'\subset H\cap K$ and $R>0$ such that $\gamma\cap H$ lies in an $R$-neighbourhood of $K'$ in $\Gamma$, where $\Gamma$ is a Cayley graph of $G$. This implies that $\gamma(\pm\infty)\in\Lambda_H(K')$.\\
	{\bf Claim 1:} $[K:K']=\infty$. Suppose not. Then $\Lambda_K(K')=\partial K$. Since $\Lambda_{G_1}(K)=\partial G_1$, we have that $\Lambda_{G_1}(K')=\partial G_1$. Now, note that $\Lambda_{G_1}(K')\subset\Lambda_{G_1}(H)$. This will imply that $\Lambda_{G_1}(H)=\partial G_1$, which is a contradiction. \\
	{\bf Claim 2:} $K'$ is not relatively quasiconvex in $G$. Suppose not. Then $\partial^{(2)}_G(K')=\emptyset$ and by \lemref{lams}, $\partial^{(2)}_H(K')=\emptyset$.  Now, by Lemma \ref{qc iff lam}, $\partial^{(2)}_G(H)\neq\emptyset$ and for any leaf $\gamma$ of $\partial^{(2)}_G(H)$, $(\gamma(-\infty),\gamma(\infty))\in\partial^{(2)}_G(K')$. This is a contradiction.	
\end{proof}
\noindent
Suppose we have the set up of Theorem \ref{alg rqc mjrafi}. As a corollary of the above result, we have the following.
\begin{theorem}\label{ct lam rqc}
	Let $Q_1<Q$ be a quasiconvex subgroup of $Q$ and let $\pi^{-1}(Q_1)=G_1$. Then any infinite index relatively quasiconvex subgroup of $G_1$ is a relatively quasiconvex subgroup of $G$. 	
\end{theorem}

\bibliography{pullback}
\bibliographystyle{amsalpha}

\end{document}